\documentclass[letterpaper, submit]{AAS}			% to submit to JAS
\usepackage[utf8]{inputenc}
\usepackage{multirow}
\usepackage{graphicx}
\usepackage{amsmath}
\usepackage{amsfonts}
\usepackage[version=4]{mhchem}
\usepackage{siunitx}
\usepackage{bm}
\usepackage{color}
\usepackage{longtable,tabularx}
\usepackage{algorithm}
\usepackage{marvosym}
\usepackage{wasysym}
\usepackage{subcaption,multicol,lipsum}
\usepackage[noend]{algpseudocode}
\usepackage{overcite}
\usepackage{tcolorbox}
\def\BState{\State\hskip-\ALG@thistlm}
\makeatletter
\def\BState{\State\hskip-\ALG@thistlm}
\makeatother
\title{How Many Impulses Redux}

\author{Ehsan Taheri \thanks{Research Assistant Professor, Department of Aerospace Engineering, Texas A\&M University, College Station, TX 77843-3141. AIAA Member.}, and John L. Junkins \thanks{Distinguished Professor, Department of Aerospace Engineering, Texas A\&M University, College Station, TX 77843-3141. AIAA Honorary Fellow.}}

\begin{document}

\maketitle

\begin{abstract}
A central problem in orbit transfer optimization is to determine the number, time, direction and magnitude of velocity impulses that minimize the total impulse. This problem was posed in 1967 by T. N. Edelbaum, and while notable advances have been made, a rigorous means to answer Edelbaum's question for multiple-revolution maneuvers has remained elusive for over five decades. We revisit Edelbaum's question by taking a bottom-up approach to generate a minimum-fuel switching surface. Sweeping through time profiles of the minimum-fuel switching function for increasing admissible thrust magnitude, and in the high-thrust limit, we find that the continuous thrust switching surface reveals the $N$-impulse solution. It is also shown that a \textit{fundamental} minimum-thrust solution plays a pivotal role in our process to determine the optimal minimum-fuel maneuver for all thrust levels. Remarkably, we find the answer to Edelbaum's question is not generally unique, but is frequently a set of equal-$\Delta v$ extremals. We further find, when Edelbaum's question is refined to seek the number of finite-duration thrust arcs for a specific rocket engine, that a unique extremal is usually found. Numerical results demonstrate the ideas and their utility for several interplanetary and Earth-bound optimal transfers that consist of up to eleven impulses or, for finite thrust, short thrust arcs. Another significant contribution of the paper can be viewed as a unification in astrodynamics where the connection between impulsive and continuous-thrust trajectories are demonstrated through the notion of \textit{optimal switching surfaces}.

\end{abstract}

\section*{Nomenclature}

% \noindent(Nomenclature entries should have the units identified)

{\renewcommand\arraystretch{1.0}
\noindent\begin{longtable*}{@{}l @{\quad=\quad} l@{}}
$\textbf{f}$ & vector function of unforced dynamics\\
$\mathbb{B}$ & control influence matrix \\
$c$ & exhaust velocity, m/s \\
$g_0$ & Earth gravitational acceleration, m/s$^2$\\
$\textbf{h}$ & specific angular momentum vector, km$^2$/s \\
$H$ & Hamiltonian\\
$I_{\text{sp}}$  & specific impulse, $s$ \\
$J$   & cost functional \\
$m$  & mass, kg\\
$N$ & number of impulses\\
$N_{\text{rev}}$ & number of revolutions \\
$N_{\text{bifr}}$ & number of bifurcations \\
$\textbf{p}$ & primer vector\\
$\textbf{r}$ & position vector, km \\
$S$  & switching function\\
$t$ & time, $s$\\
$T$ & thrust magnitude, Newtons\\
$T_{\text{min}}$ & critical minimum-thrust magnitude, Newtons\\
$\textbf{v}$  & velocity vector, km/s\\
$\textbf{x}$ & state vector \\%$[ p,f,g,h,k,l]^{\top}$ \\
$\Delta t$ & time interval, $s$\\
$\Delta v$ & impulse magnitude, km/s \\
$\rho$ & smoothing parameter \\
$\textbf{u}$ & control vector\\
$\bm{\alpha}$ & unit direction vector of thrust\\
$\delta$ & engine throttle input\\
$\bm{\psi}$ & vector of terminal constraints\\
$\bm{\lambda}$ & co-state vector associated with states\\
$\lambda_m$ & co-state associated with mass\\
$\bm{\eta}$ & vector of unknown design variables\\
$\mu$ & gravitational parameter, km$^3$/s$^2$\\
\multicolumn{2}{@{}l}{Subscripts}\\
0 & initial \\
f & final\\
$i$ & impulse \\
$\oplus$ & Earth \\
\Mars & Mars \\ 
T & target body\\
\Venus & Venus\\

\multicolumn{2}{@{}l}{Superscripts}\\
* & optimal
\end{longtable*}}
\section{Introduction} \label{sec:introduction}
Most space trajectory design algorithms make use of low-fidelity dynamical models and idealized control input assumptions to make the search space tractable [\citen{petropoulos2004shape, taheri2012shape,zuiani2012direct, avanzini2015solution,kluever2015designing,taheri2016initial,zeng2017shape,roa2016new,taheri2018shaping}]. Specifically, traditional impulsive-based trajectory analysis tools, typically with inverse-square gravity models, hold a special place for preliminary mission design [\citen
{longuski1991automated,izzo2007search,olympio2009designing,vasile20108,englander2012automated, abdelkhalik2012dynamic,chilan2013automated,ellison2018analytic,landau2018efficient}]. The output of preliminary mission design studies is the starting design for higher-fidelity optimization. Making these approximations in preliminary mission design is driven by practicality, which is natural, given the level of complexity of the overall mission design challenge. Edelbaum's question [\citen{edelbaum1967many}] was posed in the setting of inverse-square gravity models, however, the generalized question for high-fidelity force models is straightforward to ask (but not to answer).

Impulsive solutions are important since they determine both the \textit{theoretical} minimum-time and minimum-fuel extremals and also provide reachability insights. For the most part, preliminary mission design methods rely on low-fidelity dynamical models, which in turn, frequently leads to  analytical propagation of the state dynamics through Keplerian orbit models [\citen{sims1997preliminary}] or by utilizing the solution of Lambert's problem [\citen{gooding1990procedure,battin1999introduction,izzo2015revisiting,arora2015partial,bombardelli2018approximate}]. Impulsive maneuvers are also used extensively for formation flight optimal control problems [\citen{vaddi2005formation,gurfil2005relative,woffinden2007navigating,ichikawa2008optimal,kim2009hybrid,huang2012optimal,spencer2015maneuver,roscoe2015formation,sobiesiak2015impulsive, starek2016fast, chernick2018new}] and orbit reachability analyses problems [\citen{vinh1995reachable,xue2010reachable,zhang2013reachable,wen2014orbital,holzinger2014orbit,wen2014precise,chen2018new}]. 

For the impulsive thrust idealization, a fundamental quest has been to determine the optimal number, times, magnitude, and direction of the impulses, to accomplish general three-dimensional (3D) multiple-revolution orbit transfers while minimizing the total $\Delta v$. This is Edelbaum's still not rigorously answered question [\citen{edelbaum1967many}]: How many impulses? 

Optimal continuous and impulsive formulations were originally investigated by Lawden beginning in the early 1950's. In his seminal and pioneering 1963 work on optimal trajectories, Lawden derived a set of criteria [\citen{lawden1963optimal}] that define the optimality of impulsive solutions by introducing the ``primer vector'', $\textbf{p}$. The primer vector, for either continuous or impulsive thrusting, defines the instantaneous optimal direction for the thrust vector. Due to the importance of Lawden's impulsive necessary conditions, they are repeated here: 1) the primer vector and its first derivative are continuous everywhere, 2) the magnitude of the primer vector remains less than unity, i.e., $p \equiv \| \textbf{p} \| < 1$ except for the impulse times where $\| \textbf{p} \| = 1$, 3) at the impulse times, the primer vector is a unit vector along the optimal direction of impulse, and 4) at any intermediate impulse time, $dp/dt = \dot{p} = \dot{\textbf{p}}^{\top} \textbf{p} = 0$. 
Undoubtedly, Lawden's introduction of the concept of primer vector is the most fundamental breakthrough in the field of space trajectory optimization. 

Violation of the necessary conditions can be used as a measure of sub-optimality of approximate impulsive solutions, and has been used to improve sub-optimal impulsive solutions [\citen{handelsman1968primer}]. Specifically, first-order variation of a 2-impulse cost functional is derived to establish necessary conditions for a small variation that result in an improved solution through: 1) the introduction of an additional mid-course impulse, and 2) introduction of terminal coasts (either initial or final). We briefly discuss the $N$-impulse literature wherein a number of algorithms have been devised to seek minimum-$\Delta v$ impulsive solutions. The above ideas can be utilized to improve approximate impulsive solutions through classical gradient-based optimization algorithms [\citen{jezewski1968efficient}]. The time history of the primer vector (and possible violation of the optimality conditions) is frequently used in numerical algorithms to place additional impulses near maxima of $p$. The time and location of the impulses have to be finalized by direct optimization. So, we have a four-dimensional augmentation of the search space for every additional impulse. In order to minimize the cost function, the point at which $p$ takes its maximum value is usually taken as an initial iterate [\citen{handelsman1968primer}]. This approach leads to a direct method, a multivariate search problem, that has to be solved in a robust manner to result in a converged solution [\citen{prussing20102}]. 

The search for $N$-impulse solutions has been most commonly initiated from a minimum-$\Delta v$ 2-impulse (Lambert) solution for which the existence of $2 N_{\text{rev,max}}+1$ solutions is demonstrated in [\citen{prussing1993orbital,prussing2000aclass}], where $N_{\text{rev,max}}$ is the maximum number of revolutions that has to be determined and depends on the prescribed time of flight. Multiple-revolution Lambert algorithms are used to generate multiple reference trajectories, each of which are considered for multiple impulse optimization and further improvements. An initial policy is required to select those solutions that are hypothesized to lead to improvements, which in most cases, is to select among the non-unique Lambert solutions those with cheaper 2-impulse $\Delta v$ requirements. The improved solution (a 3-impulse solution) divides the problem into two new sub-arcs, each one of which can be treated similar to the original 2-impulse solution. However, this approach requires that a decision be made on which sub-arc to be optimized first. Therefore, there are $N-1$ decisions to be made, which result in $2^{N-1}$ different possibilities (analogous to branches of a tree-search problem) for values of the cost functional. Alternatives for decision making are given in [\citen{jezewski1968efficient}]. While presenting important advancements, these heuristic bootstrapping approaches with associated gradient-based solvers may get stuck in local optima since there is no guarantee of a unimodal performance surface and these methods rely on classical parameter optimization methods [\citen{prussing2000aclass,prussing1986optimal,hughes2003comparison}]. A common aspect among most of these methods is that they rely solely on impulsive-based solutions. Application of semi-infinite convex optimization using a relaxation scheme and duality theory in normed linear spaces is demonstrated in [\citen{arzelier2016linearized}] for fixed-time minimum-fuel rendezvous between close elliptic orbits without fixing a priori the number of impulses.

In order to avoid local sub-optimal convergence, a number of studies have focused on the application of evolutionary algorithms [\citen{abdelkhalik2007n, pontani2012particle,luo2007optimization,sentinella2009cooperative,luo2010interactive}] that compromise between local and global search processes to identify multiple local minima. In addition, indirect-based methods are studied in [\citen{colasurdo1994indirect}] and a homotopic-based indirect scheme is presented in [\citen{shen2015global}], which improves potential 2-impulse Lambert solutions out of the total $2N_{\textrm{rev,max}}+1$ solutions. Nevertheless, while all of the aforementioned methods have been able to find multi-impulse solutions that improve on the 2-impulse solutions with varying degrees of success, none can claim global optimality, nor can they answer Edelbaum's question with certainty. 

Under some conditions, i.e., a linear neighborhood of reference trajectories, the maximum number of impulses is shown not to exceed the number of state variables [\citen{neustadt1964optimization, edelbaum1967many}]. For a linear system, Lawden's necessary conditions are also sufficient for an optimal trajectory [\citen{prussing1995optimal}]. For both circular and elliptic orbits, the necessary and sufficient conditions for the optimal (fixed-terminal state and fixed-time) solution are derived [\citen{carter1991optimal,carter1992fuel}]. For rendezvous and transfer problems assuming a linear dynamical model, the number of impulses is at most equal to the dimension of the state space [\citen{prussing1995optimal}]. The case of planar transfer between co-planar elliptical orbits is also studied in [\citen{lawden1992optimal}]. 
Rendezvous of two spacecraft in neighboring near-circular non-coplanar orbits is reported with up to six impulses [\citen{baranov2012six}]. It is shown that the representation of the primer vector in polar coordinates leads to the separation of the in-plane and out-of-plane components of the primer vector. A complete analytic solution for the out-of-plane component of the primer vector is shown to exist, which is independent of the semi-major axis of the transfer orbit [\citen{iorfida2016geometric}]. The problem of time-fixed fuel-optimal out-of-plane elliptic rendezvous between spacecraft in a linear setting is studied with a complete analytical closed-form solution [\citen{serra2018analytical}].

On the other hand, there is a direct theoretical connection between optimal finite-thrust continuous control and an optimal sequence of velocity impulses; this connection becomes apparent in the \textit{switching surfaces} introduced and discussed herein. The fact that impulsive maneuvers constitute the limiting case of the more general finite-thrust minimum-thrust trajectory optimization problems has been stated in [\citen{neustadt1965general,colasurdo1994indirect}] not only when the thrust magnitude increases, but also when the transfer time increases [\citen{gergaud2007orbital}]. In fact, the work of Zhu, et al [\citen{zhu2017solving}], has been motivated by this fact that ``... the optimal bang-bang control and impulsive maneuvers can be obtained through continuously increasing the thrust magnitude from an optimal minimum-thrust solution.'' Implicitly, they used neighboring converged co-states to initiate the TPBVP solution for each new assigned thrust magnitude, $T$ along with appropriate changes that result in the switching function, $S$ of the indirect formulation of optimal control trajectories. 

The procedure Zhu, et al followed in [\citen{zhu2017solving}] is dependent upon beginning with a 2-impulse Lambert solution. As we discuss below, assuming a 2-impulse starting solution is generally not theoretically justified, nor does it offer any convergence guarantee, because the  optimal $N$-impulse solution we seek is obviously not generally near a Lambert 2-impulse trajectory. For sufficiently short time of flight, however, we can anticipate the 2-impulse solution will indeed be the optimal impulsive extremal and will be unique for orbit transfer maneuvers spanning a fraction of one revolution. However, for longer time intervals, as will be evident in the developments herein, one must consider fully the local behavior associated with local extrema on each feasible specification of the number of en-route revolutions along the transfer orbit.

In this investigation, we use indirect optimization methods. The most fundamental feature of indirect methods is that any trajectory satisfying the necessary conditions and all boundary conditions (BCs) is guaranteed to yield a local extremal. In space applications the equations of motion are of relatively low order, so indirect methods, when combined with reliable initialization and homotopy approaches, are attractive and lead to fast convergence to at least local extremals. These approaches are especially attractive when no state variable inequality path constraints are imposed [\citen{taheri2016enhanced}]. Indirect methods utilizing convergence enhancement homotopy techniques and control smoothing methods, while artistic, have ameliorated many of the challenges of numerically solving the TPBVPs and have been applied successfully to a number of optimal control problems [\citen{kluever1997optimal,caillau2001coplanar,bertrand2002new,caillau20033d,haberkorn2004low,bonnard2005geometric,la2006indirect,russell2007primer,silva2010smooth,jiang2012practical,dutta2012peer,taheri2016enhanced,pan2016double,zhao2016minimum,zhao2016initial,chi2017homotopy,zhu2017solving,taheri2017co,mall2017epsilon,sullo2017low,zhu2017geometric,saranathan2017relaxed,saghamanesh2018robust,perez2018fuel,junkins2018exploration,pan2018new,taheri2018generic,pan2019quadratic}].
Indirect methods exploit first-order necessary conditions, which in principle and most frequently, converge to local extremals. Therefore, methods for establishing starting co-states within the domain of attraction of the desired global extremum are desirable. As is evident herein, we have established important insights on this difficult issue. Local extrema have been found, somewhat analogous to the multi-revolution Lambert problem, to be associated with the number ($N_{\text{rev}}$) of intermediate revolutions the transfer extremal makes en-route to satisfying the final BCs. So analogous to the Lambert's problem, we specify $N_{\text{rev}}$ and seek out all local extremals associated with each $N_{\text{rev}}$. The global extremum and the associated optimal $N_{\text{rev}}$ is obtained by selecting the minimum of the minima. 
The ``fundamental'' minimum-thrust solution, i.e., the minimum of all local minima is critical in the analysis of a comprehensive approach that is devised to address global minimum-fuel solutions and its related optimal $N$-impulse solutions. Remarkably, this fundamental minimum-thrust solution belongs to the same continuous extremal field map of neighboring minimum fuel extrema that are the solutions we seek. This method is different from the previously mentioned approaches that are reviewed above. In general, the fundamental minimum-thrust solutions are obtained to establish the first profile on a \textit{switching surface} for all thrusts $T > T_{\text{min}}$. Since the minimum thrust, $T_{\text{min}}$ to reach the final state is established by this process, obviously the minimum thrust is the boundary of reachability domain. For each $N_{\text{rev}}$, the associated switching surface (for $T > T_{\text{min}}$) is an ensemble of all switching functions associated with the entire family of extremal solutions. These switching surfaces turn out to be very informative tools and provide an enhanced global understanding of the space of extremal solutions while revealing, for example, \textit{late-departure} and \textit{early-arrival} boundaries. The fundamental minimum-thrust solution is critical in constructing the switching surface that for the case of limiting thrust magnitude, $T \rightarrow \infty$ reveals the associated impulsive solution. 

We consider herein a number of test cases along with their switching surfaces and their associated impulsive solutions. An algorithm is outlined, which uses the high-thrust behavior to approximate accurately the $N$-impulse solutions for multiple-revolution trajectories. A final process is described that adjusts these impulses slightly to isolate the optimal impulsive solution. It is possible to extend the proposed process to account for the high-fidelity force models and to converge to the final solution.

The paper is organized as follows. First, a concise review of the formulation of the minimum-fuel TPBVPs using the indirect method and Pontryagin's Minimum Principle (PMP) is given. Then, we review the details of the continuation procedure when the magnitude of thrust is swept to generate the desired switching surfaces. The details of a robust algorithm that has been used to generate the $N$-impulse solutions are presented next. Then, the results are given for a number of test cases. Interpretation of the results are presented that provide insights to neighboring extremals for various thrust values, $T$ and the limiting high-thrust extremal, which are the corresponding impulsive maneuvers. In reviewing the results for these distinctly different family of optimal transfers, the versatility of the methodology to handle various unique circumstances becomes evident. Then, a discussion is given on the number of impulses, and interestingly, the non-uniqueness of the impulsive solutions and the computational effort of the method. Application of the method for solving transfer problems from the GTO to a Halo orbit around the L1 point of the Earth-Moon restricted three-body model is also presented. Finally, concluding remarks are presented.

\section{Indirect Formulation of Minimum-Fuel Problem and Continuation Procedure} \label{sec:OCP}
In this section, equations of motion and the minimum-fuel cost functionals are discussed. Optimal control theory is used to establish the necessary conditions for optimal trajectories and to define the TPBVPs to be solved numerically. In all example problems studied, it is assumed that there are no state variable path constraints other than the terminal BCs. Numerical schemes used for solving the resulting TPBVPs are discussed separately.

\subsection{Equations of Motion}\label{sec:eom}
We consider the trajectory optimization for a spacecraft moving in an inverse-square gravitational field of a central body where the spacecraft is affected also by the acceleration induced by an on-board propulsion system. Implicit in this problem, the spacecraft attitude must be controlled to steer the thrust vector, however, we ignore rotational dynamics. This usual approximation is justified, because the controlled attitude error dynamics is typically several orders of magnitude faster than the orbital dynamics and attitude errors are usually a small fraction of a degree. This approximation is well-justified to design virtually all optimal interplanetary trajectories, and most near-Earth trajectories. 

The equations of motion are expressed in terms of the modified equinoctial orbital elements (MEEs) [\citen{walker1986set}] and the variation of mass is included. MEEs are suitable for optimization of low-thrust trajectories because the most general MEE representation includes circular, elliptic and hyperbolic orbits without singularities at zero eccentricity or zero inclinations. Unlike the inertial Cartesian coordinates that are changing quickly over a revolution, the MEEs are well-behaved and varying slowly except for the true perturbed true longitude, $l$, which is a smoothly varying  function of time that reduces to Kepler's equation in the absence of perturbations [\citen{taheri2016enhanced,junkins2018exploration}]. Furthermore, prescribing the osculating final orbit's true longitude as a terminal BC permits convenient control of $N_{\text{rev}}$ in accounting for the total angular displacement through intermediate revolutions and fractions thereof by simply adding $2\pi N_{\text{rev}}$ to the final true longitude $l_f$ (in the osculating final orbit). Specifically, we replace $l_f$ by $l^*_f = l_f + 2\pi N_{\text{rev}}$, and choose integer values for $N_{\text{rev}}$ over a feasible set. This takes advantage of the evident fact that having the osculating true longitude as a coordinate permits $N_{\text{rev}}$ to be specified in the final BC. This is a key element in formulating and solving TPVBPs where multiple extremals with different number of revolutions are possible to occur [\citen{taheri2016enhanced}].
%\vspace{2mm}

\begin{tcolorbox}[colback=green!5,colframe=green!40!black,title=Why is it important to have control over $N_{\text{rev}}$?]
By removing the freedom to converge to any $N_{\text{rev}}$, we ensure that the resulting switching surfaces are unique for each specification of $N_{\text{rev}}$. This is a critical enabler since our goal is to perform a systematic study that avoids getting solutions with different number of en-route revolutions.
\end{tcolorbox}

The procedure we develop herein finds the optimal $N_{\text{rev}}$ to minimize a standard minimum-fuel performance index ($ J = \frac{T}{c} \int_{t_0}^{t_f} \delta~dt $) over all feasible $N_{\text{rev}}$ specifications where $\delta$ is engine throttling input. Let $\textbf{x} = [p,f,g,h,k,l]^{\top}$ and $m$ denote the vector of MEEs, and the spacecraft mass, respectively, and let $\textbf{u} = [u_r, u_t, u_n]^{\top} = \frac{T}{m} \delta \bm{\alpha} $ denote the thrust acceleration vector with its components expressed in the local-vertical/local-horizontal (LVLH) osculating orbital reference frame. $\bm{\alpha}$, $T$ and $\delta \in [0,1]$ denote the thrust steering unit vector, thrust magnitude and engine throttle input, respectively. The state/co-state dynamics become  
\begin{align} 
\dot{\textbf{x}} =&~\textbf{f} + \frac{T}{m} \mathbb{B} \bm{\alpha} \delta , \label{eq:elementsdiffeq}\\ 
\dot m =& -\frac{T}{c}\delta, \label{eq:massdiffeq} \\
\dot{\bm{\lambda}} = &- \left [ \frac{\partial H}{\partial \textbf{x}}\right ]^{\top}, \label{eq:costatedynamics}\\
\dot{\lambda}_m =& - \frac{\partial H}{\partial m} \label{eq:masscostate},
\end{align}
where $c=I_{\text{sp}}g_0$ is the exhaust velocity, $I_{\text{sp}}$, $g_0$ are the engine's specific impulse and the gravitational acceleration at sea level, respectively. The $\textbf{f} = \textbf{f}(\textbf{x})\in \mathbb{R}^6$ is the unforced part of the state dynamics and $\mathbb{B} = \mathbb{B}(\textbf{x}) \in \mathbb{R}^{6 \times 3}$ denotes the control influence matrix 
\begin{align} \label{eq:AB}
\textbf{f} =\begin{bmatrix}
 0 \\
 0 \\
 0 \\
 0 \\
 0 \\
 \sqrt{\mu p}(\frac{w}{p})^2
 \end{bmatrix},
\mathbb{B} = \begin{bmatrix}
0                          & \frac{2p}{w}\sqrt{\frac{p}{\mu}}               & 0\\
\sqrt{\frac{p}{\mu}}\sin(l) & \sqrt{\frac{p}{\mu}}\frac{1}{w}[(w+1)\cos(l)+f] & -\sqrt{\frac{p}{\mu}}\frac{g}{w}[h\sin(l)-k\cos(l)]\\
-\sqrt{\frac{p}{\mu}}\cos(l)& \sqrt{\frac{p}{\mu}}\frac{1}{w}[(w+1)\sin(l)+g]            & \sqrt{\frac{p}{\mu}}\frac{f}{w}[h\sin(l)-k\cos(l)]\\
0                          & 0                                              & \sqrt{\frac{p}{\mu}}\frac{s^2\cos(l)}{2w}\\
0                          & 0                                              & \sqrt{\frac{p}{\mu}}\frac{s^2\sin(l)}{2w}\\
0                          & 0                                              & \sqrt{\frac{p}{\mu}}\frac{1}{w}[h\sin(l)-k\cos(l)]\\
\end{bmatrix}.
\end{align}
In these equations, two intermediate positive variables are $w=1+f\cos(l)+g\sin(l)$, $s^2=1+h^2+k^2$, and $\mu$ is the gravitational mass parameter of the central body. The co-state vector associated with the equinoctial orbit state vector is denoted by $\bm{\lambda} = [\lambda_p,\lambda_f,\lambda_g,\lambda_h,\lambda_k,\lambda_l]^{\top}$ and $\lambda_m$ is the co-state associated with the mass. $H$ is the Hamiltonian corresponds to the minimum-fuel cost functional, $J = \frac{T}{c}\int_{t_0}^{t_f} ~\delta~dt$, where $t_0$ and $t_f$ are fixed. The Hamiltonian becomes
\begin{equation}\label{eq:Hamiltonian}
 H = \frac{T}{c} \delta +\bm{\lambda}^{\top} \left [\textbf{f}+\frac{T}{m} \mathbb{B} \bm{\alpha} \delta \right ] - \lambda_m \frac{T}{c}\delta.
\end{equation}
The optimal control direction and throttling input that minimize $H$ are 
\begin{align}
\bm{\alpha}^* = & \frac{\textbf{p}}{||\textbf{p}||},  \label{eq:optimaldirection} \\
\delta^{*}(S) = & \left \{
  \begin{tabular}{cc}
  1 & \text{if $S>0$}  \\
  0 & \text{if $S<0$}
  \end{tabular}
\right \} = \frac{1}{2} \left [ 1 + \text{sign}(S) \right], \label{eq:delta}
\end{align}
where the \textit{primer vector}, $\textbf{p} \equiv -\mathbb{B}^{\top} \pmb{\lambda}$ and the switching function, $S$ is defined as
\begin{equation}\label{eq:sfm-p1}
S=\frac{c ||\mathbb{B}^{\top} \pmb{\lambda}||}{m}+\lambda_m - 1.
\end{equation}
In addition, the hyperbolic tangent smoothing (HTS) method [\citen{taheri2018generic,taheri2018genericconf}] is used as a means for smoothing the otherwise jump-discontinuous engine throttle input as
\begin{equation} \label{eq:smoothedcontrolp1}
\delta^*(S) = \frac{1}{2} \left [ 1 + \text{sign}(S) \right] \cong  \delta^*(S,\rho) = \frac{1}{2} \left [ 1 + \tanh \left( \frac{S}{\rho} \right ) \right],
\end{equation}
where $\rho > 0$ is the smoothing level (and is used as the continuation parameter for the numerical continuation procedure). The HTS method accurately approximates the optimal throttle, $\delta^*(S)$ by a smooth, differentiable function, $\delta^*(S,\rho)$; this smoothed throttle is quite effective since it enlarges the domain of convergence such that, in most cases, even a moderate number of random sets of initial co-state guesses leads to convergence to the local extrema. For the rare cases that $S = 0$ for a finite time interval, we may have a singular control ($0 < \delta < 1$). In the event a singular sub-arc is encountered, $S = \dot{S} = \ddot{S} = 0$ can be investigated as functions of ($t,\textbf{x}(t),\bm{\lambda}(t)$), along with the remaining necessary conditions (including Kelly condition [\citen{kelley1964second})], to see if there exists a singular throttle function $0 < \delta^* <1$.

For a fixed-time rendezvous problem, the final BCs (seven equality constraints) can be written as
\begin{equation}\label{eq:posvelconp1}
\bm{\psi}(\textbf{x}(t_f),t_f) = \left [ [\textbf{x}(t_f) - \textbf{x}_T]^{\top},\lambda_m(t_f) \right ]^{\top} = \textbf{0},
\end{equation}
where $\textbf{x}_T$ denotes the final target state values associated with the target body. The final value of the mass co-state has to be zero since the final mass is free. For multi-revolution transfers, the number of revolutions, $N_{\text{rev}}$ is an unknown and has to be determined. Computational experience indicates, for $N_{\text{rev}}$ greater than some problem-dependent minimum integer, there is one local extremum for each $N_{\text{rev}}$ choice. 

While we have no theoretical proof that there is only one extremal for each $N_{\text{rev}}$ choice (for \textit{fixed terminal BCs and time of flight}), we believe this to be true for inverse-square force fields based on extensive computations. This is an important point, because our method presently hypothesizes this to be true. If there is only one minimum-fuel local extrema per $N_{\text{rev}}$, then the minimum of all local minima will identify $N^*_{\text{rev}}$ and the global minimum-fuel trajectory.

Let $\textbf{z} = [\textbf{x}^{\top},m,\bm{\lambda}^{\top},\lambda_m]^{\top}$ denote the state-costate vector, then, we can write,
\begin{equation} \label{eq:F}
\dot{\textbf{z}} = \textbf{F} = [\dot{\textbf{x}}^{\top}, \dot{m}, \dot{\bm{\lambda}}^{\top}, \dot{\lambda}_m]^{\top},
\end{equation}
where $\bm{\alpha} = \bm{\alpha}^*$ and $\delta = \delta^*(S,\rho)$ (Note $\dot{\textbf{x}}$, $\dot{m}$, $\dot{\bm{\lambda}}$, $\dot{\lambda}_m$ is shorthand for the RHS of Eqs.~\eqref{eq:elementsdiffeq},\eqref{eq:massdiffeq},\eqref{eq:costatedynamics},\eqref{eq:masscostate}). Once, these values are substituted into $\textbf{F}$, the equations of motion can be integrated numerically, if initial BCs are fully specified. However, only the initial state $\textbf{x}(t_0) = \textbf{x}_0$ and $m(t_0) = m_0$ are specified. The final state $\textbf{x}(t_f)$ as well as the final co-states are a function of the initial co-state $\bm{\eta}(t_0)$ where $\bm{\eta}(t_0) = [\bm{\lambda}^{\top}(t_0),\lambda_m(t_0)]^{\top}$ is the vector of unknowns to be determined such that Eq.~\eqref{eq:posvelconp1} is satisfied. Thus, we have a TPBVP that requires a starting estimate $\bm{\eta}(t_0)$ within the domain of convergence of the algorithm used to satisfy the prescribed BCs. There are seven constraints in Eq.~\eqref{eq:posvelconp1} and seven unknown elements in $\bm{\eta}(t_0)$. 
\begin{figure}[htbp!]
\begin{multicols}{2}
\centering
\includegraphics[width=0.52\textwidth]{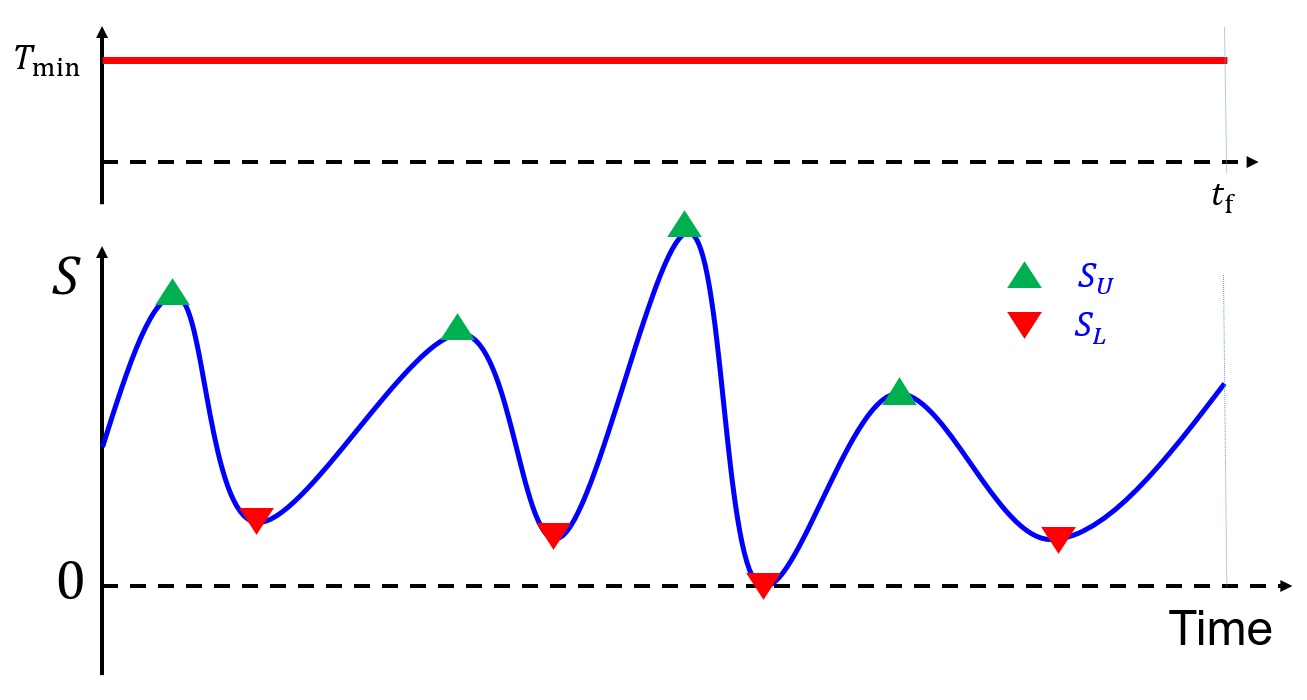}
\caption{Typical time histories of thrust magnitude and switching function for $T_{\text{min}}$.}
\label{fig:Tminshow}
\hfill
\centering
\includegraphics[width=0.47\textwidth]{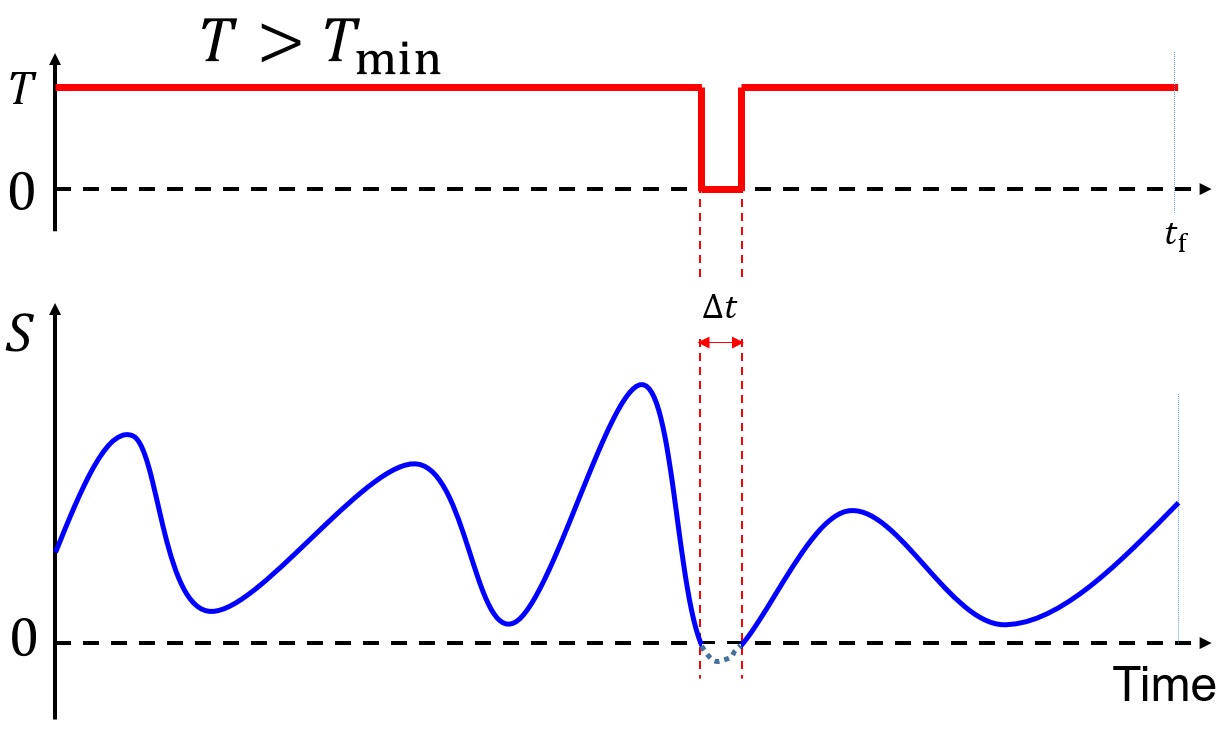}
\caption{Time histories of thrust magnitude and switching function for $T$ slightly greater than $T_{\text{min}}$.}
\label{fig:TbigerTmin}
\end{multicols}
\end{figure}

\section{Thrust Magnitude Sweeping for Construction of Switching Surfaces} \label{sec:continuationprocedure}
The continuation procedure over the thrust magnitude is considered to generate a family of switching function profiles that constitute the switching surface. It is noteworthy that the application of switching surfaces (extremal field maps) for analysis of globally optimal coplanar time-free orbit transfers has been investigated in [\citen{breakwell1964minimum, small1971minimum, small1982globally}]. 

Figure \ref{fig:Tminshow} depicts the time histories of a typical switching function and its associated (constant) thrust profile for $T_{\text{min}}$. Local extrema of the switching function are denoted by triangles. A slight increase in the thrust magnitude leads to a downward shift and distortion of the switch function and the appearance of a coast arc for a finite time interval, $\Delta t$ as is shown in Figure \ref{fig:TbigerTmin}. Ultimately, a procedure can be devised to sweep over increasing values of the thrust magnitude until the zeros of the switching function occur in pairs that occur a small $\Delta t < \epsilon$ apart. If the thrust duration satisfies $\Delta t = \epsilon < (t_f-t_0)/1000$, one can usually approximate the short thrust arcs as impulses.

% \begin{figure}[htbp!]
% \begin{multicols}{2}
% \centering
% \includegraphics[width=0.4\textwidth]{figures/HighThurstSwitchThrust.jpg}
% \caption{Time histories of switching function and thrust for $T \gg T_{\text{min}}$.}
% \label{fig:HighT}
% \hfill
% \centering
% \includegraphics[width=0.6\textwidth]{figures/EDSS_2.eps}
% \caption{Depiction of a representative switching surface that consists of six thrust ridges and seven coast canyons; thrust ridges narrow to approach six impulses at high thrust, $T_{\text{max}}$.}
% \label{fig:SS}
% \end{multicols}
% \end{figure}

\begin{figure}[htbp!]
\centering
\includegraphics[width=3.0in]{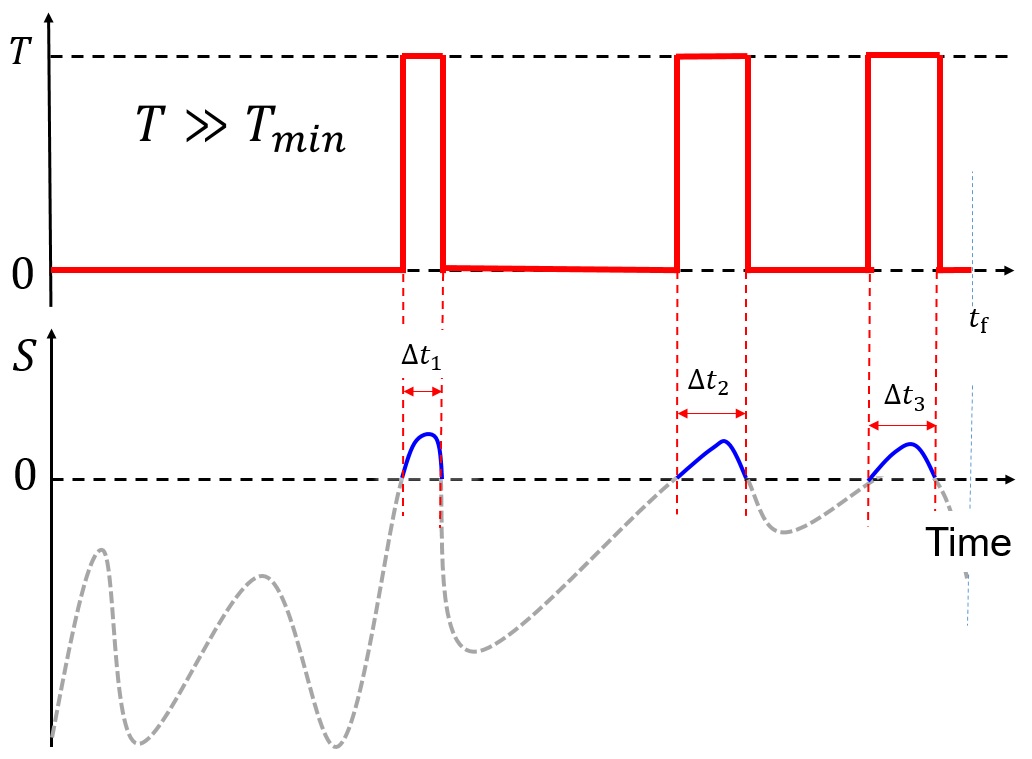}
\caption{Time histories of thrust magnitude and switching function for $T \gg T_{\text{min}}$.}
\label{fig:HighT}
\end{figure}

We find, for a sufficiently large thrust that the time duration of all thrust arcs becomes shorter than $\epsilon$ (a prescribed threshold based on the mission time). Then, we can approximate the thrust as impulsive with near-negligible error. Figure \ref{fig:HighT} shows representative profiles of the switching function and thrust magnitude versus time for very high thrust values, i.e., $T \gg T_{\text{min}}$. As the thrust magnitude is assigned increasingly large values, we observe that the time duration of all the thrust arcs shrink while the thrusting time sequence remains approximately unchanged and the solution becomes closely approximated by isolated impulses. At this stage, we seek to replace the continuous thrust by a finite number of impulsive thrusts by formulating and solving an $N-$impulse trajectory optimization problem, with the number, times, magnitudes and direction of the thrust impulses known approximately. 

\begin{figure}[htbp!]
\centering
\includegraphics[width=5.0in]{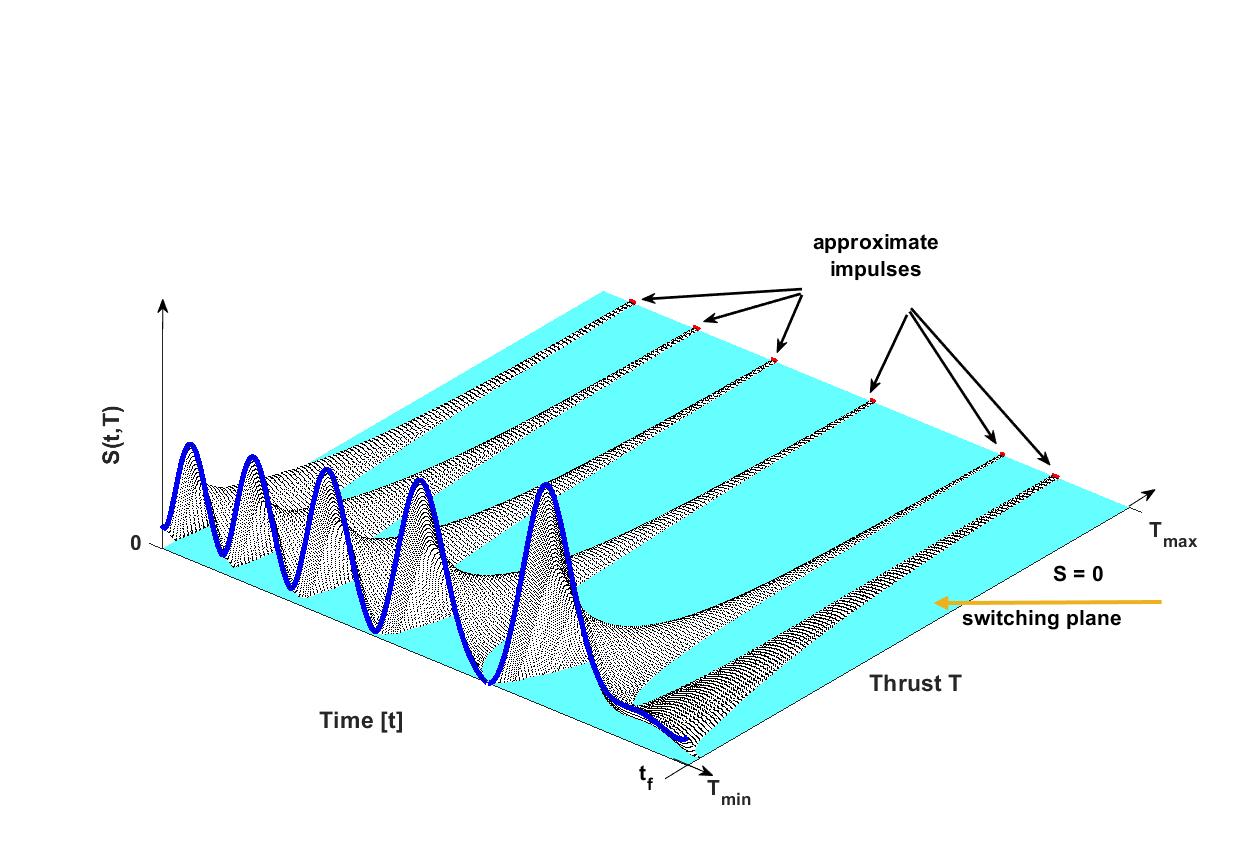}
\caption{A representative switching surface that consists of six thrust ridges and seven coast canyons; thrust ridges narrow to approach six impulses at high thrust, $T_{\text{max}}$.}
\label{fig:SS}
\end{figure}

Our goal is to generate a surface that is formed by sweeping the variable of interest, in this case thrust magnitude, $T$ and concatenating the switching functions. Figure \ref{fig:SS} depicts a representative switching surface where a solid blue curve denotes the switching function associated with the minimum-thrust extremal. In practice, we may require logarithmic scales to reveal sufficient details of these surfaces. Using a topographic analogy, increasing $T$ leads to the $S < 0$ coast ``canyons'' being wider while the $S > 0$ thrust ``ridges'' have lower peak $S$-values and become more narrow (in time). This switching surface has six thrust ridges and seven coast canyons when thrust magnitude is swept in its defined bound $T \in [T_{\text{min}}, T_{\text{max}}]$ where $T_{\text{max}} \gg T_{\text{min}}$. With increasing thrust, the time duration of all thrust arcs become smaller. Qualitatively, it is useful to consider the plane defined by $S = 0$ to represent the surface of an $S = 0$ ``lake'' defined by its shore lines (contour) intersection defining the boundary with the $S > 0$ topography. The thrust ridges at high thrust magnitudes approach six impulses of negligible time duration as the ``coast lakes'' become wider. The considered switching surface in Figure \ref{fig:SS} has a well-behaved topography, however, for many orbit transfers, the width of thrust ridges may not always decrease monotonically as the thrust magnitude is increased as in this illustration, and in some cases there are surprising and counter-intuitive features. Specifically, there are cases in which a thrust ridge will be created at some critical thrust value (either as an independent ``island'' in the middle of one of the coast ``lakes'' or as a ``peninsula'' that ``breaks off'' from a thrust ridge). These cases are associated with bifurcations that occasionally occur in the switching surface. 
%Figure \ref{fig:EMSS} and \ref{fig:ETo1989MLSS} below show the $S = 0$ contour boundaries of actual thrust ridges (dark blue) and $S < 0$ ``lakes'' in the coast canyons (light blue), for two specific families of optimal orbit transfers.

The individual switching surface topography for each orbit transfer is a function of the two sets of orbital BCs (including especially, relative phasing, inclination, size, shape and orientation) and the force model assumptions. The high dimensionality of the space that underlies each switching surface makes it difficult to predict the fine structure meandering of these surfaces, especially at low thrust levels. The optimal control switching surface is a fundamental attribute of controls associated with the family of extremals and does not depend, for example, on the choice of coordinates, although nonlinearity and efficient convergence to the solution underlying TPBVP does indeed depend on coordinate choice [\citen{taheri2016enhanced,junkins2018exploration}]. A discussion on this point is given later in a separate section. 

Study of the topography of the particular generated switching surface associated with a family of extremals, as will be shown, is very useful for trajectory and mission design purposes. Note that every time slice (constant $T$ profile) of the surface of Figure \ref{fig:SS} corresponds to the on-off switching function for a particular extremal trajectory, i.e., a minimum-fuel optimal transfer between the prescribed initial and final states over $t \in [t_0,t_f]$.

\section{Optimization Scheme For N-Impulse Solutions} \label{sec:Nimpulseop}
Given our ability to use the switching surface large thrust limiting behavior to approximate accurately the time, direction and magnitude of all velocity impulses, only slight adjustments are required to achieve final convergence. Our experience is that for multiple-revolutions, the convergence of the final direct optimization problem can be made more efficient if the whole trajectory is divided into several segments with the boundary of each segment defined by each impulse time and a forward-backward numerical integration scheme is adopted. 
\begin{figure}[htbp!]
\centering
\includegraphics[width=5.0in]{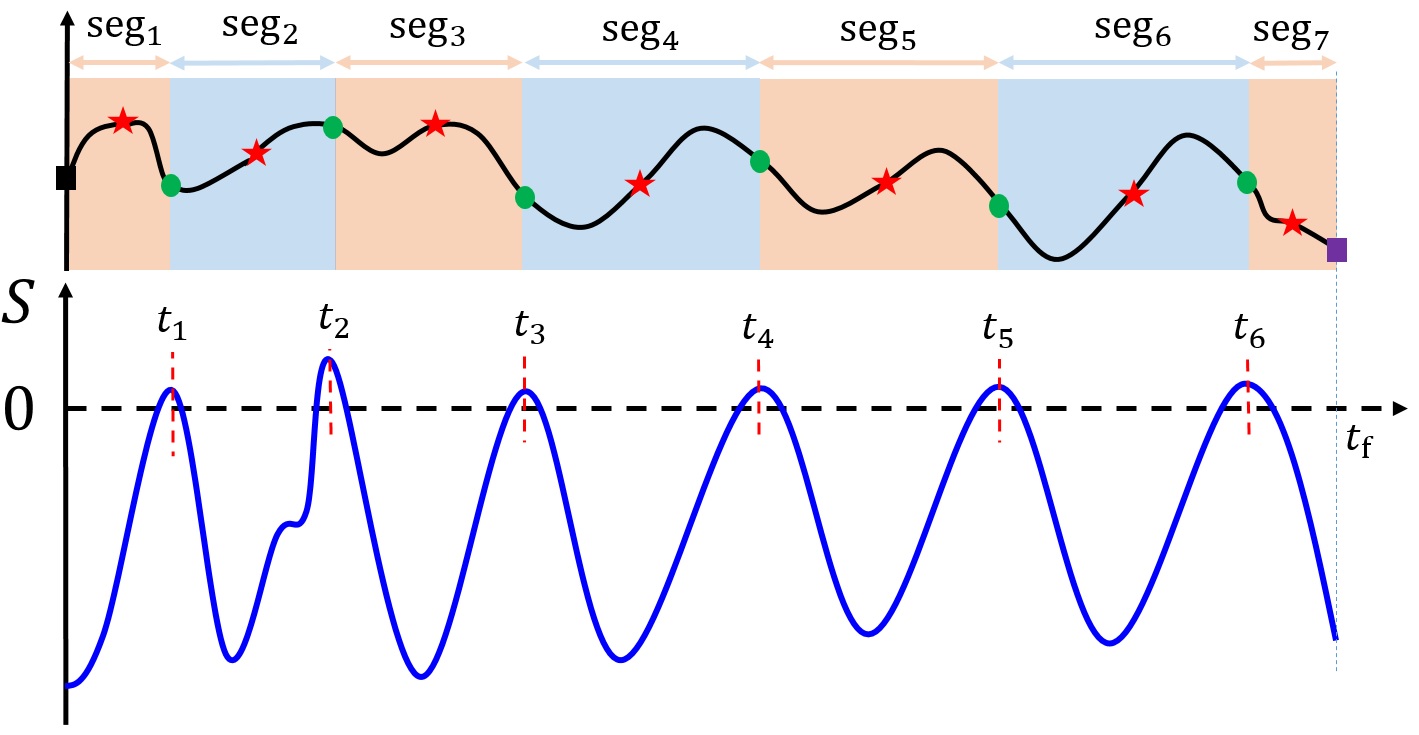}
\caption{Depiction of the switching function of multiple-revolution problems with $T >> T_{\text{min}}$ and numerical forward-backward scheme.}
\label{fig:NImpulseMany}
\end{figure}

Figure \ref{fig:NImpulseMany} shows the switching function of a representative multiple-revolution trajectory that consists of six impulses. The main reason for adopting such a scheme for impulsive optimization is the observation that impulse times (the position and velocity vectors at the associated impulse times) change near negligibly due to the use of switching surface to estimate these times precisely. In the majority of the test cases, qualitatively, these impulses are found by the converged indirect solution to be applied near peripasis/apoapsis (of the intermediate elliptical orbits) and/or the ascending/descending nodes of the initial and final orbits. Therefore, the time interval between consequent impulses (time duration of each segment), at most, corresponds to a complete revolution around the central body. This scheme allows us to utilize parallel computation, and improves the convergence. Also, we frequently invoke a high-fidelity force model at this stage, and the use of parallel computation is facilitated to improve wall-clock computational efficiency for final convergence. 

The trajectory is, therefore, divided into $M$ segments. In this example, seven segments $M=7$ (colored differently) in the upper part of the Figure \ref{fig:NImpulseMany} are considered. The times of intermediate impulses are denoted by $t_i$, $i = 1, \cdots, 6$. The discussion herein is for a solution in which the trajectory consists of intermediate impulses only (and no impulses occur at $t_0$ or $t_f$ denoted by squares in Figure \ref{fig:NImpulseMany}). However, the methodology is general and can handle impulses that are applied at the initial and final times, as well, should the switching function (at $T_{\text{max}}$) indicate terminal impulses. Except for the first and last segments, the beginning and the end of each segment consists of an impulse denoted by green circles.

Let $t^{-}_i$ and $t^{+}_i$ denote the time instants immediately before and after the $i^{\text{th}}$ impulse, respectively, the velocity vectors are similarly denoted as $\textbf{v}^{-}_i$ and $\textbf{v}^{+}_i$. At the moment of impulse, the position vector remains the same, i.e., $\textbf{r}^{-}_i = \textbf{r}_i = \textbf{r}^{+}_i$. Consider the fourth segment where its time duration is $\Delta t_{s,4} = t_4 - t_3$. States (position and velocity) are propagated forward (subscript `F') from $ t = t_3$ to $t = t_3 + \frac{\Delta t_{s,4}}{2}$ to get $\textbf{r}_F$ and $\textbf{v}_F$. Similarly, the states are propagated backward (subscript `B') from $t = t_4$ to $t = t_4 - \frac{\Delta t_{s,4}}{2}$ to get $\textbf{r}_B$ and $\textbf{v}_B$. 

The error between the states is used to from a residual vector (at the mid-point of the segment marked by a red star) 
\begin{align}
\bm{\Delta} = [\bm{\Delta}_{s,1}^{\top}, \cdots, \bm{\Delta}_{s,N}^{\top}]^{\top},
\end{align}
where $\bm{\Delta}_{s,i} \in \mathbb{R}^6$, $i=1, \cdots, N$ denotes the state residual vector at the mid-point of the $i^{\text{th}}$ segment and is defined as
\begin{align}
\bm{\Delta}_{s,i} = \begin{bmatrix}
\textbf{r}_F-\textbf{r}_B\\
\textbf{v}_F-\textbf{v}_B
\end{bmatrix}.
\end{align}
The matrix of decision variables is denoted by $\textbf{X}$ becomes
\begin{align} \label{eq:OpX}
\textbf{X} = \begin{bmatrix}
\textbf{r}(t_0), & \textbf{r}(t_1), & \cdots,  &\textbf{r}(t_6), &\textbf{r}(t_f)\\
\textbf{v}(t_0), & \textbf{v}(t_1), & \cdots,  &\textbf{v}(t_6), &\textbf{v}(t_f)\\
\Delta \textbf{v}(t_0), & \Delta \textbf{v}(t_1), & \cdots, & \Delta \textbf{v}(t_6), &\Delta \textbf{v}(t_f)\\
t_0, & t_1, & \cdots, & t_6, &t_f
\end{bmatrix},
\end{align}
where each impulse consists of ten decision variables (i.e., position, velocity, velocity impulse vectors and time of impulse). The initial and final variables can be fixed by setting the lower and upper bounds of them to be equal to the desired parameters. If impulses at the initial and final times have to be considered, the lower and upper bounds of the decision variables are modified accordingly. Ultimately, an optimization problem for the impulsive solution can be formulated as
\begin{equation} \label{eq:impulsiveopt}
\begin{aligned}
& \underset{\textbf{X}}{\text{minimize}}
& & \sum_{k=1}^{N} \| \Delta \textbf{v}_k \|, \\
& \text{subject to}
& & \bm{\Delta} = \textbf{0}.
\end{aligned}
\end{equation}
Any NLP solver (we have used MATLAB's \textit{fmincon}) chosen for minimizing the cost defined in Eq.~\eqref{eq:impulsiveopt} benefits from a good approximation of the time, direction and magnitude of the impulsive thrusts, which accelerates the convergence performance. These are precisely the information that we extract from the extremal field map, i.e., from the extremal associated with the high thrust limit of the switching surface thrust ridges. Moreover, due to the ensured quality of our starting estimate of the optimal maneuver, we do not have to guess the number of impulses. For instance, in Figure \ref{fig:NImpulseMany}, we can readily see that there exist six thrust arcs where good estimates of the velocity impulses can be obtained by using simple formula as
\begin{equation} \label{eq:approxformula}
\Delta v_i \cong \frac{T \Delta t_{s,i}}{m_i}, \; i = 1, \ldots, N,
\end{equation}
where $\Delta t_i$ denotes the thrust time interval, $T$ denotes the thrust level and $m_i$ denotes the mass at the mid-point, $t_i$ of the respective $i^{\text{th}}$ thrust interval. At each impulse time, $t_i$ the direction of the thrust is also known from the direction of the primer vector, $\textbf{p}(t_i) = - \mathbb{B}(\textbf{x}(t_i),t_i)^{\top} \bm{\lambda}(t_i)$. It is straightforward to calculate the impulse vector as $\Delta \textbf{v}_i = \Delta v_i \bm{\alpha}_i$ where $\bm{\alpha}_i = \textbf{p}(t_i)/\|\textbf{p}(t_i) \|$. All of the required values are retrievable from the extremal which is the solution of the TPBVP for $T \gg T_{\text{min}}$. Note that it is possible to parameterize the impulsive optimization problem by position-formulation (also known as Feasible Iterate Appraoch (FIA) [\citen{hughes2003comparison}]). The FIA parameterization uses Lambert problem and satisfaction of the position boundary conditions are guaranteed when two-body dynamics govern the motion. It reduces the number of design variables significantly. However, the method proposed in this paper is a general method applicable to beyond two-body dynamics. 

The analysis of these switching surfaces are best explained in the context of specific orbit transfers, while there are several generalized points of view that emerge, there are also specific features and behaviors that may or may not arise in the switching surface associated with a particular orbit transfer. So, we consider different cases to permit the diversity of behaviors to be explained. 

There are two key points that require explanation. First, the above construction is dependent on a reliable method to solve the underlying family of optimal control problems. For spacecraft trajectories, indirect methods are critical elements of the proposed procedure since, based on our experience, these methods can be significantly faster. However, the key point, they provide more rigorous and accurate optimal trajectories than the corresponding direct methods. Therefore, a detailed discussion is devoted to an enhanced process to solve the TPBVPs that arise in the indirect formulation of optimal control problems. It is also possible to construct these surfaces using any type of direct optimization method. The second important point is related to the fact that the minimum-thrust trajectory (also, because a duality exists, this minimum-thrust trajectory is also a minimum-time trajectory if the corresponding thrust is judiciously specified) is the base solution from which the computation of the switching surface is initiated. 

On the other hand, minimum-time and minimum-fuel trajectories typically have one local extrema for each of a number of $N_{\text{rev}}$ en route revolutions in the orbit transfer (we emphasize that we are dealing with rendezvous-type maneuvers), and the number of revolutions (we will see) affects the structure of the switching surface. Therefore, a reliable strategy is required to find the fundamental minimum-thrust solution, which implicitly requires us to find the local extremals associated with multiple revolutions. Note, when the constant thrust is always `on', the minimum thrust extremal is also the minimum-fuel extremal. This process is somewhat analogous to finding all solutions in the multiple-revolution $2$-impulse Lambert problems [\citen{prussing1970optimal,battin1999introduction,woollands2017multiple}]. The details of an algorithm for finding fundamental minimum-thrust solutions are explained in the next section. The final observation is that the methods of this section are well suited for parallel computation, which will facilitate efficient computations when high-fidelity force models are used for final convergence. 

\section{Procedure for Finding Minimum-Thrust Solution} \label{sec:minthrustalgorithm}
As discussed earlier, the minimum-thrust solution is, in fact, nothing but the minimum-time solution for the prescribed boundary conditions, and time of flight. However, minimum-time solution is not unique. We seek the thrust $T_{\text{min}}$ for which not only the minimum time ($t^*_f - t_0$) is equal to the desired maneuver time ($t_f - t_0$) but also requires the least amount of propellant. Therefore, a procedure is devised to find the solution to the minimum-thrust solution, which is based on the formulation of minimum-time trajectory optimization problem. For minimum-time problem the state/co-state dynamics is the same as those in Eqs.~\eqref{eq:elementsdiffeq}-\eqref{eq:masscostate}. The optimal control vector is known [\citen{taheri2017co}] and is characterized by
\begin{align}
\bm{\alpha}^* = & \frac{\textbf{p}}{||\textbf{p}||}, \\
\delta^{*} = & 1.
\end{align}
Note that the switching function of the minimum-time problem has a different mathematical expressions and is known to remain non-negative along the entire trajectory, i.e., $S=\frac{c ||\mathbb{B}^{\top} \pmb{\lambda}||}{m}+\lambda_m -1> 0$. In practice, mass and its co-state associated can be omitted from the numerical analysis; however, we kept this formulation since we can use the vector of converged solution to start the minimum-fuel continuation procedure. Since the terminal time, $t_f$ is free, optimality conditions require the following condition on the final value of the Hamiltonian, $H^*(t_f)=0$. On the other hand, neither the state equations, cost functional, nor the terminal constraints depend on time explicitly which means that the Hamiltonian is a constant along the optimal trajectory, i.e., $H^*(t)=0$. Therefore, for a rendezvous-type maneuver and its associated boundary conditions (that we have considered in this paper), the vector of terminal constraints become 
\begin{equation}\label{eq:terminalconmintime}
\bm{\psi}(\textbf{x}(t_f),t_f) = \left [ [\textbf{x}(t_f) - \textbf{x}_T]^{\top},\lambda_m(t_f),  H(t_f) \right ]^{\top} = \textbf{0}.
\end{equation}
The TPBVP associated with the minimum-time problem consists of vector $\Theta = [\bm{\lambda}^{\top}(t_0),\lambda_m(t_0),t_f]^{\top}$ with eight unknown values and the vector of terminal constraints is given in Eq.~\eqref{eq:terminalconmintime}.

Our goal is to determine the minimum-thrust solution. In order to find a solution, the above TPBVP is augmented with thrust magnitude as one additional variable, $T$ and we also augment the vector of final constraints with an additional equality constraint, i.e., $t_f^* - t_f = 0$, where $t_f^*$ is the minimum time of flight. The inclusion of this equality constraint is crucial to guide the solution toward the minimum-thrust magnitude, $T_{\text{min}}$. Therefore, the augmented (subscript `a') design vector of the optimization problem is $\Theta_a = [\Theta^{\top},T]^{\top}$ and the augmented vector of final constraints becomes
\begin{equation}
\bm{\psi}_a(\textbf{x}(t_f),t_f) = \left [ [\textbf{x}(t_f) - \textbf{x}_T]^{\top},\lambda_m(t_f),  H(t_f),  t_f^* - t_f \right ]^{\top} = \textbf{0}.
\end{equation}
In the above optimization problem, a good estimate for the time of flight is known! In fact, the prescribed time of flight is the desired minimum time solution, corresponding to the $T_{\text{min}}$ for which the sought extremal is the minimum time maneuver. However, a good estimate of the thrust magnitude will enhance the convergence performance of any chosen solver. \\
A simple numerical procedure is outlined to provide an estimate for the thrust, which is based on work-energy principle. The work/energy principle states that that for a particle the work done by all forces equals the change in the kinetic energy, which in our problem can be written as
\begin{equation}\label{eq:workenergy}
\frac{1}{2}m_f v^2_f-\frac{1}{2}m_0 v^2_0 = \int_{\textbf{r}_0}^{\textbf{r}_f} -\frac{\mu m}{r^3} \textbf{r} \cdot d \textbf{r} + T \int_{t_0}^{t_f} \bm{\alpha} \cdot \textbf{v}~dt,  
\end{equation}
where $v_f = \| \textbf{v}_f \|$ and $v_0 = \| \textbf{v}_0 \|$. In addition, the final mass, $m_f$ is related to the initial mass, $m_0$ through $m_f = m_0 - \frac{T}{c} (t_f-t_0)$. The second integral on the right-hand side of Eq.~\eqref{eq:workenergy} is not straightforward to evaluate because it depends on the unknown path and the associated optimal steering direction vector for thrust. However, it is known that the maximum change in the kinetic energy is achieved when the thrust is aligned along or against the velocity vector, i.e., $\bm{\alpha} = \pm \textbf{v}/\| \textbf{v}\|$. This fact is frequently used as a warm start for low-thrust trajectory optimization [\citen{taheri2012shape,taheri2015fast}]. Therefore, the second integral can be approximated as 
\begin{equation} \label{eq:secondint}
\int_{t_0}^{t_f} \frac{\textbf{v}}{\| \textbf{v}\|} \cdot \textbf{v}~dt = \int_{t_0}^{t_f} \| \textbf{v}\|~dt = \pm \bar{r} \dot{\bar{\theta}} (t_f-t_0),
\end{equation}
where $\bar{r} = \frac{r_0+r_f}{2}$ and $\dot{\bar{\theta}} = \sqrt{\frac{\mu}{\bar{r}^3}}$ are mean radius and mean angular velocities, respectively, and it is assumed that $\|\textbf{v}\| \approx r \dot{\theta}$, which neglects the non-planar, and radial components of velocity, which are usually small but definitely not negligible. So this assumption will typically give a smaller than optimal variation in the orbit due to thrust, or put another way, would lead to a similar large thrust to accomplish the change in kinetic energy. Over-estimating the thrust is preferred, because larger than minimum thrust still leads to feasible solutions and thrust can be reduced until the switching function just touches zero at one point to identify the desired thrust $T_{\text{min}}$. After all, the result of this simplifying approach will be used as an initial guess for the actual minimum-thrust optimization problem, which justifies ``reasonable'' simplifications. 
%However, the assumptions of constant mass ($= m_0$) in the first integral effect with the result that the computed thrust solved from Eq.~\eqref{eq:workenergy} will be too small or too large, in the examples tested. We find that 30 \% errors are not unusual.
Upon substitution of Eq.~\eqref{eq:secondint} into Eq.~\eqref{eq:workenergy} and evaluating the first integral (that leads to $\frac{\mu m_f}{r_f}-\frac{\mu m_0}{r_i}$ assuming crudely that $m$ is evaluated at the terminal points), one can solve for an estimated value for the thrust using the following relation 
\begin{equation} \label{eq:thrustest}
T_{\text{estimate}} = \frac{\frac{m_0}{2} \left ( v^2_f - v^2_0 \right )+\mu m_0 (\frac{1}{r_0}-\frac{1}{r_f})}{\frac{t_f-t_0}{c} \left ( \frac{v^2_f}{2}+\frac{\mu}{r_f} \right ) \pm \bar{r} \dot{\bar{\theta}} (t_f-t_0)}.
\end{equation}
The proper sign of $\pm$ is determined by the fact that for trajectories to more (less) energetic orbits, the energy has to increase (decrease). The thrust obtained through Eq.~\eqref{eq:thrustest} can be used as a good initial guess for minimum-thrust optimization problem. Table \ref{tab:minthrust} shows the details of the results of the above optimization algorithms for finding minimum-thrust magnitude for the first two test cases under two-body dynamical model. MATLAB \textit{fsolve} is used for solving the TPBVP associated with the minimum-thrust optimization problem. 
\begin{table}[h!] 
\begin{center} 
		\caption{Details of minimum-thrust solutions for the first two test cases.}\label{tab:minthrust}
		{\small%\scriptsize
		\begin{tabular}{c c c c}
        \hline
        \hline
         Problem & $T_{\text{estimate}}$ & $T_{\text{min}}$ & CPU time \\
                 &          [N]        & [N]            & [s] \\
         \hline
         Earth-to-1998ML & 0.141        &  0.1265  & 0.25 \\
         Earth-to-Venus  & 0.111        &  0.0663  & 8.42\\
        \hline
        \hline
        \end{tabular}
		}
	\end{center}
\end{table}
We should add that we have also tried Edelbaum's method [\citen{edelbaum1961propulsion}] that establishes a relation between $\Delta v$, time of flight and thrust level for circular to circular orbit transfers. Since the considered test cases are rendezvous maneuvers, the amount of required $\Delta v$ that has to be used in Edelbaum's relation has to be increased to take into account the additional required energy. Our experience shows that $1.5 \times \Delta v$ leads to convergence for the considered problems under two-body dynamics.

The number of revolutions is a factor that has to be considered during the procedure of solving minimum-thrust optimization problem. The projection of the initial and target position vectors onto the $x-y$ plane of an inertial frame make angles, $\theta_0 \in [0, 2 \pi]$ and $\theta_T \in [0, 2\pi]$, respectively, with respect to the $x$ axis. Without loss of generality and assuming that $\theta_T > \theta_0$, the difference between the angles is denoted by $\Delta \theta = \theta_T - \theta_0$. The number of revolutions is considered to update the final angles through $\theta_f = \Delta \theta+ 2\pi \times N_{\text{rev}}$. A rigorous way to define $N_{\text{rev}}$ is to count the number of successive piercing of the trajectory through the plane defined by the initial orbit radius and orbit normal. Our experience shows that the global ``optimal'' minimum-fuel (and also minimum-thrust) solution corresponds to a unique number of revolutions required to achieve the maneuver. Therefore, the minimum-thrust problem has to be sought starting from a lower number of revolutions, $N_{\text{rev},l}$. 

For two-body problems, let $\tau_l$ and $\tau_u$ denote the smaller and larger values of the orbital periods of the involved bodies bounded from below and above according to the following relation
\begin{equation} \label{eq:Nrevbound}
N_{\text{rev},l} = \max \{ \text{floor}(\frac{t_f - t_0}{\tau_u}-1),0 \} \leq N_{\text{rev}} \leq N_{\text{rev},u} = \text{ceil}(\frac{t_f - t_0}{\tau_l}+1),
\end{equation}
where $N_{\text{rev},l}$ and $N_{\text{rev},u}$ denote the lower and upper bounds for the number of revolutions, respectively, and \textit{ceil} and \textit{floor} operators return the next larger (next smaller) integer number from their arguments. This formula neglects that truth that thrusting, even low thrusting, alters the two body period of the transfer vehicle from the starting orbital period. We specify this in the perturbed true anomaly desired for the optimal maneuver, but as we near convergence, the physical number of revolutions is defined more rigorously as the number of actual passages of the transfer orbit through the reference plane defined by the initial orbit normal and the initial position vector. Except where the starting and arrival position vectors are nearly co-linear, Eq.~\eqref{eq:Nrevbound} holds. Put another way, $N_{\text{rev}}$ revolutions using the period of the starting orbit is not guaranteed to correspond to $N_{\text{rev}}$ revolutions of the perturbed transfer orbit, but any discrepancy is easily cleared up by counting revolutions.
\begin{figure}[htbp!]
\centering
\includegraphics[width=5.0in]{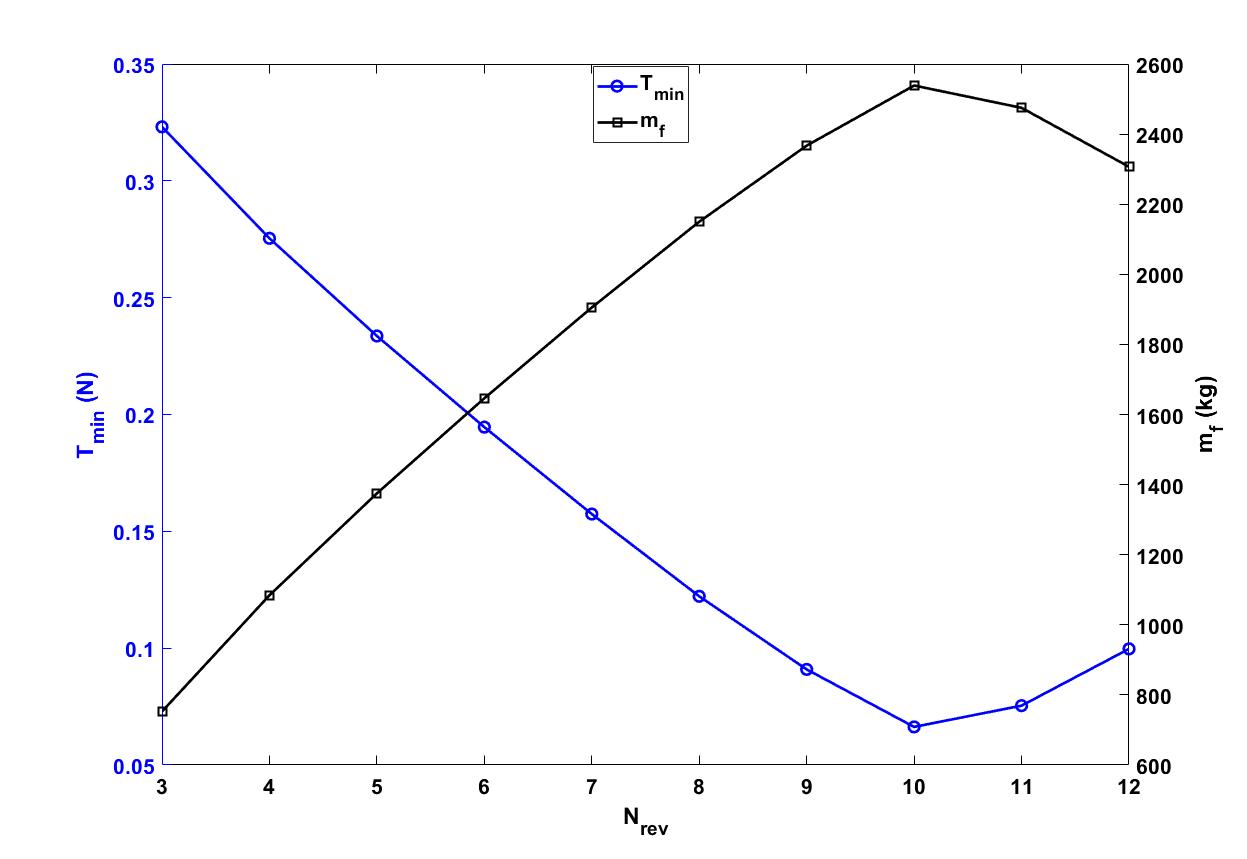}
\caption{Changes in $T_{\text{min}}$ and $m_f$ vs. $N_{\text{rev}}$ for Earth-to-Venus problem.}
\label{fig:EVMinThrustNrevs_early}
\end{figure}
Starting from $N_{\text{rev},l}$, if convergence to satisfaction of the optimal control necessary conditions is achieved, the solution is the minimum-thrust value; otherwise, the value of the selected number of en route revolutions should be incremented by one. The process is repeated until convergence is achieved. Alternatively, all of these extremals, can be solved simultaneously, and the one with the smallest value of the propellant mass consumed will be the global extremal, i.e., the fundamental minimum-thrust solution. Note that numerical results indicate the there exist solutions with $N_{\text{rev}}>N_{\text{rev},u}$, but they are typically not of practical interest since they lead to sub-optimal solutions with large fuel consumption. On the other hand, there exist no solutions when $N_{\text{rev}}<N_{\text{rev},l}$; reachability of the given target state requires a specific minimum number of revolutions which may be $N_{\text{rev}} = 0,1,2,\cdots.$ Technically, of course, the solution is generally reached in $N_{\text{rev}}$ plus a fraction of the $N_{\text{rev}} + 1$ revolution.

For instance, for the Earth-to-Venus problem (which is one of the test cases studied in this paper), Figure \ref{fig:EVMinThrustNrevs_early} shows the changes in minimum-thrust and its associated final mass for different number of en route revolutions. Clearly, the fundamental minimum-thrust solution corresponds to $N_{\text{rev}} = 10$. For $N_{\text{rev}} \geq 11$ the trajectories make un-necessary revolutions by getting closer to the Sun. Note that all of these solutions are minimum-time solutions ($\delta^*(t) = 1$) but with different thrust levels. Figure \ref{fig:EVMinThrustNrevs_early} is obviously useful for sizing and mission planning purposes.
\section{Results}\label{sec:results}
Five minimum-fuel trajectory optimization problems (under two-body dynamics) are considered where the strength and nonlinearity of the gravitational field varies significantly. As a consequence, distinct topological features in the respective switching surfaces appear, in particular, for multi-revolution trajectories. %In addition, a GTO to Halo orbit around L1 point in restricted three-body dynamics of the Earth-Moon systems is also studied.

For interplanetary cases, the canonical units are adopted to normalize the state and control inputs where one distance unit (DU) is equal to the astronomical unit (AU), and $2\pi \times$ Time Unit (TU) is 1 year. In the numerical simulations, the gravitational parameter of the Sun is set to $\mu_{\odot} = 132712440018$ km$^3$/s$^2$, $g_0 = 9.8065$ m/s$^2$, whereas the gravitational parameter of the Earth is set $\mu_{\oplus} = 398600$ km$^3$/s$^2$. For Earth-bound problems, one distance unit (DU) is equal to the Earth radius at the equator, $R_{e} = 6378$ km, and Time Unit (TU) is 806.78557 seconds, $\text{TU} = R_e^{3/2}/\sqrt{\mu_{\oplus}}$. 
In both cases, the TU is the time for a satellite in the circular reference orbit to move through one radian of true anomaly. 

\subsection{Interplanetary Rendezvous From Earth to Mars}
The first test case is chosen similar to the first interplanetary Earth-to-Mars transfer case in [\citen{zhu2017solving}] in order to validate the results and to present more insights into the ensemble of solutions as intended, in the light of the switching surface. The same BCs are taken along with the specified time of flight, $t_f-t_0 = 793$ days. The following values are considered for the parameters of the spacecraft and its low-thrust propulsion system: $m_0 = 2000$ kg,  and $I_{\text{sp}} = 3000$ s. Some fraction of $m_0$ is propellant. 

For a fixed dry mass, by maximizing the final spacecraft mass, we minimize the fuel required and maximize the useful payload we can deliver to the final state. The value of thrust is considered as the sweeping parameter, to generate an infinite family of optimal maneuvers, i.e., $T \in [T_{\text{min}},T_{\text{u}}]$, where $T_{\text{u}}$ is some upper value of thrust and is set to $T_{\text{u}} = 10$ N. Several slices of the switching surface will be studied in more details to introduce particular concepts the may re-appear in other switching surfaces. 

The family of state and co-state variables that underlie the switching surface constitute an extremal field map. Obviously, this extremal field map, for an infinite family of maximum thrust values has immediate utility in sizing of propulsion systems for mission design purposes, which will be explained. The Earth position and velocity vectors at the departure time, $t_0$ are 
\begin{align*}
\textbf{r}_{\text{E}} = & [58252488.0107,135673782.5313, 2845.0581]^{\top} \text{km},\\
\textbf{v}_{\text{E}} = & [-27.8445,11.6599,0.0003]^{\top} \text{km/s}.
\end{align*}
As is usual for preliminary design of solar missions, we assume we are just outside the Earth's sphere of influence at departure and Mars sphere of influence at arrival, i.e., we ignore Earth and Mars gravity. The position and velocity vectors of Mars at the final time, $t_f$ are 
\begin{align*}
\textbf{r}_{\text{M}} = & [36216277.8004, -211692395.5225, -5325189.0499]^{\top}\text{km},\\ 
\textbf{v}_{\text{M}} = & [24.7988,6.1682,-0.4800]^{\top} \text{km/s}. 
\end{align*}
The minimum-$\Delta v$ 2-impulsive solution can be obtained by solving the corresponding Lambert problem. The magnitude of the impulses at the initial and final time instants are $\Delta v(t_0) = 3.0157$ km/s and $\Delta v(t_f) = 3.0318$ km/s, which correspond to a solution with one revolution around the Sun, i.e., $N_{\text{rev}} = 1$. 

\begin{figure}[htbp!]
\centering
\includegraphics[width=4.0in]{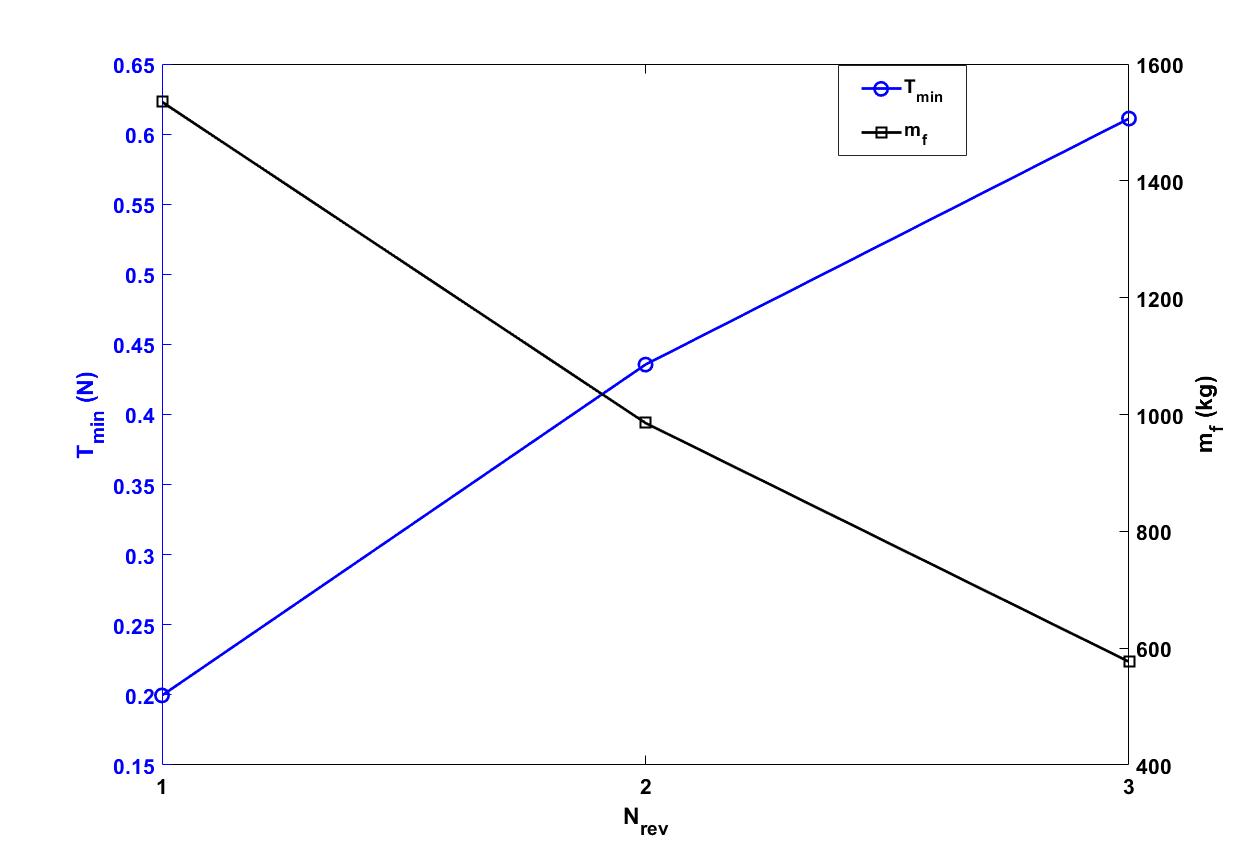}
\caption{Changes in $T_{\text{min}}$ and $m_f$ vs. $N_{\text{rev}}$ for the Earth-to-Mars problem.}
\label{fig:EMMinThrustNrevs}
\end{figure}
Before generating the switching surface, we need to determine the optimal number of en-route revolutions (and its associated minimum-thrust solution) based on the algorithm outlined earlier. Figure \ref{fig:EMMinThrustNrevs} shows the changes in $T_{\text{min}}$ and $m_f$ vs. the feasible values for $N_{\text{rev}}$. The critical value of thrust for the fundamental minimum-thrust (or minimum-time) for the given BCs is found to be $T_{\text{min}} = 0.1996$ N. Thus, a relatively low thrust can send a significant payload to Mars, but the time of flight is significant. The fundamental minimum-thrust solution corresponds to $N^*_{\text{rev}}=1$, which indicates that the maneuver completes one plus a fraction of the second revolution along the transfer trajectory. The initial phase angle between the position vectors is small, but $N_{\text{rev}} = 0$ does not lead to a solution. Put it another way, with the given parameters of the propulsion system and BCs, the amount of propellant required for $N_{\text{rev}} = 0$ is larger than the initial mass of the spacecraft.

\begin{figure}[htbp!]
\centering
\includegraphics[width=5.0in]{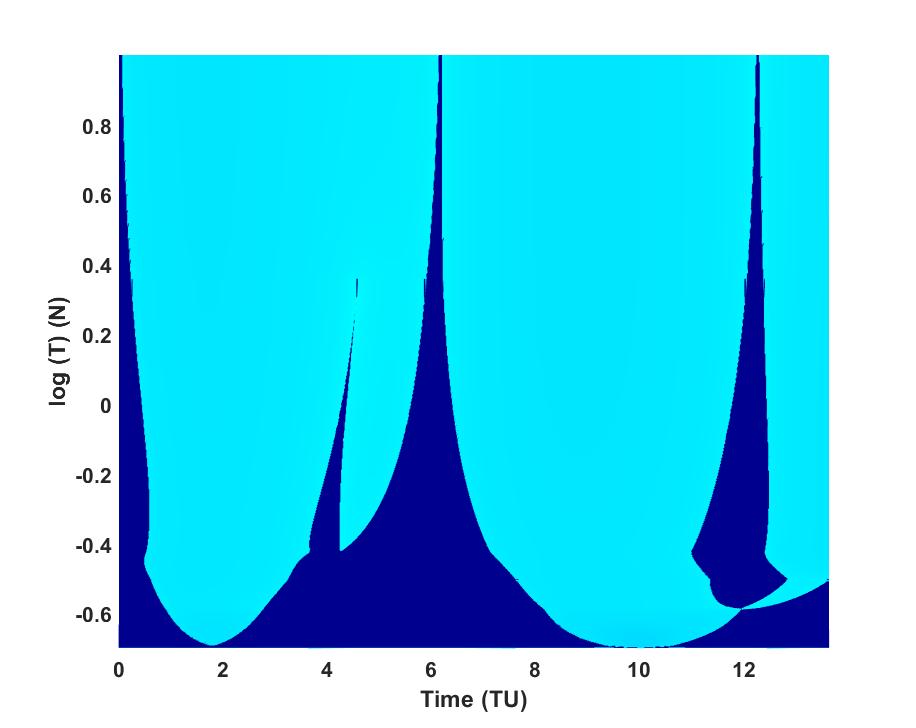}
\caption{Switching surface $S = 0$ contour map for the Earth-to-Mars problem for $N^*_{\text{rev}} = 1$.}
\label{fig:EMSS}
\end{figure}

\begin{figure}[htbp!]
\centering
\includegraphics[width=4.5in]{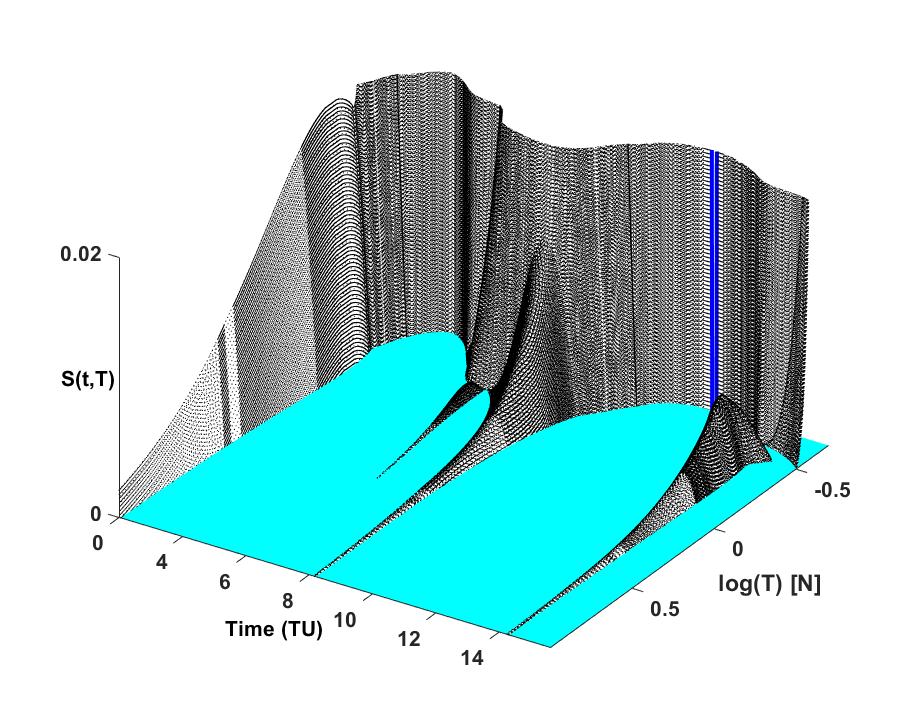}
\caption{Enlarged view of the switching surface for the Earth-to-Mars problem with $N^*_{\text{rev}} = 1$.}
\label{fig:EMSS_log_2}
\end{figure}

We found the same minimum-thrust magnitude for the given boundary conditions as in [\citen{zhu2017solving}], i.e., $T_{\text{min}} = 0.1996$ N. The thrust value is then swept in the given range $T \in [0.1996, 10]$ N.  Figure \ref{fig:EMSS} shows the top-view of the switching surface generated by sweeping the thrust magnitudes where dark blue regions denote thrust ridges, whereas the light blue regions denote coast canyons ($S<0$). At first glance, the contour plot seems to consists of three main thrust ridges ($S>0$) with wide bases (for very low-thrust) that all ridges have a tapering trend up to the top of the curve as thrust magnitude increases. 

This surface can be viewed in many ways, which reveals a number of interesting and illuminating facts. There is a significant number of changes in the topology of the switching surface that occur at the lower part of the plot, in the region of very low-thrust magnitudes. These changes are due to the local changes in the switching function and gradual passages of interesting local features through $S = 0$ as explained in the previous section.  

\begin{figure}[htbp!]
\centering
\includegraphics[width=5.0in]{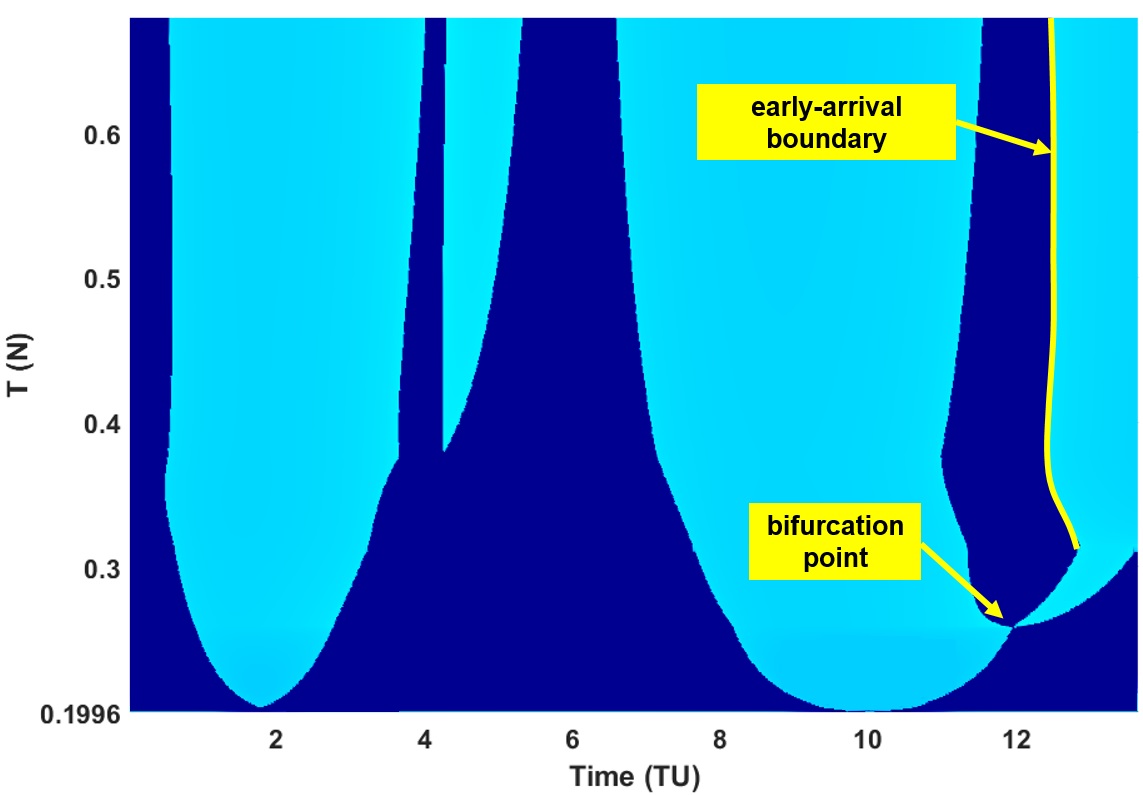}
\caption{Low-thrust region of the switching surface for the Earth-to-Mars problem along with the bifurcation point and the early-arrival boundary.}
\label{fig:SSEM_lower}
\end{figure}

The second region is associated with the medium thrust values where there are only four thrust ridges. Any given thrust value corresponds to a horizontal profile (slice) of this surface, which is the switch function for that particular maximum thrust level. As the thrust magnitude is increased, there is a slender dagger-like thrust ridge between the first two main thrust ridges that vanishes as thrust is increased. Beyond this thrust magnitude, $T>2.2682$ N the whole family of optimal minimum-fuel trajectories are characterized by three thrust ridges and the time duration (width) of these thrust ridges keep shrinking as thrust magnitude increases. It is a trivial observation that these three thrust ridges are tending toward three impulsive thrusts for this case. Figure \ref{fig:SSEM_lower} shows a ``blow-up'' enlarged view of the lower thrust region so we can discuss the changes in the topology of the switching surface.

We mention without proof that the center-line location of the ``persistent high thrust ridges'' $S = 0$ contour is very insensitive to $T/m$ and $I_{\text{sp}}$ as $T$ becomes large. Figure \ref{fig:EMSS_log_2} shows the 3D \textit{fundamental} switching surface for the Earth-to-Mars problem. Note that an opposite-angle view is chosen so that the high (``mountainous'') region of the surface does not hide the interesting features. On the other hand, the difference of the values of the switching function at low- and high-thrust levels is sufficiently small that it is difficult to see the actual 3D features of the switching surface. Therefore, an enlarged view ($S < 0.02$) is given in Figure \ref{fig:EMSS_log_2} for demonstration purposes so the low thrust details are more visible. The switching function of the fundamental minimum-thrust solution (with $T_{\text{min}} = 0.1996$ N) is shown using solid blue line (only a portion of the fundamental minimum-thrust solution is visible in the enlarged view).

%The second region is associated with the medium thrust values where there are only four thrust ridges. Any given thrust value corresponds to a horizontal profile (slice) of this surface, which is the switch function for that maximum thrust level. As the thrust magnitude is increased, there is a slender dagger-like thrust ridge between the first two main thrust ridges that vanishes as thrust is increased. Beyond this thrust magnitude, $T>2.2682$ N the whole family of optimal minimum-fuel trajectories are characterized by three thrust ridges and the time duration (width) of these thrust ridges keep shrinking as thrust magnitude increases. It is a trivial observation that these three thrust ridges are tending toward three impulsive thrusts for this case. Figure \ref{fig:SSEM_lower} shows a ``blow-up'' enlarged view of the lower thrust region so we can discuss the changes in the topology of the switching surface. 

The switching surface $S = 0$ contour map is also shown in Figure \ref{fig:EMSS} where dark blue regions denote ($S > 0$) thrust ridges, whereas the light blue regions denote ($S<0$) coast canyons. We draw your attention to a bifurcation that occurs at $t = 12$ TU and log(T) $\approx -0.587$ N. This bifurcation phenomenon (``peninsula'' form) results in the creation of a thrust ridge if thrust is slightly increased. These thrust ridges are not due to a branching phenomenon of a core thrust ridge (see, for instance, the dagger-like thrust ridge near $t = 4$ TU, which appears near log(T) $\approx -0.39$ N and vanishes around log(T) $\approx 0.39$ N). We should mention logarithmic scales are used for the thrust axis to magnify the lower region of the switching surface where a number of significant changes occur. Another important point is the fact that beyond a critical thrust magnitude, all of the minimum-fuel trajectories consist of a final coast phase. 
\begin{figure}[htbp!]
\begin{multicols}{2}
\centering
\includegraphics[width = 0.44\textwidth]{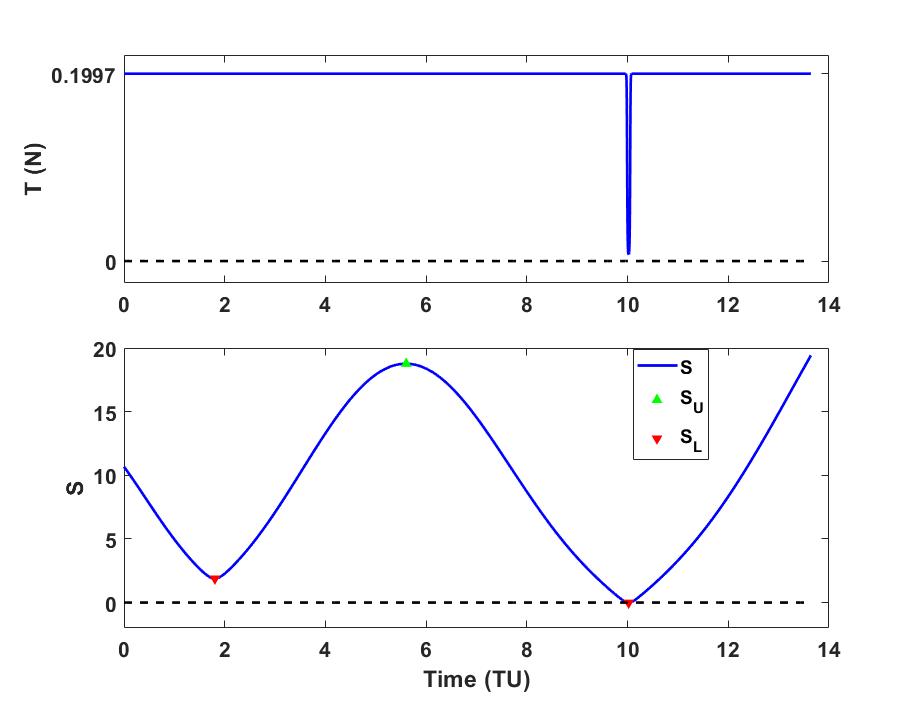}
\caption{Earth-to-Mars switching function with \\$T_{\text{min}} \approx 0.1997$ N.}
\label{fig:EMCT1_SF}
\centering
\includegraphics[width = 0.45\textwidth]{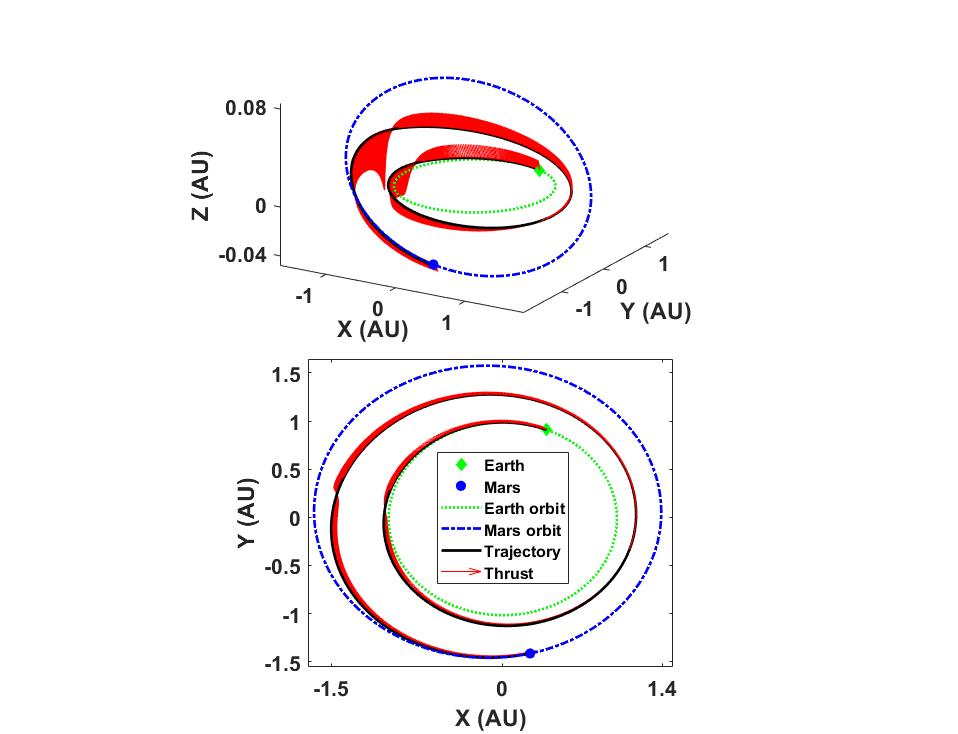}
\caption{Earth-to-Mars heliocentric trajectory with $T_{\text{min}}$.}
\label{fig:EMCT1_Traj}
\end{multicols}
\end{figure}
In other words, the right-most points of the last thrust ridge (see Figure \ref{fig:SSEM_lower}) define what we describe as an \textit{early-arrival} boundary. Clearly, the early-arrival boundary has an interesting profile especially near the bifurcation point. Although this condition does not exist in the switching surface of the Earth-to-Mars problem, in general, we should anticipate the existence of \textit{late-departure} boundaries that correspond to the existence of an initial coast phase. 

We proceed by inspecting the solution and switching function of a number of particular slices of the switching surface. Figure \ref{fig:EMCT1_SF} shows the switching function for minimum-thrust magnitude where the switching function is non-negative except for an individual point where it osculates (kisses) the $S = 0$ line. 

\begin{figure}[htbp!]
\begin{multicols}{2}
\centering
\includegraphics[width = 0.5\textwidth]{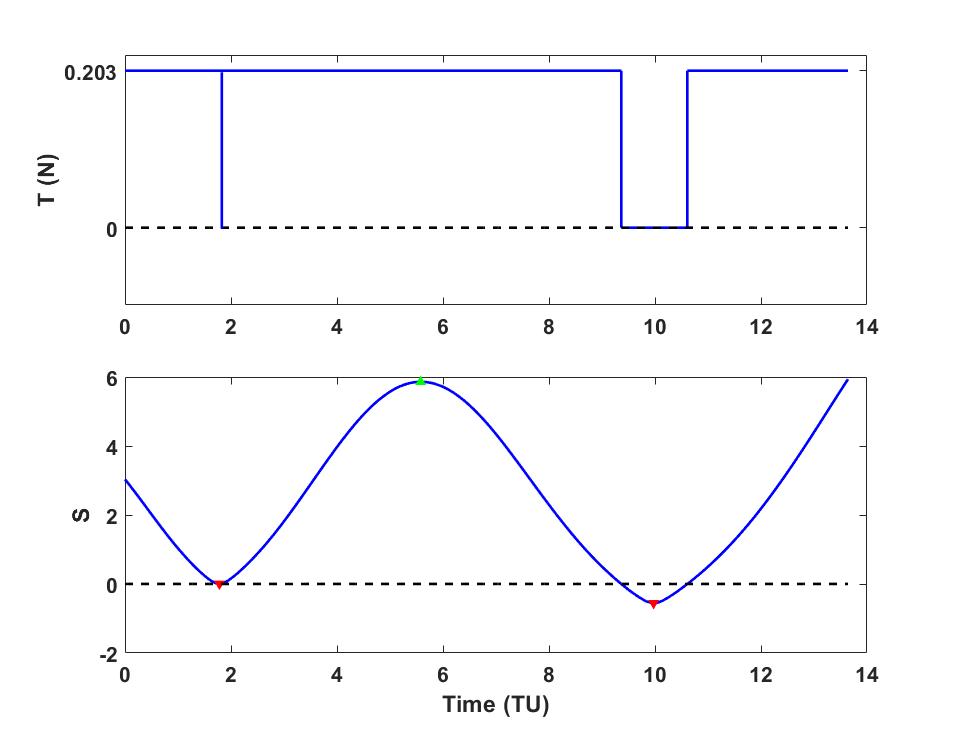}
\caption{Earth-to-Mars switching function with \\$T \approx 0.203$ N.}
\label{fig:EMCT2_SF}
\centering
\includegraphics[width = 0.5\textwidth]{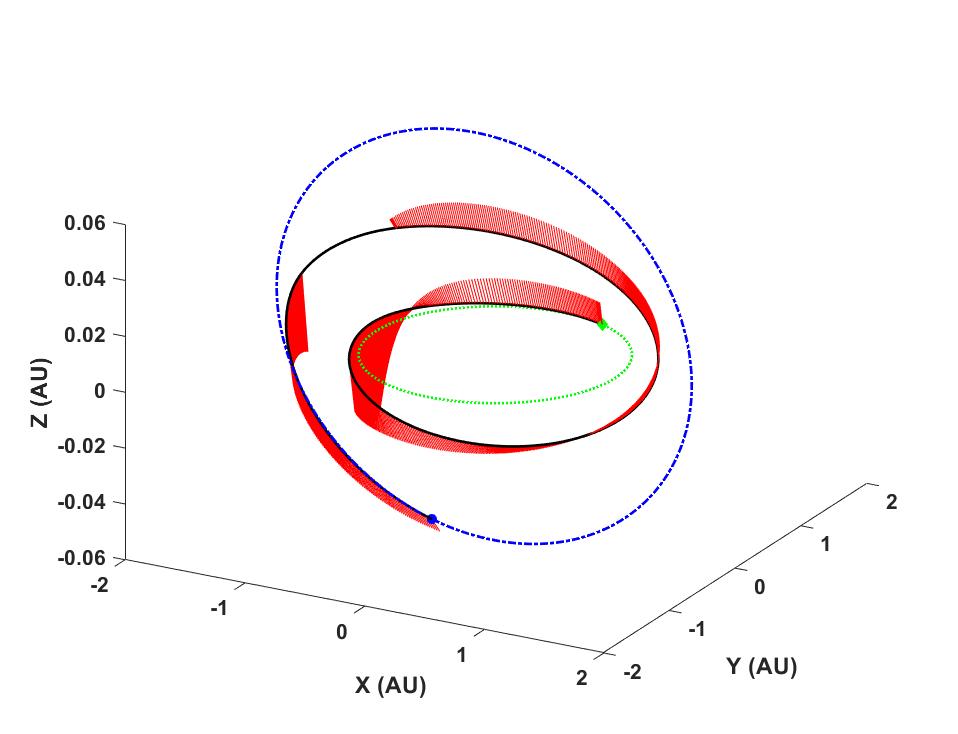}
\caption{Earth-to-Mars heliocentric trajectory with $T \approx 0.203$ N.}
\label{fig:EMCT2_Traj}
\end{multicols}
\end{figure}

Figure \ref{fig:EMCT1_Traj} shows the heliocentric trajectory of the minimum-thrust case. Vertical scale of Figure \ref{fig:EMCT1_Traj} greatly exaggerates the out-of-plane motion to emphasize that the maneuver is fully three dimensional. It reveals the existence of a number of interesting phenomena. The first coast arc, obviously, is centered around the time instance of 10 TU. However, as the thrust magnitude increases a second coast appears close to 1.9 TU. Figures \ref{fig:EMCT2_SF} and \ref{fig:EMCT2_Traj} show the switching function, thrust profile and heliocentric trajectory for this critical thrust, $T \approx 0.203$ N. 
There are two interesting phenomena that occur as the thrust magnitude is increased. The first one is a bifurcation, i.e., creation of a new thrust arc. Note, that this is the third thrust ridge (the third main thrust ridge among the total three thrust ridges) that does not disappear as the thrust magnitude is increased. Unlike the other two main thrust ridges that have a wide root at the low thrust beginning, this thrust ridge is created at a higher thrust magnitude where a corner is easy to distinguish. 
%This act of birth of a new thrust region points also at the impending disappearance of the very last thrust arc and a ``terminal coast'' region. 

\begin{figure}[htbp!]
\begin{multicols}{2}
\centering
\includegraphics[width = 0.5\textwidth]{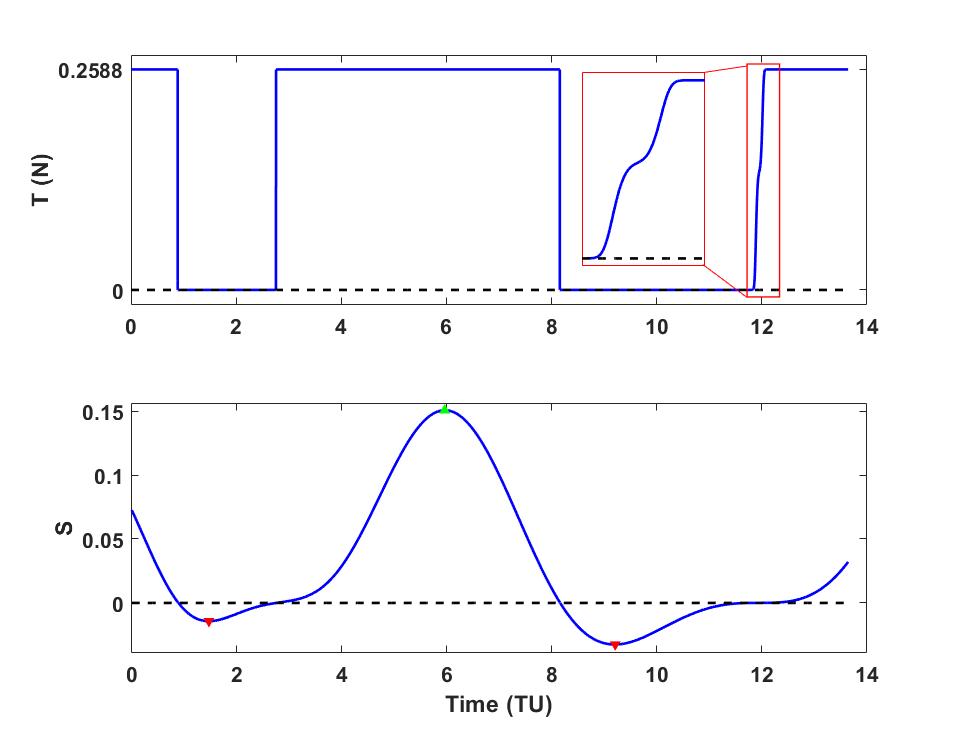}
\caption{Earth-to-Mars switching function with \\$T \approx 0.2588$ N.}
\label{fig:EMCT3_SF}
\centering
\includegraphics[width = 0.50\textwidth]{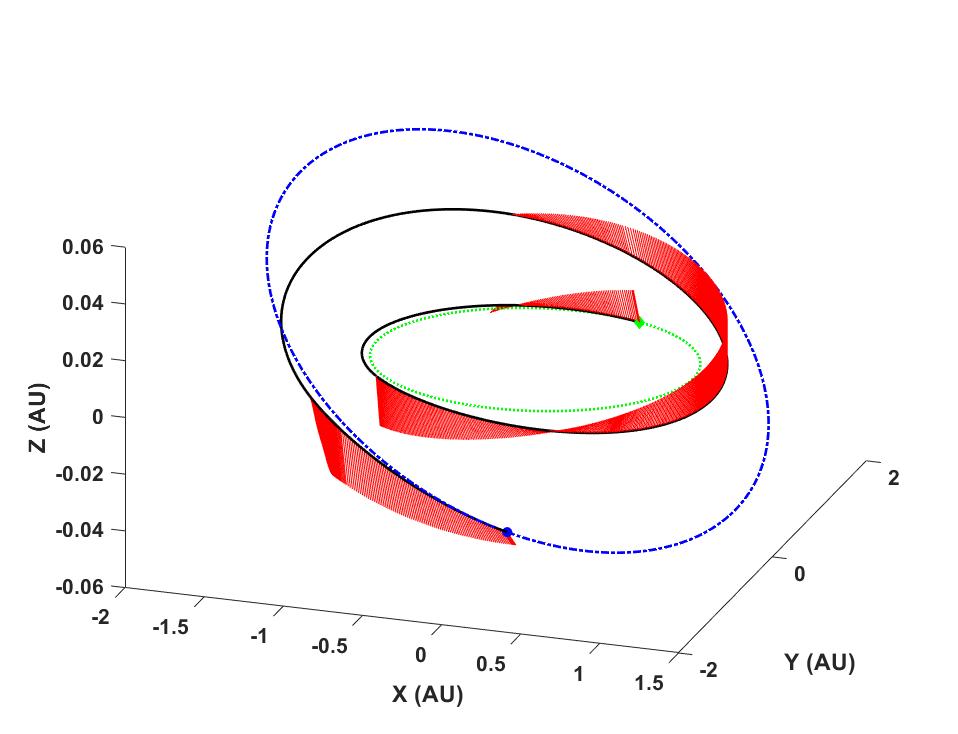}
\caption{Earth-to-Mars heliocentric trajectory with $T \approx 0.2588$ N.}
\label{fig:EMCT3_Traj}
\end{multicols}
\end{figure}

Figure \ref{fig:EMCT3_SF} appears to show a relatively flat profile of $S$ over a finite time interval, which triggers an alarm for the possible presence of a singular arc. In fact, this bifurcation point can be characterized by confirming three identities: $S = 0$, $\dot{S} = 0$, and $\ddot{S} = 0$ over a small finite time interval. Figures \ref{fig:EMCT3_SF} and \ref{fig:EMCT3_Traj} show the switching function, thrust profile and heliocentric trajectory for this \text{unique} critical thrust, $T \approx 0.2588$ N. 
Note that the thrust profile takes intermediate values over a finite but small time interval. This fact is directly related to the fact $S = 0$ over a finite time interval in the vicinity of 12 TU, i.e, we have detected the presence of a singular arc. This is the singularity reported in [\citen{zhu2017solving}], but not identified as a near-singular thrust arc. Insofar as is known, this is the first singular thrust arc found for an Earth-to-Mars low-thrust transfer. The singular control is characterized through an intermediate thrust level. We have verified that the optimal maneuver is, however, weakly sensitive to the intermediate thrust value, so this finding is mainly of academic interest in the current computations. 
%It is also interesting to note that hyperbolic tangent smoothing has been able to capture singular thrust profiles. 

%Had it not been for this bifurcation, the last thrust arc would have remained attached to the final boundary and only shrink in its length as the thrust magnitude is increased (similar to the very first thrust arc that never detaches from the initial time). 
%The second interesting phenomenon is, in fact, an immediate consequence of the bifurcation event (at $T \approx 0.2588$ N).

The second interesting phenomenon is that the terminal phase of flight becomes a pure ballistic coast arc. As the maximum thrust magnitude is increased above $T \approx 0.2588$ N, the singular thrust arc takes on an off-bang-off optimal structure; however, beyond another critical thrust magnitude, $T \approx 0.3146$ N, there is no final thrust arc at time $t_f$. This is the thrust magnitude at which the switching function at the final time crosses $S = 0$ line and takes a negative value (thrust off). 
%Note that the width of the thrust ridge decreases for $ 0.2588 < T < 0.3146$ N as the thrust magnitude is increased. In fact, the last point of thrusting is at the point of the rendezvous. 
\begin{figure}[htbp!]
\begin{multicols}{2}
\centering
\includegraphics[width = 0.50\textwidth]{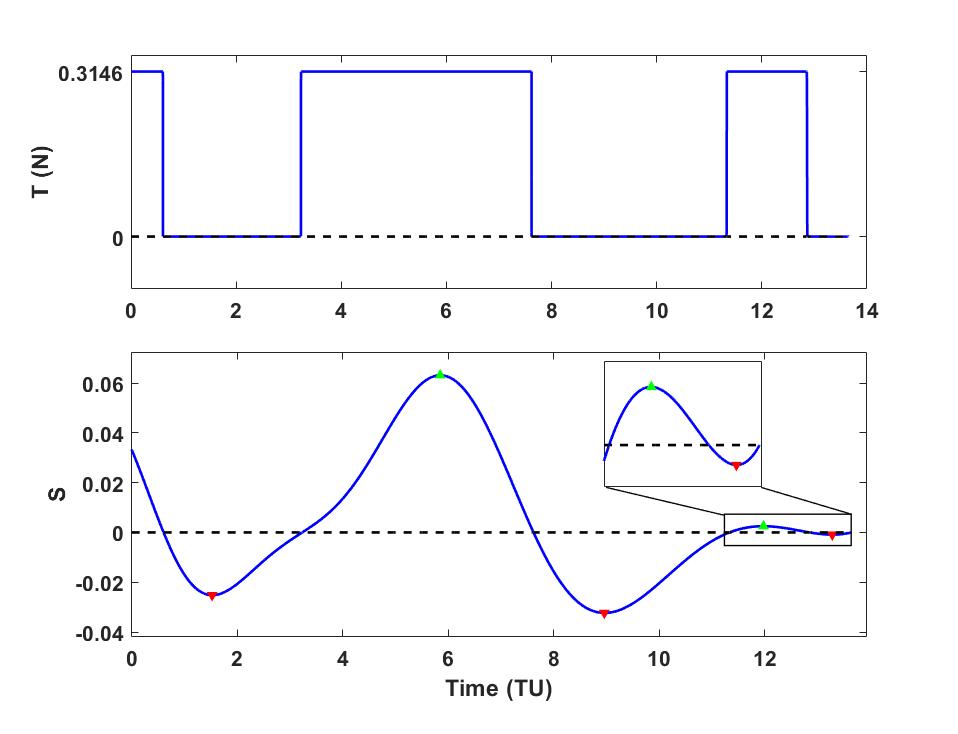}
\caption{Earth-to-Mars switching function with \\$T \approx 0.3146$ N.}
\label{fig:EMCT4_SF}
\centering
\includegraphics[width = 0.48\textwidth]{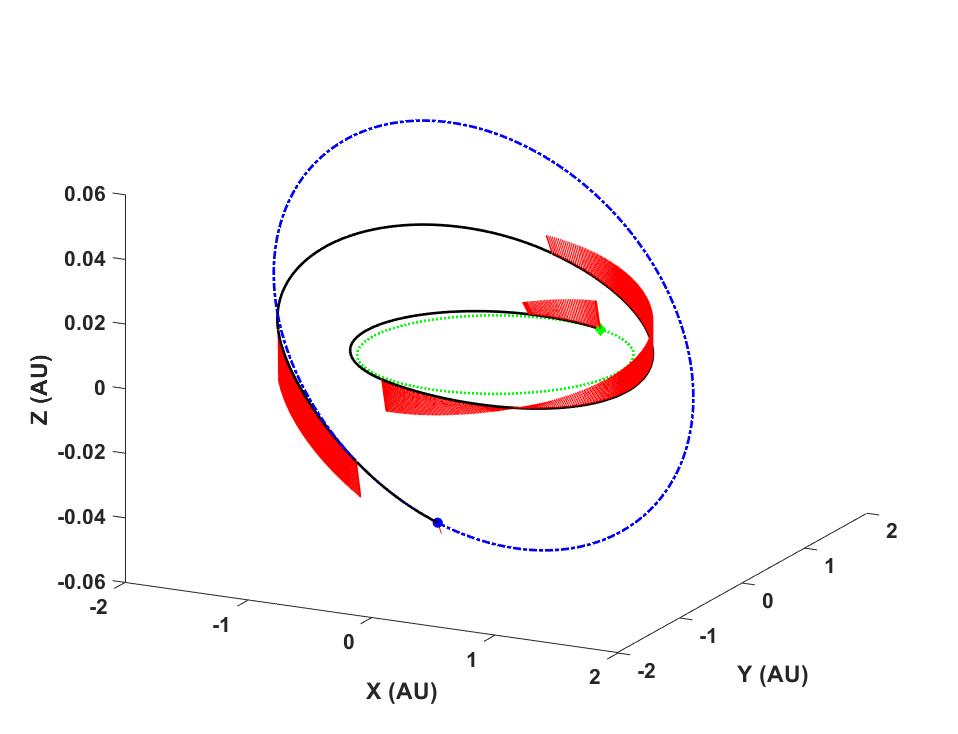}
\caption{Earth-to-Mars heliocentric trajectory with \\$T \approx 0.3146$ N.}
\label{fig:EMCT4_Traj}
\end{multicols}
\end{figure}

In other words, with a propulsion system that creates a higher thrust, it is possible to reduce the time of mission, and rendezvous with Mars at an earlier time than the chosen $t_f$. This time coincides with the right-most point of the last thrust ridge created due to the bifurcation. The locus of these points is coined as ``early-arrival'' boundary. Had we formulated the optimal maneuver with final time free (but confined within a range), this early arrival time would have been found and along this free final time boundary, the Hamiltonian would vanish (since the Hamiltonian is not an explicit function of time in our formulation). Thus, the fixed initial and final time switching surface clearly reveals the boundaries for an infinite set of free initial and final time maneuvers.

A slight increment in the thrust magnitude permits a shift in the time of flight from $t_f$ to a point on the early-arrival boundary corresponding to that particular thrust magnitude. This point denotes the beginning of the early-arrival boundary (see Figure \ref{fig:SSEM_lower}). The physical interpretation of such a trajectory is that the spacecraft rendezvous with the Mars on its orbit at an earlier time and coasts with the Mars for the remainder of time until $t_f$. 

\begin{figure}[htbp!]
\begin{multicols}{2}
\centering
\includegraphics[width = 0.5\textwidth]{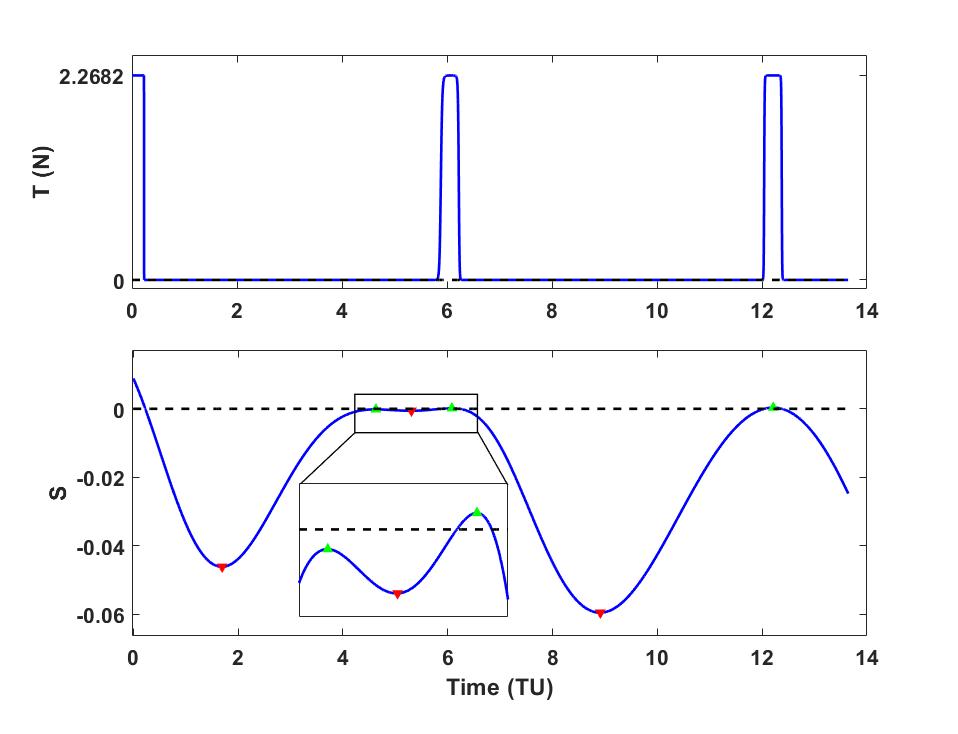}
\caption{Earth-to-Mars switching function with \\$T \approx 2.2682$ N.}
\label{fig:EMCT6_SF}
\centering
\includegraphics[width = 0.5\textwidth]{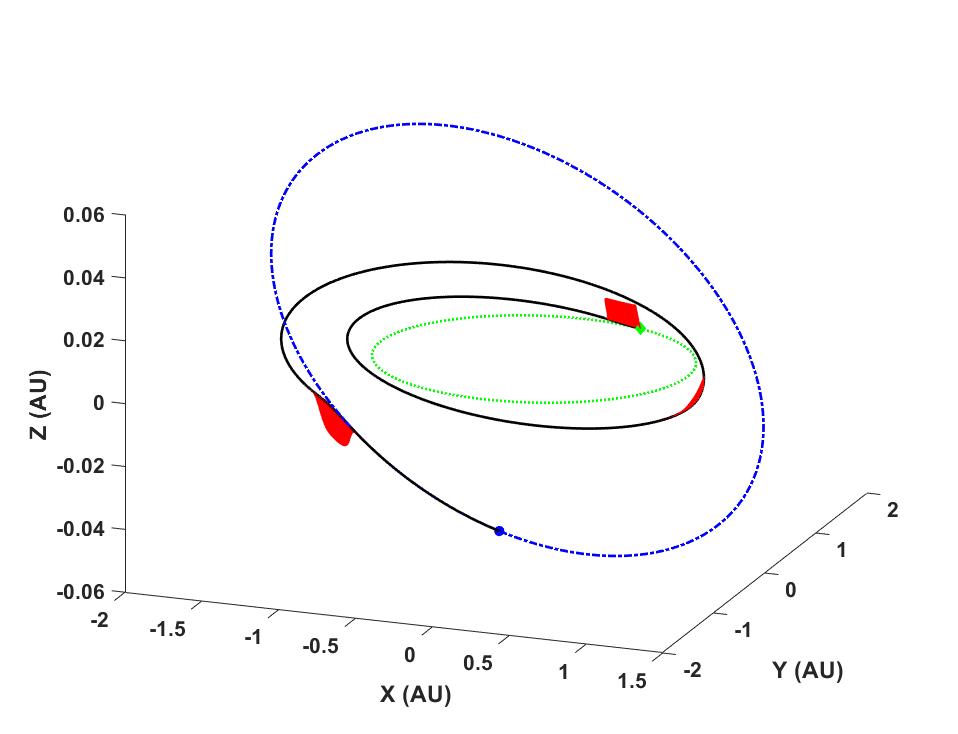}
\caption{Earth-to-Mars heliocentric trajectory with $T \approx 2.2682$ N.}
\label{fig:EMCT6_Traj}
\end{multicols}
\end{figure}

Note that, in this problem, the early-arrival boundary has its own profile. While the thrust magnitude is increased in a monotonic manner, the time of flight (governed by the profile of the early-arrival boundary) exhibits a decrease-increase-decrease profile. Therefore, the shortest time of flight in the lower part of the surface corresponds to the left-most point of this early-arrival boundary (see Figure \ref{fig:SSEM_lower}). A weak local minimum in the time of flight occurs in the vicinity of $T \approx 0.3758$ N. This local minimum in time of flight occurs in the the lower thrust region. For $T > 0.44$ N, the yellow line in Figure \ref{fig:SSEM_lower} shows the boundary monotonically shifts to the left of this local minimum. It is interesting that the analysis of the switching surface leads to transfer-type solutions while the original problem has been formulated as a rendezvous-type maneuver with a fixed time of flight. 
\begin{figure}[htbp!]
\begin{multicols}{2}
\centering
\includegraphics[width = 0.5\textwidth]{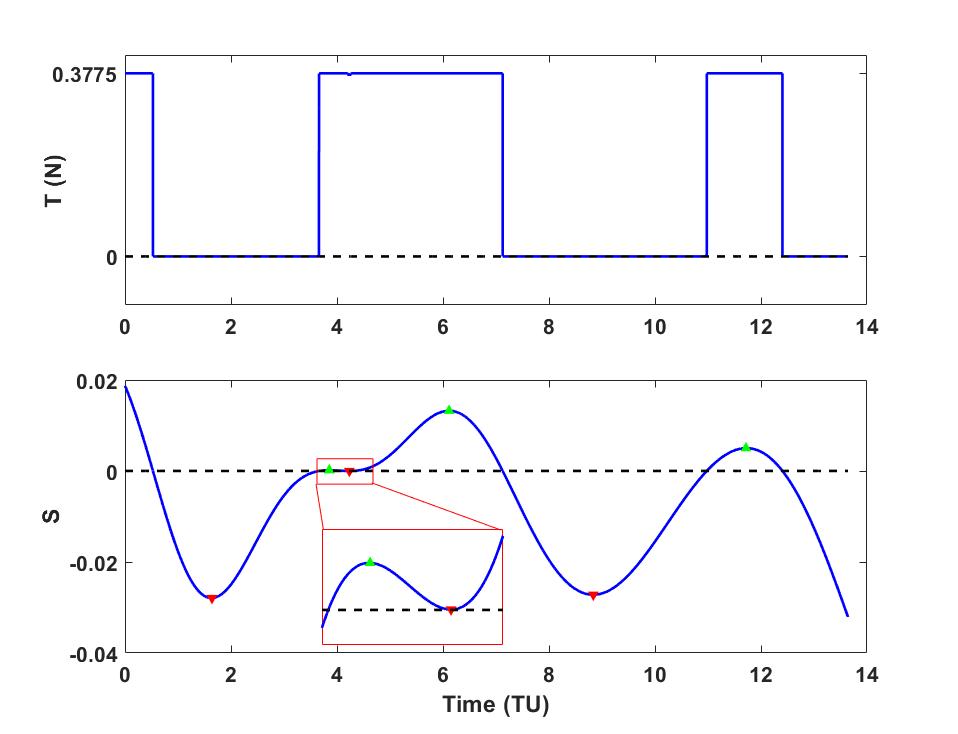}
\caption{Earth-to-Mars switching function with \\$T \approx 0.3775$ N.}
\label{fig:EMCT5_SF}
\centering
\includegraphics[width = 0.5\textwidth]{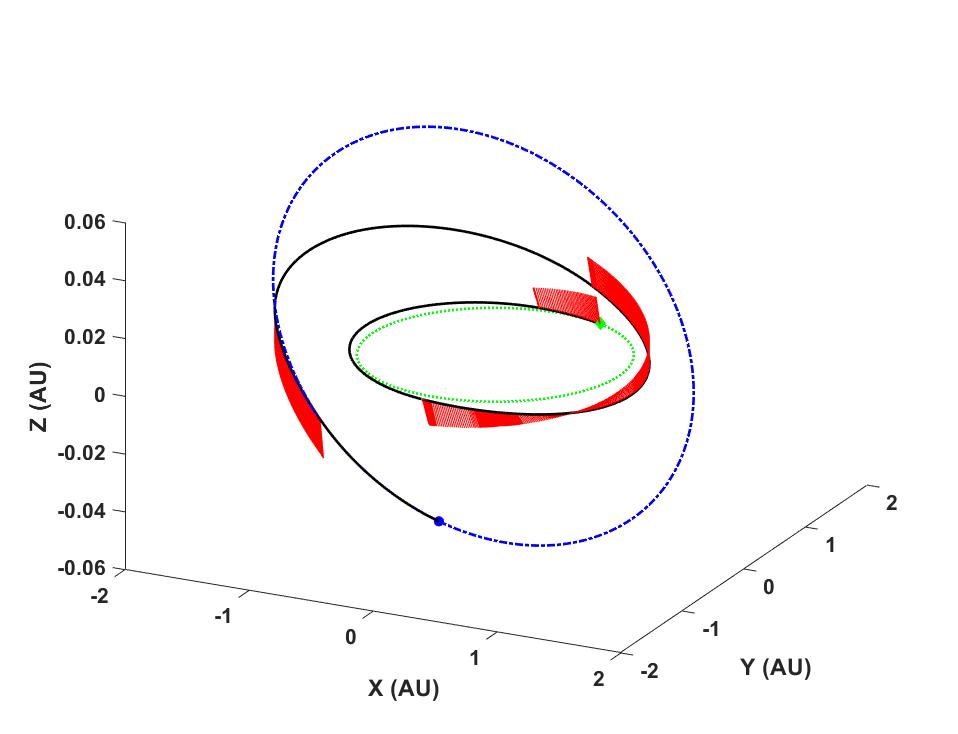}
\caption{Earth-to-Mars heliocentric trajectory with $T \approx 0.3775$ N.}
\label{fig:EMCT5_Traj}
\end{multicols}
\end{figure}
By increasing the thrust magnitude a branching occurs in the middle thrust ridge and divides it into two thrust ridges, i.e., one small thrust ridge is created (branched out) to the left of the main middle thrust ridge. Figures \ref{fig:EMCT5_SF} and \ref{fig:EMCT5_Traj} show the switching function, thrust profile and heliocentric trajectory for this critical thrust, $T \approx 0.3775$ N. At first glance, one might anticipate the possibility of another singular arc near $t = 4.1$ TU, however, a zoom on this region reveals distinct local zeros at the switch function. However, the width of the small (dagger-like) thrust ridge shrinks (to a double zero of $S$, $\dot{S} = 0$ and $\ddot{S} > 0$ at a point), where the thrust ridge eventually disappears at some specific value of thrust, $T \approx 2.2682$ N. Figures \ref{fig:EMCT6_SF} and \ref{fig:EMCT6_Traj} show the switching function, thrust profile and heliocentric trajectory for this critical thrust, $T \approx 2.2682$ N. Another important point worthy of mentioning is that this narrow thrust ridge lasts over a large range of thrust parameters, $T \in \{0.3775,2.2682\}$ N. 

For thrust magnitudes beyond this critical magnitude, application of the traditional control smoothing methods (we tried both logarithmic [\citen{bertrand2002new}] and hyperbolic tangent smoothing) encounter difficulties due to the fact that the length of the thrust ridges becomes so small that smoothing parameter, $\rho_{\text{min}}$ has to take very small values to capture the optimal off-bang-off thrust profile; in fact, the typical continuation procedures become very sensitive to the changes in the continuation parameter. Therefore, the methodology outlined in [\citen{zhu2017solving}] becomes extremely helpful for those very high-thrust regions of the switching surface. Beyond this critical thrust magnitude, $T = 2.2682$ N optimal trajectories consist of only three thrust ridges separated by coast canyons and the width of thrust ridges shrink as the thrust magnitude is increased. 
On the other hand, the number of remaining thrust ridges can be viewed as an early indication of the number of impulses in a pure $N$-impulse solution, which corresponds to when the y axis goes to infinity in such a fashion that the thrust magnitude times delta time remains finite, i.e., $T \times \Delta t \Rightarrow \Delta v$. For finite thrust, the product of the thrust and the time duration is approximately the $\Delta v$ of each impulse (see Eq.~\eqref{eq:approxformula}). Note that from the switching surface, the shortest transfer time coincides with the early arrival moment of the last impulse. Therefore, the impulsive solution, in general, is found to establish the lower limit of the time of flight.

The switching function at a high-trust magnitude ($T \approx 3$ N) is used for generating the impulsive solution; we see $S>0$ only in 3 small time intervals. On the other hand, the minimum-$\Delta v$ 2-impulse solution can be obtained by solving the corresponding Lambert problem. The magnitude of the impulses at the initial and final time instants are $\Delta v(t_0) = 3.0157$ km/s and $\Delta v(t_f) = 3.0318$ km/s, respectively, which correspond to a solution with one revolution around the Sun, i.e., $N_{\text{rev}} = 1$. Table \ref{tab:impulsiveEM} summarizes the initial (unoptimized) solution and the optimized one and the minimum-$\Delta v$ 2-impulse solution obtained by solving the Lambert problem. The middle column summarizes the approximate values obtained from Eq.~\eqref{eq:approxformula}, which denote the initial values before performing $N$-impulse optimization. 
\begin{table}[h!] 
	\begin{center} 
		\caption{Summary of minimum-$\Delta v$ 2- and 3-impulse solutions for the Earth-to-Mars problem.}\label{tab:impulsiveEM}
		{\small%\scriptsize
		\begin{tabular}{c c c c c c c}
        \hline
        \hline
         \multirow{2}{*}{impulse \#} &    \multicolumn{2}{c}{\textbf{Lambert}}&       \multicolumn{2}{c}{\textbf{3-impulse}}&     \multicolumn{2}{c}{\textbf{Optimized 3-impulse}}\\
                    &     $\textbf{t}_{i}$ \textbf{(days)} &    $\bm{\Delta} v$ \textbf{(km/s)}  &        $\textbf{t}_{i}$ \textbf{(days)} &    $\bm{\Delta} v$ \textbf{(km/s)} &   $\textbf{t}_{i}$ \textbf{(days)} &    $\bm{\Delta}v $ \textbf{(km/s)} \\
                    \cline{2-7}
         %\hline
         1       &       0.0       &     3.0157            &         0           &       1.356          &       0       &        1.417 \\
         2       &       793       &     3.0318            &         354.27      &       2.029          &     358.99    &        1.925 \\
         3       &        -        &       -               &         710.78      &       2.168          &     711.72    &        2.268 \\
         \hline
         $\sum \Delta v$        &    -    &     6.047      &           -         &       5.554          &              &        \textbf{5.611}\\
        \hline
        \end{tabular}
		}
	\end{center}
\end{table}
\begin{figure}[htbp!]
\centering
\includegraphics[width=5.0in]{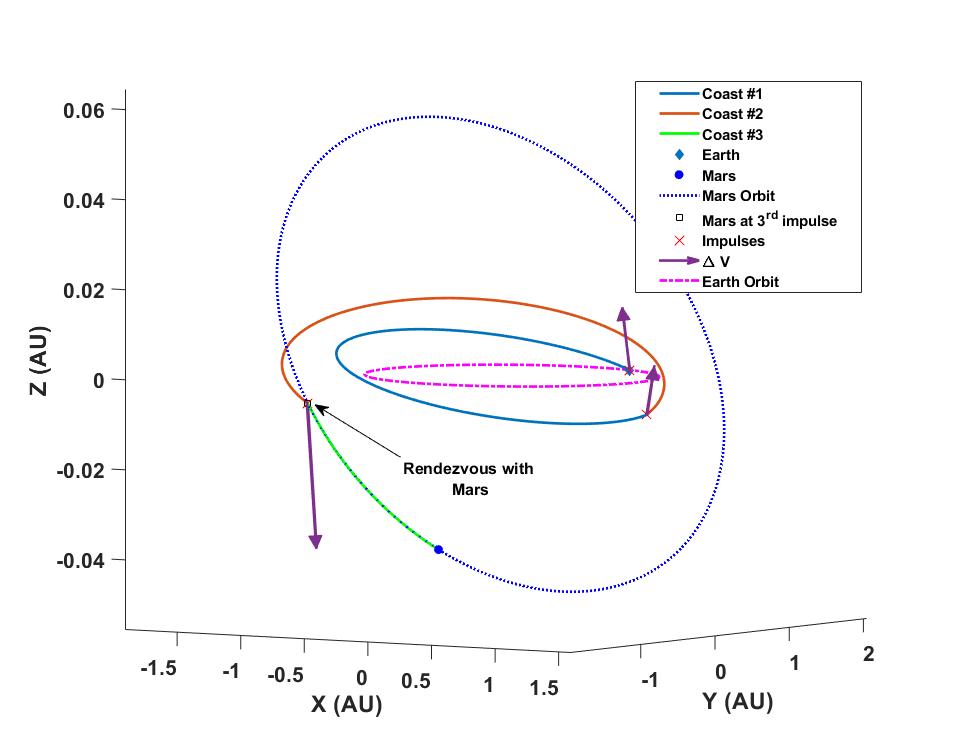}
\caption{Minimum-$\Delta v$ 3-impulse trajectory for the Earth-to-Mars problem.}
\label{fig:EM_3ImpulseTraj}
\end{figure}
Figure \ref{fig:EM_3ImpulseTraj} shows the details of the minimum-$\Delta v$ 3-impulse trajectory for the Earth-to-Mars problem. In all trajectory plots, initial and final locations correspond to the positions at the departure date, $t_0$ and the arrival date, $t_f$. The green arc denotes the last coast arc while on Mars orbit if the time of flight has been set to $t_f = 793$ days. The optimal number of impulses for this problem is not large, but it serves as the prelude to study problems with a greater number of impulses. 

The numbers indicate that the optimal minimum-$\Delta v$ 3-impulse solution corresponds to a 7.2\% reduction in total $\Delta v$ compared to the minimum-$\Delta v$ 2-impulse solution obtained by solving Lambert problem. In addition, time of flight of the 3-impulse trajectory, $t_f = 711.72$ days shows a 10.24\% reduction compared to the time of flight of minimum-$\Delta v$ 2-impulse solution, $t_f = 793$ days. Note that the spacecraft rendezvous with Mars at the time instant of the third impulse (about 81.28 days earlier, i.e., nearly three months).

One may ask the following question: is it possible to recover the original 2-impulse solution by increasing thrust magnitude to infinity, $T = \infty$? The answer lies in the fact that optimality conditions have been guaranteed to establish that three thrust impulses are required to minimize propellant consumption. In other words, it is not possible to recover 2-impulse minimum-$\Delta v$ solution (where the impulses are at prescribed initial and final times) since the 2-impulse solution is obviously not \textit{optimal}! This is a known fact and is used for designing many-impulse solutions. Put it another way, the time history of the magnitude of the primer vector for the initial bi-impulsive solution violates the optimality conditions of Lawden; hence, it is possible, in principle, to improve the solution [\citen{handelsman1968primer}].
\subsection{Interplanetary Rendezvous From Earth to Asteroid 1989ML}
\begin{figure}[htbp!]
\begin{multicols}{2}
\centering
\includegraphics[width=0.50\textwidth]{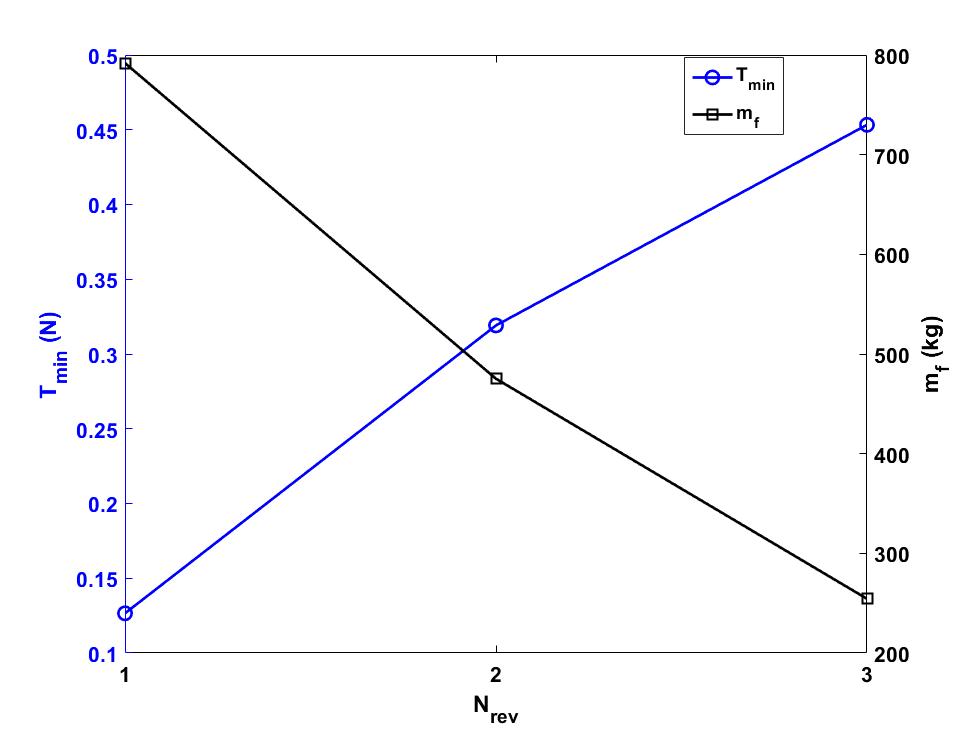}
\caption{Changes in $T_{\text{min}}$ and $m_f$ vs. $N_{\text{rev}}$ for the Earth-to-1989ML problem.}
\label{fig:ETo1989MLMinThrustNrevs}
\hfill
\centering
\includegraphics[width=0.5\textwidth]{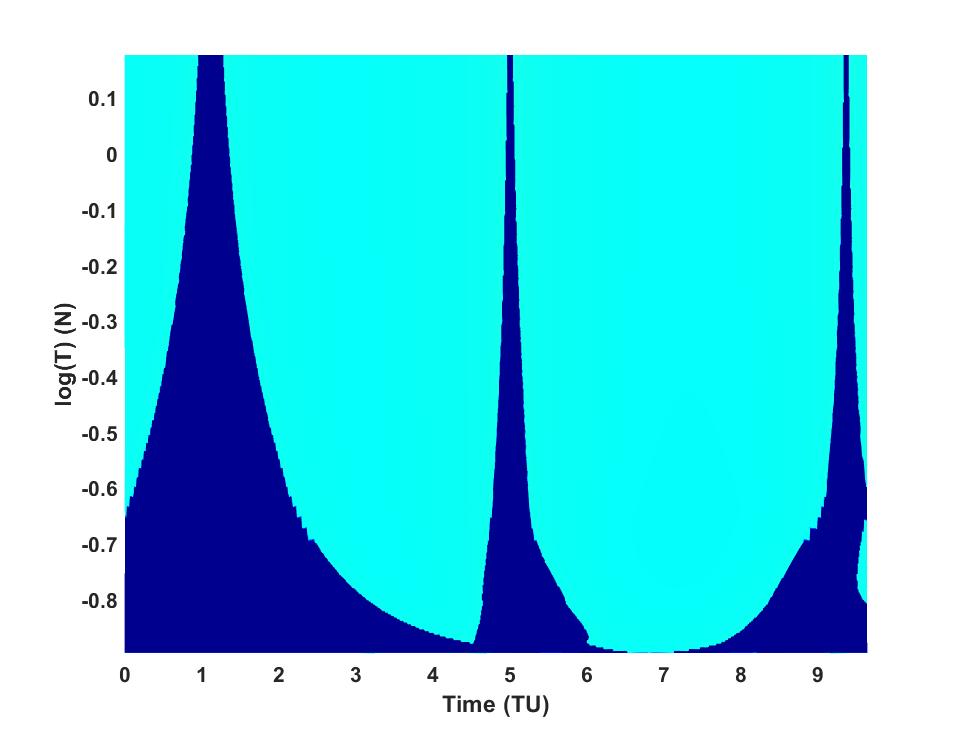}
\caption{Switching surface $S = 0$ contour map for the Earth-to-19989ML problem for $N^*_{\text{rev}} = 1$.}
\label{fig:ETo1989MLSS}
\end{multicols}
\end{figure}

This problem is taken from [\citen{taheri2016initial}] and Table.~\ref{tab:OE_1898ML} gives the classical orbital elements of the asteroid 1989ML in which the epoch date is given as the Modified Julian Date (MJD).
\begin{table}[h!] 
\begin{center} 
		\caption{Keplerian orbital elements of the asteroid 1989ML wrt the Sun.}\label{tab:OE_1898ML}
		{\small%\scriptsize
		\begin{tabular}{c c c c c c c}
        \hline
        \hline
         $a$ & $e$ & $i$& $\Omega$& $\omega$& $M$& Epoch \\
         $[\text{AU}]$ &  & [deg]& [deg]& [deg]& [deg]& [MJD] \\
         \hline
         1.2721 & 0.13649 & 4.3782& 104.3571& 183.3249& 117.36689& 53900 \\
        \hline
        \hline
        \end{tabular}
		}
	\end{center}
\end{table}
The eccentricity and inclination of this asteroid make it a moderately difficult-to-reach target (see Table \ref{tab:OE_1898ML}), due to the $\approx 4.4^{\circ}$ inclination. The following values are considered for the parameters of the spacecraft and its low-thrust propulsion system: $m_0 = 1000$ kg,  and $I_{\text{sp}} = 3000$ s. The specified time of flight is $t_f-t_0 = 560$ days. The Earth position and velocity vectors at the departure time, $t_0$ are $\textbf{r}_{\oplus} = [-109310123.96,-103935506.96, 1736.32]^{\top}$ km, $\textbf{v}_{\oplus} =  [20.0414,-21.7003,0.000309]^{\top}$ km/s, respectively. The position and velocity vectors of the asteroid 1989ML at the final time, $t_f$ are given as $\textbf{r}_{\text{T}} = [81709931.65, -143042471.97, -3344947.036]^{\top}$ km, $\textbf{v}_{\text{T}} = [26.5207,14.3234,-2.2390]^{\top}$ km/s, respectively.
\begin{figure}[htbp!]
\begin{multicols}{2}
\centering
\includegraphics[width=0.5\textwidth]{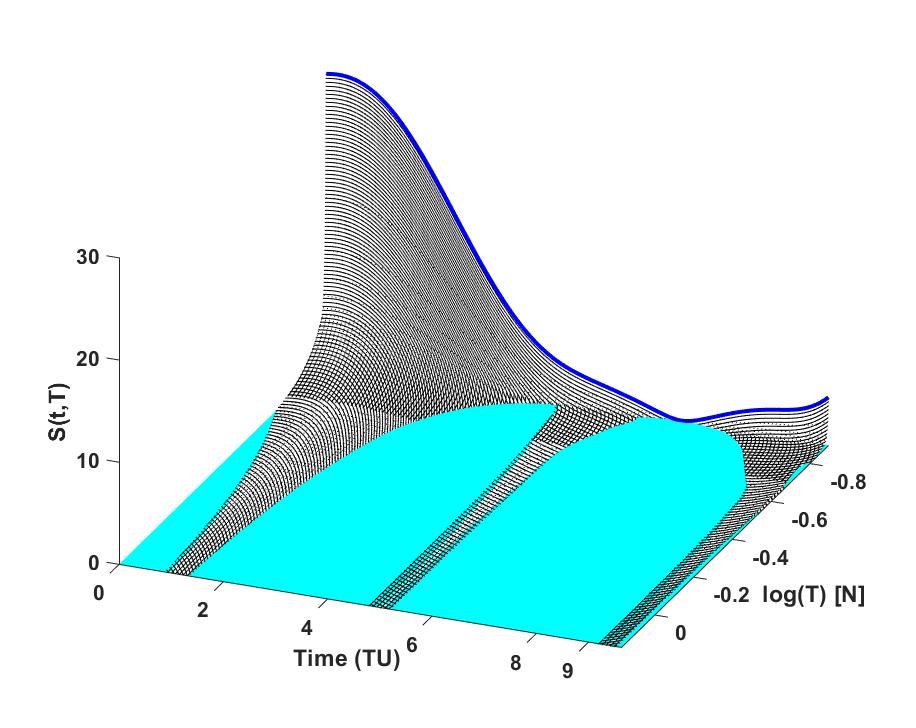}
\caption{Switching surface for the Earth-to-1989ML problem with $N^*_{\text{rev}} = 1$.}
\label{fig:ETo1989ML_SS_log_2}
\hfill
\centering
\includegraphics[width=0.50\textwidth]{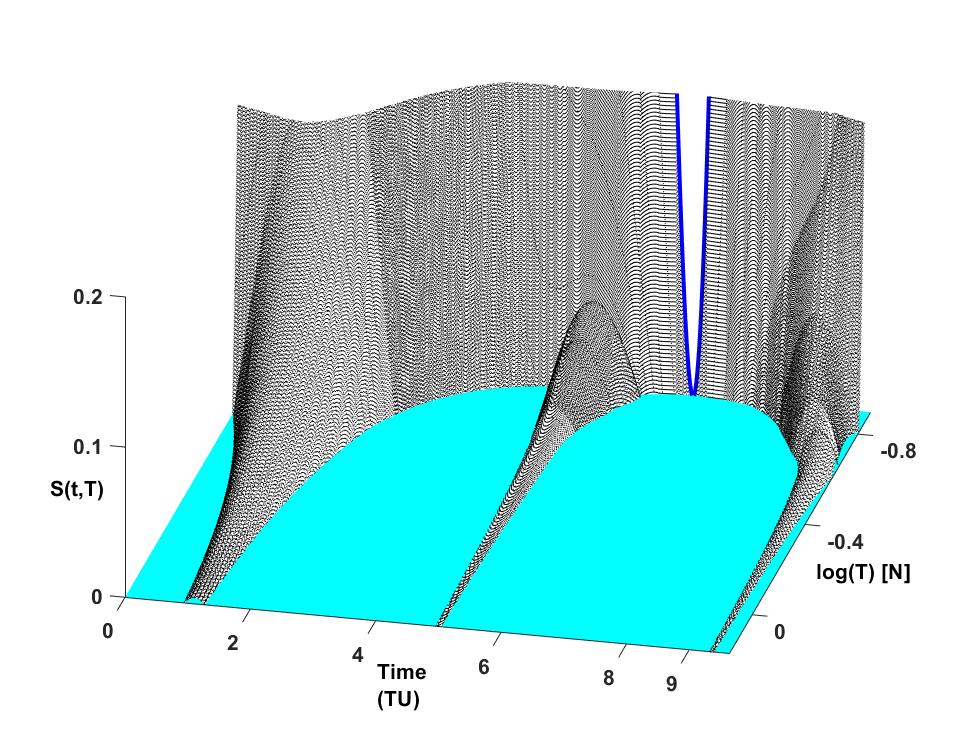}
\caption{Enlarged view of the switching surface for the Earth-to-1989ML problem with $N^*_{\text{rev}} = 1$.}
\label{fig:ETo1989ML_SS_log_3}
\end{multicols}
\end{figure}
The value of thrust is swept over $T \in [T_{\text{min}},T_{\text{max}}]$, where $T_{\text{max}} = 1.5$ N. 
Figure \ref{fig:ETo1989MLMinThrustNrevs} shows the changes in $T_{\text{min}}$ and final mass $m_f$ versus three different feasible values of $N_{\text{rev}}$; there is no feasible minimum-fuel solution for $N_{\text{rev}}=0$ due to the adverse initial phasing, which requires more-than-available propellant. The fundamental minimum-thrust solution corresponds to $N^*_{\text{rev}} = 1$. The fundamental minimum-thrust magnitude for the given BCs is $T_{\text{min}} = 0.12659$ N. The thrust value is then swept in the given range $T \in [0.1266, 1.5]$ N. 

\begin{figure}[htbp!]
\centering
\includegraphics[width=4.0in]{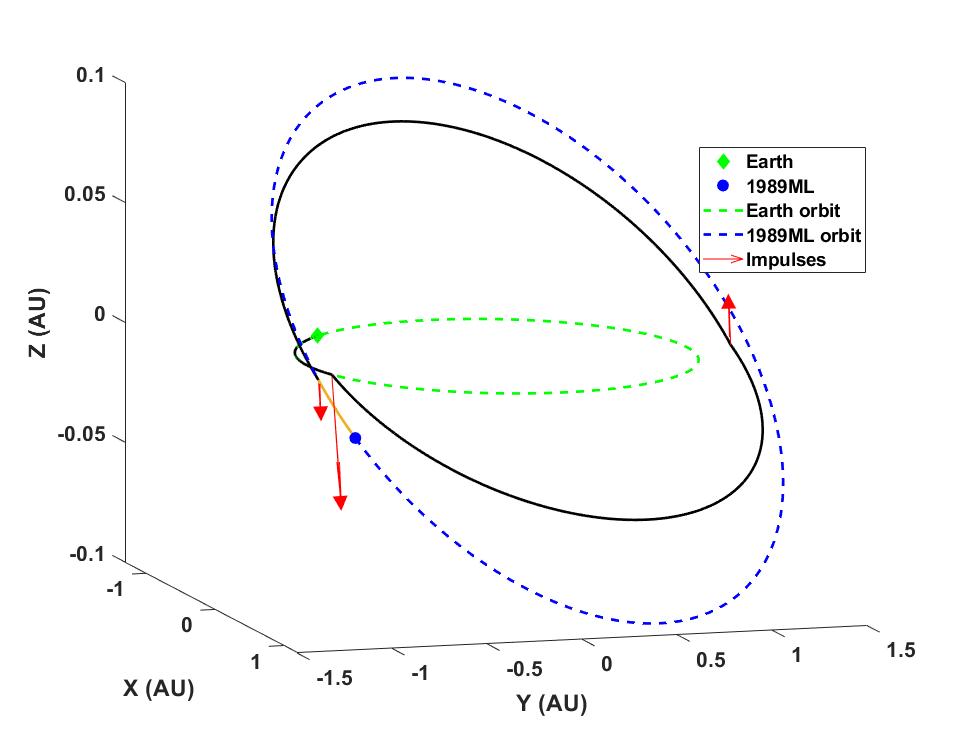}
\caption{Minimum-$\Delta v$ 3-impulse trajectory for the Earth-to-1989ML problem.}
\label{fig:ETo1989ML_3ImpulseTraj}
\end{figure}

Figure \ref{fig:ETo1989ML_SS_log_2} shows the fundamental switching surface for the Earth-to-1989ML problem and Figure \ref{fig:ETo1989ML_SS_log_3} shows an enlarged view of the surface for $S<0.2$. Figure \ref{fig:ETo1989MLSS} shows the switching surface $S=0$ contour map where it is easy to distinguish late-departure and early-arrival boundaries. In addition, the early-arrival boundary consists of two segments. The first part lies in the low thrust region. Then, the final coast phase is replaced by a final thrust ridge. Then, with a further increase in the thrust magnitude, the final phase of flight becomes a pure coast and the duration of the coast phase increases as the thrust magnitude is increased. 
Clearly, to establish a rendezvous with the target body, greater thrust and longer time is needed which corresponds to the rightmost point on the last thrust ridge as it remains on the $t = t_f$ boundary. The peculiar profile of the early-arrival boundary is a indication of the difficulty encountered in reachability analysis in astrodynamics. Eventually, only three thrust ridges remain and become increasingly narrow for increasing thrust. This indicates that the optimal number of impulses is three at the high-thrust limit. Table \ref{tab:impulsiveETo1989ML} summarizes the impulsive solution and Figure \ref{fig:ETo1989ML_3ImpulseTraj} shows the details of the minimum-$\Delta v$ 3-impulse trajectory for the Earth-to-1989ML problem.
\begin{table}[h!] 
	\begin{center} 
		\caption{Summary of minimum-$\Delta v$ impulsive solutions for the Earth-to-1989ML problem.}\label{tab:impulsiveETo1989ML}
		{\small%\scriptsize
		\begin{tabular}{c c c c c c c}
        \hline
        \hline
         \multirow{2}{*}{impulse \#} &    \multicolumn{2}{c}{\textbf{Lambert}}&       \multicolumn{2}{c}{\textbf{3-impulse}}&     \multicolumn{2}{c}{\textbf{Optimized 3-impulse}}\\
                    &     $\textbf{t}_{i}$ \textbf{(days)} &    $\bm{\Delta} v$ \textbf{(km/s)}  &        $\textbf{t}_{i}$ \textbf{(days)} &    $\bm{\Delta} v$ \textbf{(km/s)} &   $\textbf{t}_{i}$ \textbf{(days)} &    $\bm{\Delta} v$ \textbf{(km/s)} \\
                    \cline{2-7}
         %\hline
         1       &       0.0       &     2.7891            &         64.465       &       2.6033          &     64.4932    &        2.5999 \\
         2       &       560       &     4.08984           &         290.347      &       0.6551          &     290.347    &        0.7082 \\
         3       &        -        &       -               &         544.185      &       0.5786          &     544.272    &        0.61077 \\
         \hline
         $\sum \Delta v$        &         &      6.879             &                    &       3.837          &              &        \textbf{3.9189}\\
        \hline
        \end{tabular}
		}
	\end{center}
\end{table}
% \begin{figure}[htbp!]
% \centering
% \includegraphics[width=5.0in]{figures/ETo1989ML_SS.eps}
% \caption{Switching surface for the Earth-to-19989ML problem.}
% \label{fig:ETo1989MLSS}
% \end{figure}
% \begin{table}[h!] 
% 	\begin{center} 
% 		\caption{Summary of the 3-impulse solution for the Earth-to-1989ML problem.}\label{tab:impulsiveETo1989ML}
% 		{\small%\scriptsize
% 		\begin{tabular}{ c c c c c}
%         \hline
%         \hline
%          impulse \# &  1 & 2 & 3 & $\sum_1^3 \Delta v_i$\\
%           \hline                         % &    & \cline{3-8}\\
%          $t_i$ (days)        & 64.4932 & 290.347 & 544.272 & - \\
%          $\Delta v_i$ (km/s) & 2.5999  & 0.70820 & 0.61077 &  \textbf{3.9189} \\
%          \hline
%          \hline
%         \end{tabular}
% 		}
% 	\end{center}
% \end{table}

\subsection{Interplanetary From Earth and Rendezvous With Asteroid Dionysus}
Here, an interplanetary minimum-fuel mission from the Earth to asteroid Dionysus is studied, where the optimal solution is known from [\citen{taheri2016enhanced}]. This asteroid is a fairly difficult and expensive target to reach due to the orbit's high eccentricity and inclination values of 0.542 and $13.54$ degrees, respectively. The BCs and parameters are chosen to match to those reported in [\citen{taheri2016enhanced}]. The solution to this problem consists of multiple revolutions around the Sun with intermediate thrust and coast arcs, so we can expect this family of orbit transfers to result in a very different switching surface than the one obtained in the Earth-to-Mars and Earth-to-1989ML transfer cases. The following values are considered for the parameters of the spacecraft and its low-thrust propulsion system: $m_0 = 4000$ kg, and $I_{\text{sp}} = 3000$ s. 
\begin{figure}[htbp!]
\centering
\includegraphics[width=5.5in]{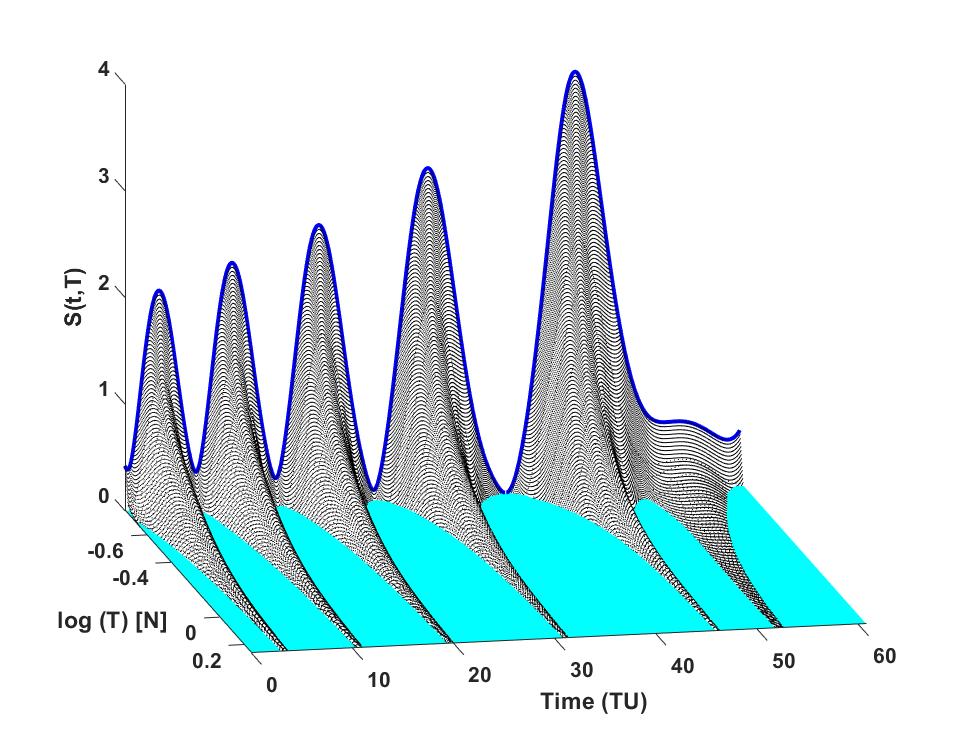}
\caption{Switching surface for the Earth-to-Dionysus problem with $N^*_{\text{rev}} = 5$.}
\label{fig:EDSS_3}
\end{figure}
Table.~\ref{tab:OE} gives the classical orbital elements of the asteroid Dionysus in which the Epoch date is given as the Modified Julian Date (MJD).
\begin{table}[h!] 
	\begin{center} 
		\caption{Keplerian orbital elements of asteroid Dionysus wrt the Sun.}\label{tab:OE}
		{\small%\scriptsize
		\begin{tabular}{c c c c c c c}
        \hline
        \hline
         $a$ & $e$ & $i$& $\Omega$& $\omega$& $M$& Epoch \\
         (AU) &  & [deg]& [deg]& [deg]& [deg]& [MJD] \\
         \hline
         2.2 & 0.542 & 13.6& 82.2& 204.2& 114.4232& 53400 \\
        \hline
        \end{tabular}
		}
	\end{center}
\end{table}
The spacecraft departs from the Earth on December 23, 2012 and the mission time of flight is 3534 days. So we have a relatively large spacecraft and we wish to use low thrust - it is not surprising that the time of flight will be thousands of days. Any low-thrust trajectory from the Earth to asteroid Dionysus must undergo large changes in eccentricity and inclination values and the optimal solutions with low thrust un-surprisingly requires several revolutions around Sun. The Earth position and velocity vectors at the departure are $\textbf{r}_{\oplus} = [-3637871.081, 147099798.784, -2261.441]^{\top}$ km, $\textbf{v}_{\oplus} = [-30.265097, -0.8486854, 0.0000505]^{\top}$ km/s, respectively.

\begin{figure}[htbp!]
\begin{multicols}{2}
\centering
\includegraphics[width=0.50\textwidth]{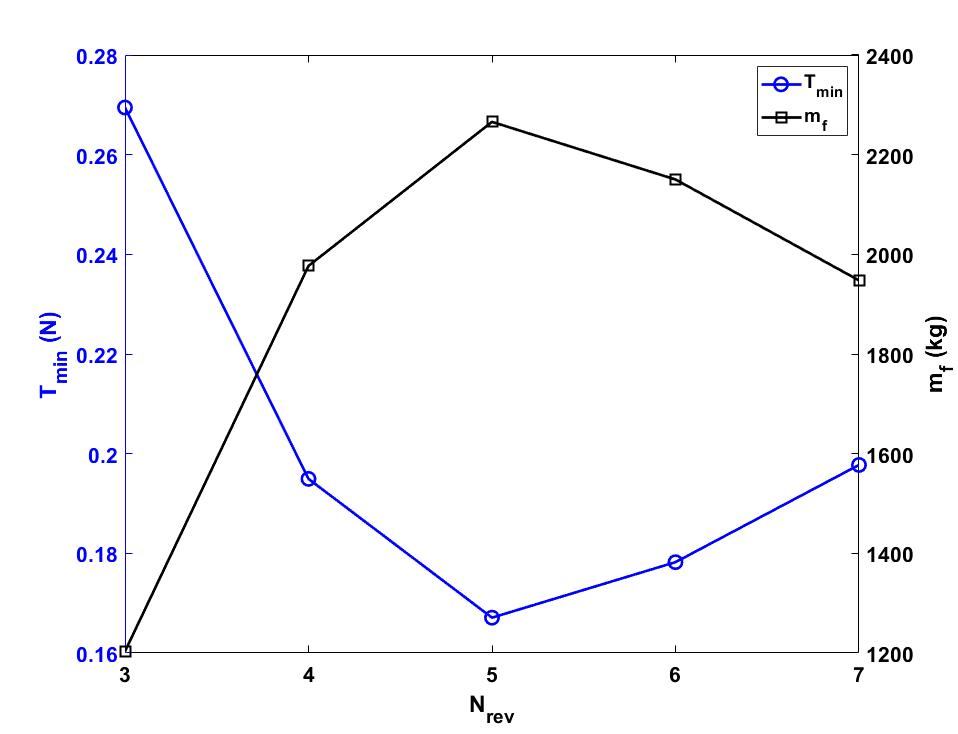}
\caption{Changes in $T_{\text{min}} $ and $m_f$ vs. $N_{\text{rev}}$ for the Earth-to-Dionysus problem; $m_0 = 4000$ kg.}
\label{fig:EDMinThrustNrevs}
\hfill
\centering
\includegraphics[width=0.5\textwidth]{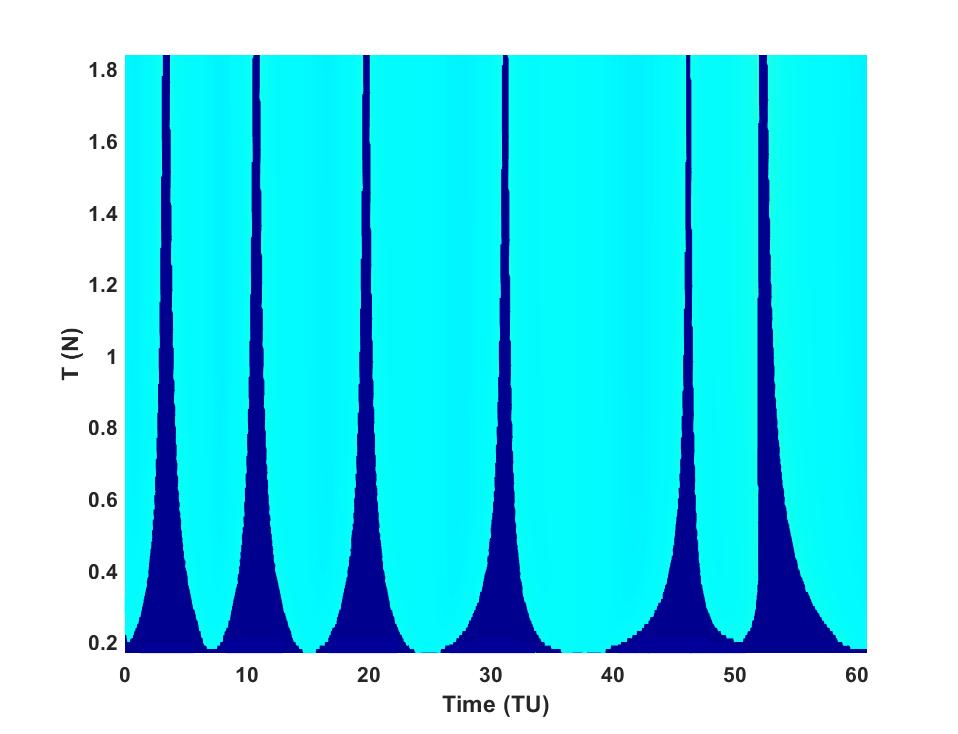}
\caption{Switching surface $S = 0$ contour map for the Earth-to-Dionysus problem for $N^*_{\text{rev}} = 5$.}
\label{fig:EDSS}
\end{multicols}
\end{figure}

\begin{figure}[htbp!]
\begin{multicols}{2}
\centering
\includegraphics[width=0.52\textwidth]{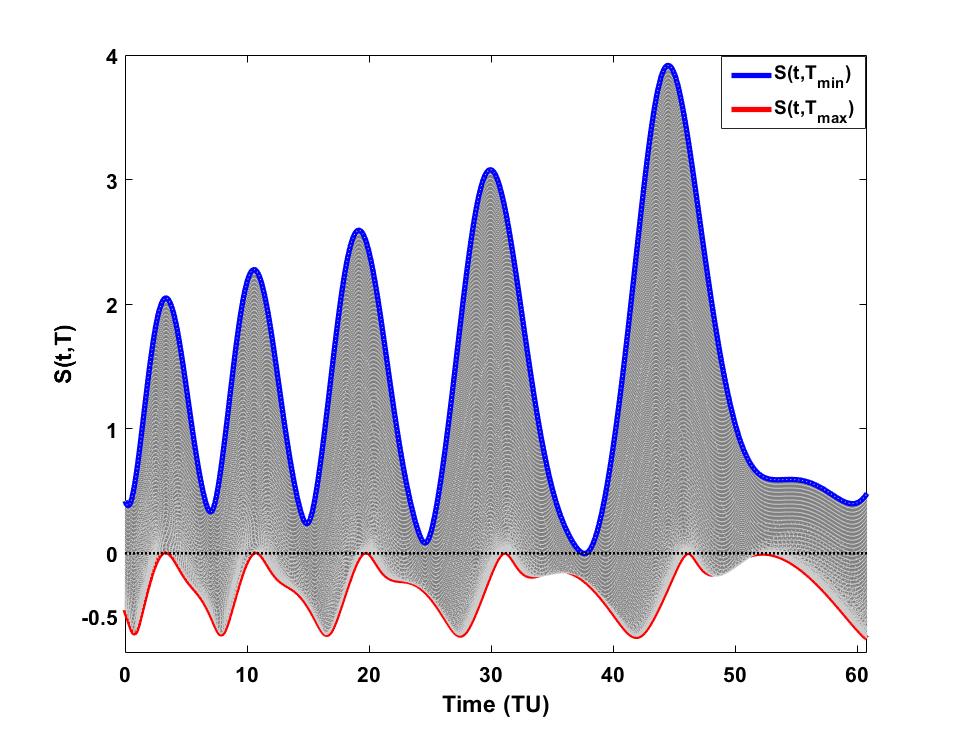}
\caption{Envelope of the switching functions for the Earth-to-Dionysus problem with $T_{\text{min}} = 0.1671$ N and $T_{\text{max}} = 1.8$ N with $N^*_{\text{rev}} = 5$.}
\label{fig:EDSS_Side}
\hfill
\centering
\includegraphics[width=0.50\textwidth]{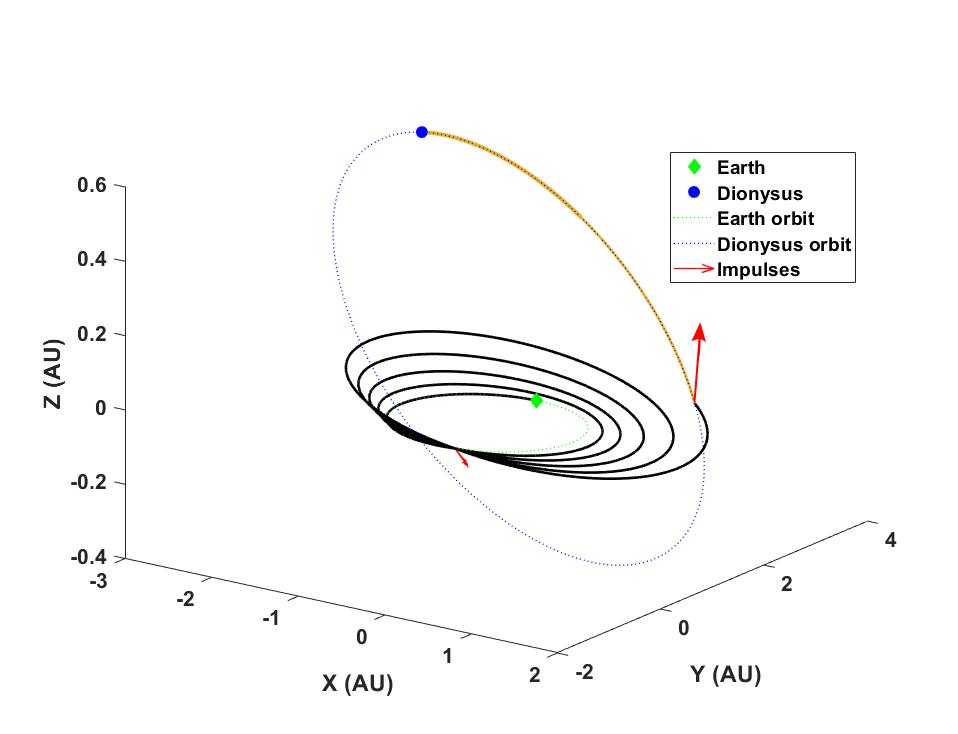}
\caption{Minimum-$\Delta v$ 6-impulse trajectory for the Earth-to-Dionysus problem with $N^*_{\text{rev}} = 5$.}
\label{fig:ED_6ImpulseTraj}
\end{multicols}
\end{figure}

Figure \ref{fig:EDSS_3} shows the fundamental switching surface for the Earth-to-Dionysus problem over the given range of the thrust magnitudes. Unlike the previous two test cases, Figure \ref{fig:EDMinThrustNrevs} shows that the fundamental minimum-thrust solution, $T_{\text{min}} = 0.1671$ N corresponds to an intermediate number of revolutions, $N^*_{\text{rev}} = 5$ when $T_{\text{min}} = 0.1673$ N. Figure \ref{fig:EDSS} shows the switching surface $S = 0$ contour map for the Earth-to-Dionysus problem. The problem consists of six thrust ridges and the switching surface reveals the existence of late-departure and early-arrival boundaries. Figure \ref{fig:EDSS_Side} shows an envelope of the switching functions for the Earth-to-Dionysus problem. 

%The value of thrust is considered as the sweeping parameter, i.e., $T \in [T_{\text{min}},T_{\text{u}}]$, where $T_{\text{u}}$ is some upper value of thrust and is set to $T_{\text{max}} = 1.8$ N. 
Figure \ref{fig:ED_6ImpulseTraj} shows the details of the minimum-$\Delta v$ 6-impulse trajectory for the Earth-to-Dionysus problem. The first five impulses are applied at the perihelion of the intermediate osculating elliptical orbits, which coincides with the line of nodes of the orbits of the Earth and asteroid Dionysus. The final impulse is also applied at the osculating ascending node and changes the inclination and energy of the spacecraft at the time of (early) rendezvous. The spacecraft arrives early and coasts with the asteroid on its orbit for nearly 501.808 days! Had we left $t_f$ in the optimal control formulation, we would have found that the free final time transversality condition (Hamiltonian = 0) would have converged to the free final time on the early-arrival boundary. A different color is used to denote the final coast phase. This shows that the time of flight can be reduced significantly, in this case, when impulsive maneuvers are performed.
% \begin{figure}[htbp!]
% \centering
% \includegraphics[width=4.0in]{figures/EDSS_Side.eps}
% \caption{Envelope of the switching functions for the Earth-to-Dionysus problem with $T_{\text{min}} = 0.1673$ N and $T_{\text{max}} = 1.8$ N.}
% \label{fig:EDSS_Side}
% \end{figure}

% \begin{figure}[htbp!]
% \centering
% \includegraphics[width=4.0in]{figures/ED_6ImpulseTraj.eps}
% \caption{Minimum-$\Delta v$ 6-impulse trajectory for the Earth-to-Dionysus problem.}
% \label{fig:ED_6ImpulseTraj}
% \end{figure}
Table \ref{tab:impulsiveED} summaries the time and magnitude of the six impulses. Minimum-$\Delta v$ Lambert solutions require significant amount of propellant and are not reported since they are 2-impulse solutions. This is a challenging problem for traditional methods used for generating impulsive trajectories since the number of impulses and the dimension of search space gets large. 
\begin{table}[h!] 
	\begin{center} 
		\caption{Summary of the minimum-$\Delta v$ 6-impulse solution for the Earth-to-Dionysus problem for $N^*_{\text{rev}} = 5$.}\label{tab:impulsiveED}
		{\small%\scriptsize
		\begin{tabular}{ c c c c c c c c}
        \hline
        \hline
         impulse \# &  1 & 2 & 3 & 4 & 5 & 6 & $\sum_1^6 \Delta v_i$\\
          \hline                         % &    & \cline{3-8}\\
         $t_i$ (days) & 193.246 &   624.164 &  1147.393 & 1810.202 & 2683.730 & 3032.192 &  -\\
         $\Delta v_i$ (km/s) & 1.61045 &   1.58896 &  1.594121 & 1.496731 & 1.231273 & 2.385884 &  \textbf{9.90742} \\
         \hline
         \hline
        \end{tabular}
		}
	\end{center}
\end{table}
\subsection{GTO-to-GEO Test Problem}
\begin{figure}[htbp!]
\centering
\includegraphics[width=4.0in]{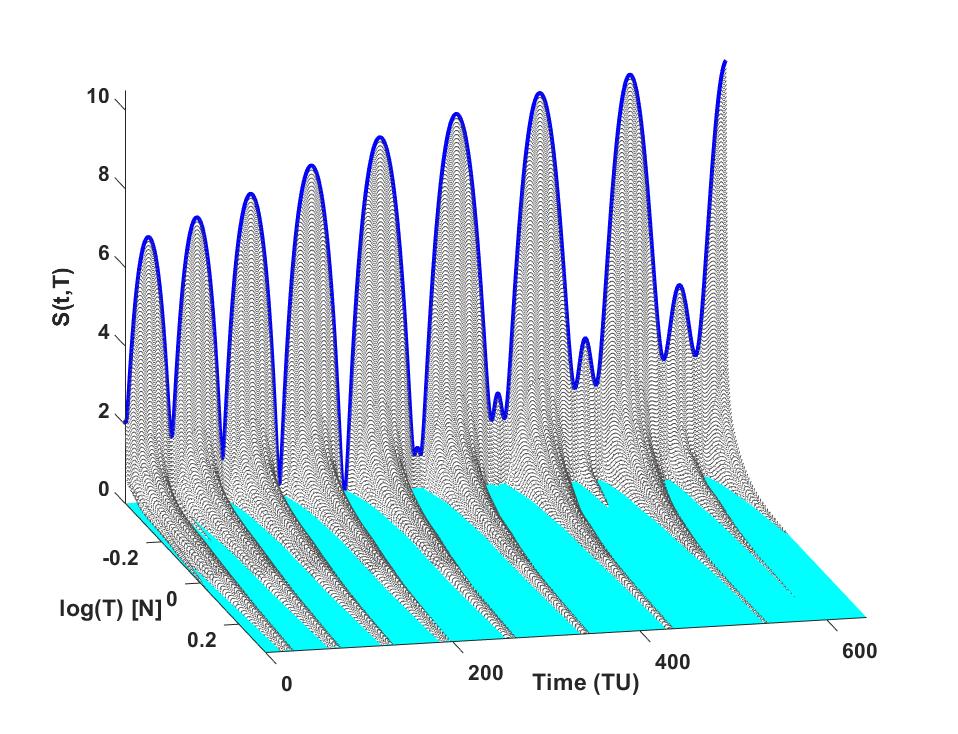}
\caption{Switching surface for the GTO-to-GEO problem with $N^*_{\text{rev}} = 8$.}
\label{fig:GTO2GEOSS_3}
\end{figure}

The fourth test problem that we considered is an orbit raising problem from a geostationary transfer orbit (GTO) to a geostationary Earth orbit (GEO) with their parameters defined in Table \ref{tab:OE-GTO}. Sending a spacecraft to GEO is of practical use since its orbital period is identical to the Earth's rotation period, which makes GEO an ideal orbit for placing communication satellites. It is assumed that the true anomaly of the spacecraft at the departure on the GTO is zero so that the initial position vector is aligned along the positive x-axis of the Earth-centered inertial equatorial frame. In addition, the target point lies along the negative x-axis, i.e., the x-y plane phase angle between the initial and final position vectors is 180 degrees. The initial and final position and velocity vectors are $\textbf{r}_{0} = [6738.9, 0.0, 0.0]^{\top}$ km, $\textbf{v}_{0} = [0.0, 10.0258, 1.231]^{\top}$ km/s, $\textbf{r}_{\textbf{T}} = [-42165, 0.0, 0.0]^{\top}$ km, $\textbf{v}_{\textbf{T}} = [0.0, -3.0746, 0.0]^{\top}$ km/s. 
\begin{table}[h!] 
	\begin{center} 
		\caption{Keplerian orbital elements of the LEO and GTO orbits wrt the Earth.}\label{tab:OE-GTO}
		{\small%\scriptsize
		\begin{tabular}{c c c c c c}
        \hline
        \hline
        Orbit & $a$  & $e$ & $i$   & $\Omega$ & $\omega$ \\
              & (km) &     & [deg] & [deg]    & [deg]    \\
         \hline
        GTO &24505 & 0.725 & 7.0 & 0.0 & 0.0 \\
        GEO &42165 & 0.0   & 0.0 & - & - \\
        \hline
        \end{tabular}
		}
	\end{center}
\end{table}
\begin{figure}[htbp!]
\begin{multicols}{2}
\centering
\includegraphics[width=0.51\textwidth]{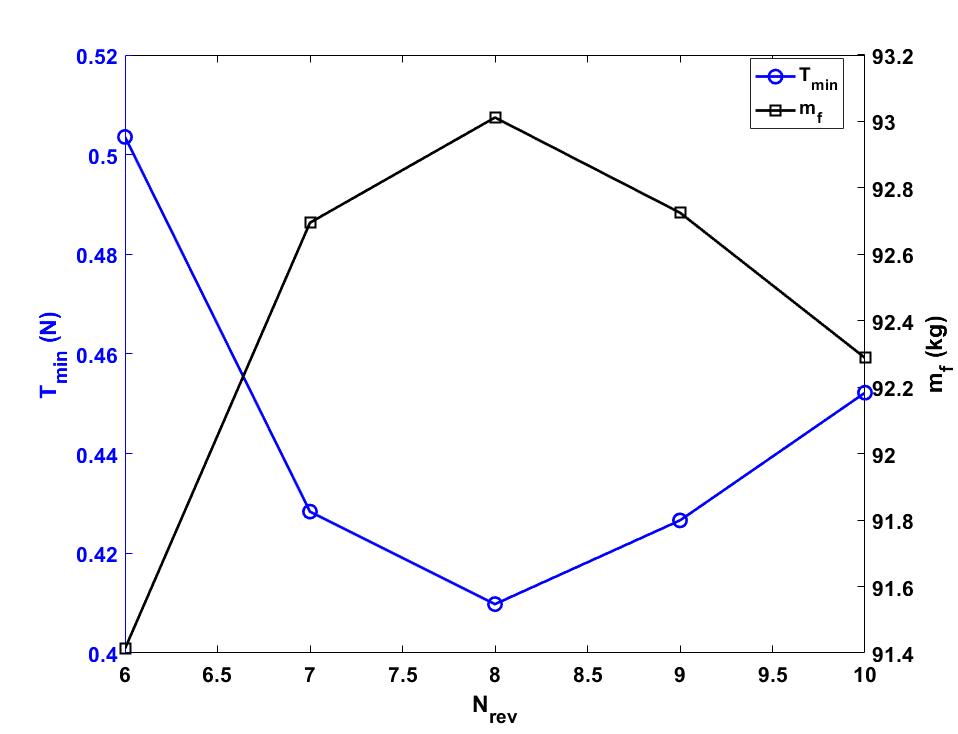}
\caption{Changes in $T_{\text{min}}$ and $m_f$ vs. $N_{\text{rev}}$ for the GTO-to-GEO problem.}
\label{fig:GTO2GEOMinThrustNrevs}
\hfill
\centering
\includegraphics[width=0.49\textwidth]{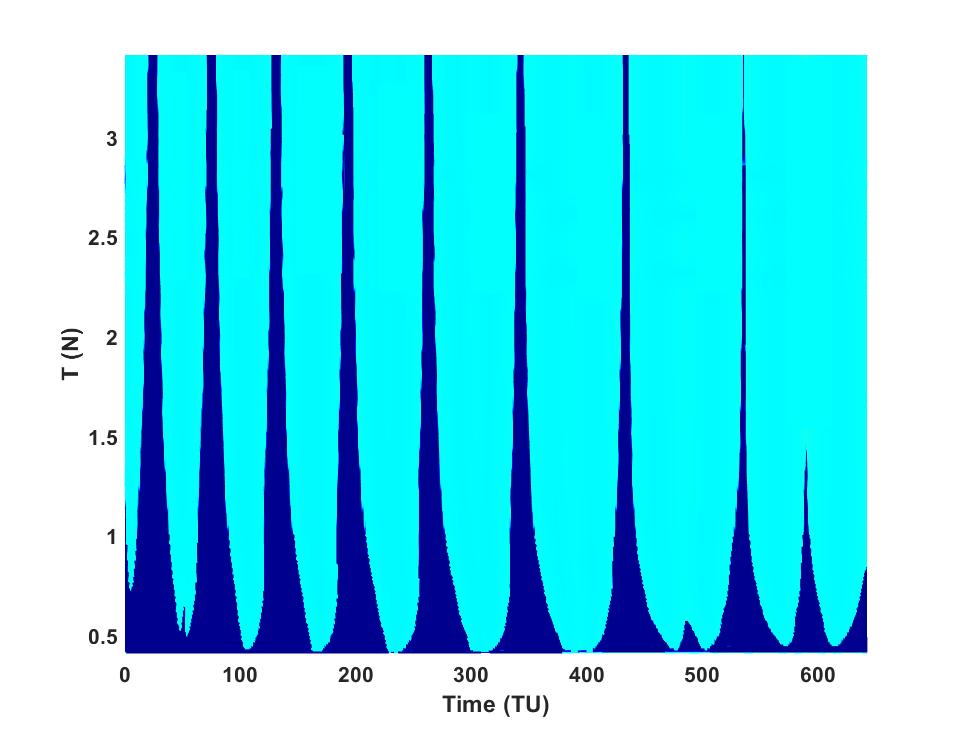}
\caption{Top view of the switching surface for the GTO-to-GEO problem with $N^*_{\text{rev}} = 8$.}
\label{fig:GTO2GEOSS}
\end{multicols}
\end{figure}
In this problem, the initial mass of spacecraft is $m_0 = 100$ kg, and the propulsion system is a low-thrust engine with a specific impulse of $I_{\text{sp}} = 3100$ seconds. It is assumed that the constant time of flight is $t_f = 6$ days. Figure \ref{fig:GTO2GEOMinThrustNrevs} shows the changes in minimum-thrust and its associated final mass for different number of en route revolutions. The fundamental minimum-thrust solution corresponds to $N^*_{\text{rev}} = 8$ with $T_{\text{min}} = 0.4098$ N. The thrust is varied within a given range of $T \in [T_{\text{min}},T_{\text{max}}]$ where $T_{\text{max}} = 3.4$ N. Figure \ref{fig:GTO2GEOSS_3} shows the fundamental switching surface and its topology over the range of thrust magnitudes from an opposite view so that the details at the low-thrust region can be seen. The switching function associated with the fundamental minimum-thrust solution is shown by a solid blue line.  This scenario is envisioned to provide a low-cost means of inserting multiple small spacecraft in GEO, by ``dropping them off'' on a GTO orbit and using electric propulsion in view of much more expensive chemical propulsion to circularize them at apogee of the GTO transfer. 

Figure \ref{fig:GTO2GEOSS} shows the switching surface $S=0$ contour map. Similar to the Earth-to-Dionysus problem, there are a number of changes and topographic features that occur near the low-thrust region of the surface. There are eight thrust ridges that remain at high thrust magnitudes all of which correspond to apogee thrusting arcs. At the lower region of the switching surface, a very small needle-like thrust arc around $t = 50$ TU is a different feature. 

\begin{figure}[htbp!]
\begin{multicols}{2}
\centering
\includegraphics[width = 0.45\textwidth]{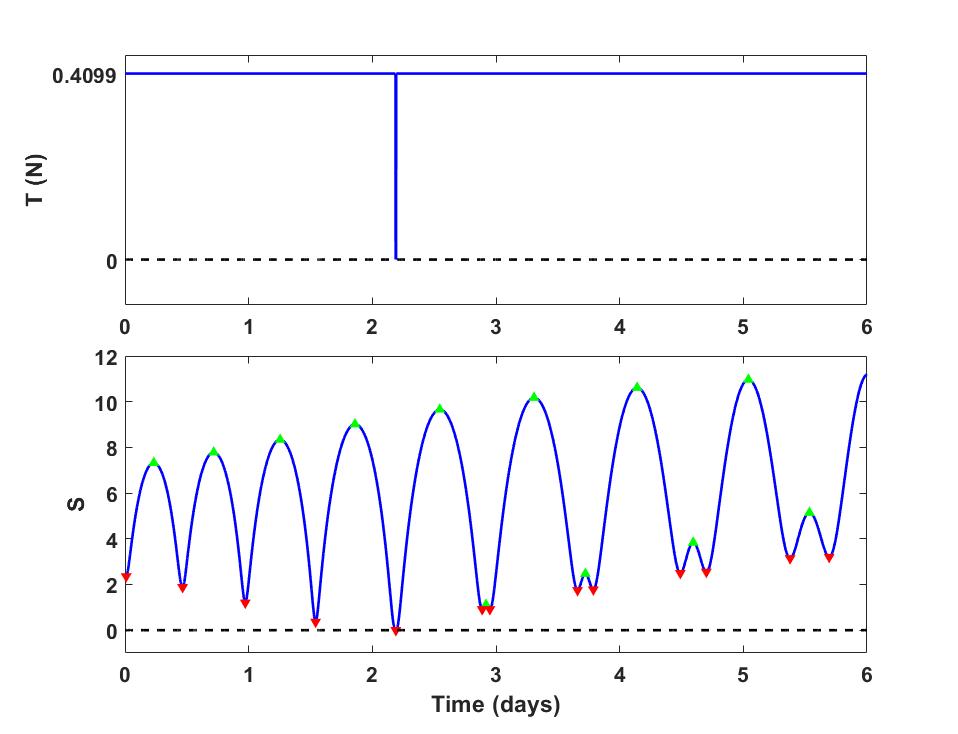}
\caption{GTO-to-GEO switching function with $T_{\text{min}} \approx 0.4098$ N.}
\label{fig:GTO2GEOCT1_SF}
\centering
\includegraphics[width = 0.45\textwidth]{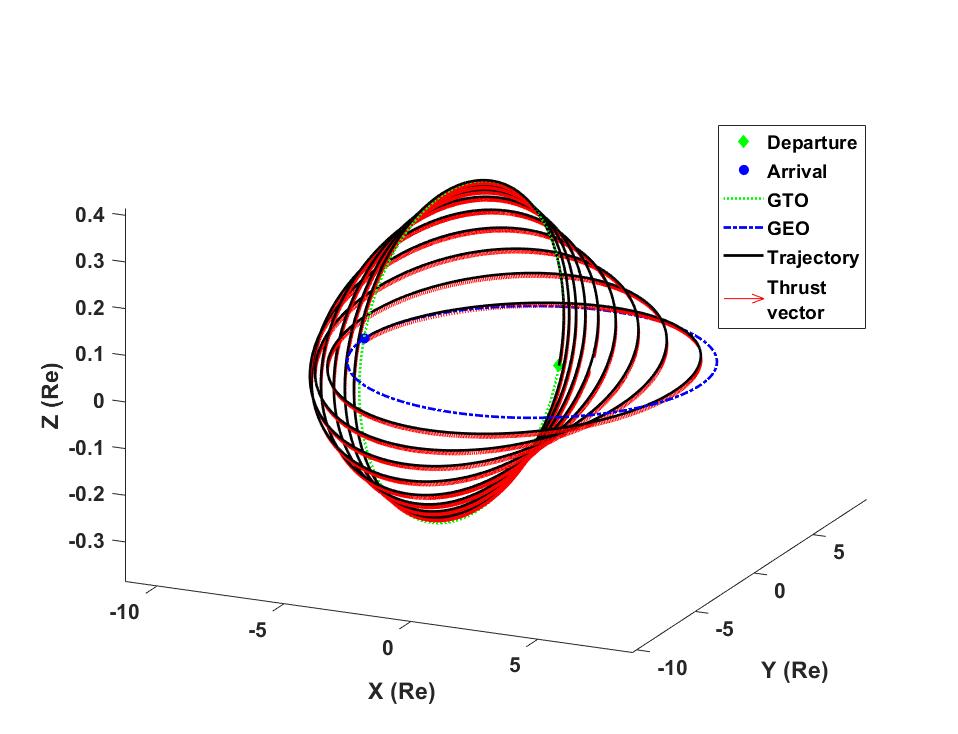}
\caption{GTO-to-GEO trajectory with $T_{\text{min}} \approx 0.4098$ N.}
\label{fig:GTO2GEOCT1_Traj}
\end{multicols}
\end{figure}

Figures \ref{fig:GTO2GEOCT1_SF} and \ref{fig:GTO2GEOCT1_Traj} show the switching function, thrust profile and trajectory for minimum-thrust solution $T_{\text{min}} \approx 0.40985$ N. Although the overall topology of the switching function seem not to possess any significant difference from the previous switching function, the surfaces features in the right-most portion of Figures \ref{fig:GTO2GEOCT1_SF} is indeed quite revealing after the explanation below is reviewed. Several of the distinct features of this switching function are due to the strength and nonlinearity of the gravitational field of the Earth compared to the previous two interplanetary problems, especially due to the eccentricity of the GTO orbit.  

One can quickly distinguish between two types of minima. Starting from left of the plot, the first five minima are individual ones, whereas the rest of the switching function contains ``tooth-like'' double minima. The former indicate the existence of perigee thrusting that might be ``believed'' (incorrectly) to vanish by increasing the thrust value. These local minima correspond to the inner-most revolutions where perigee thrusting becomes increasingly ``inefficient'' as the thrust value gets larger. The later tooth-like minima, however, represent a peculiar feature associated with them, i.e, there is a local maximum of $S$ flanked by two local minima.

What do these tooth-like features point to? Actually, as we seek to answer this question, we see now there is more to this switching function than first meets the eyes in part because the unusual phenomena lie in the region below $T \approx 0.42$ N, not shown in Figure \ref{fig:GTO2GEOSS}. The answer to what these features point to, lies, in the types of the orbit and the regions where it is most efficient to apply thrust. 

These features denote thrusting at the perigee of quasi-elliptical-shaped intermediate orbits and their characteristic symmetry (of the local minima around a local maximum) is associated with the fact that the thrusting occurs during perigee passages. These types of perigee thrusting arcs appear at higher perigee altitudes where perigee becomes less well-defined. We will show this through the numerical results, but, for now, it is possible to imagine that at some thrust magnitude, both local minima touch the $S = 0$ line at the same time. Then, a slight increase in the thrust magnitude implies that these local minima assume negative value, i.e., there is a short interval of thrusting surrounded by short coast arcs before and after the thrust arc. At some higher thrust magnitude, the local maximum will also cross $S = 0$ line and assumes a negative value, i.e., this perigee thrusting arc will vanish completely and there will be no local perigee thrust arc for large values of thrust. 

\begin{figure}[htbp!]
\begin{multicols}{2}
\centering
\includegraphics[width = 0.5\textwidth]{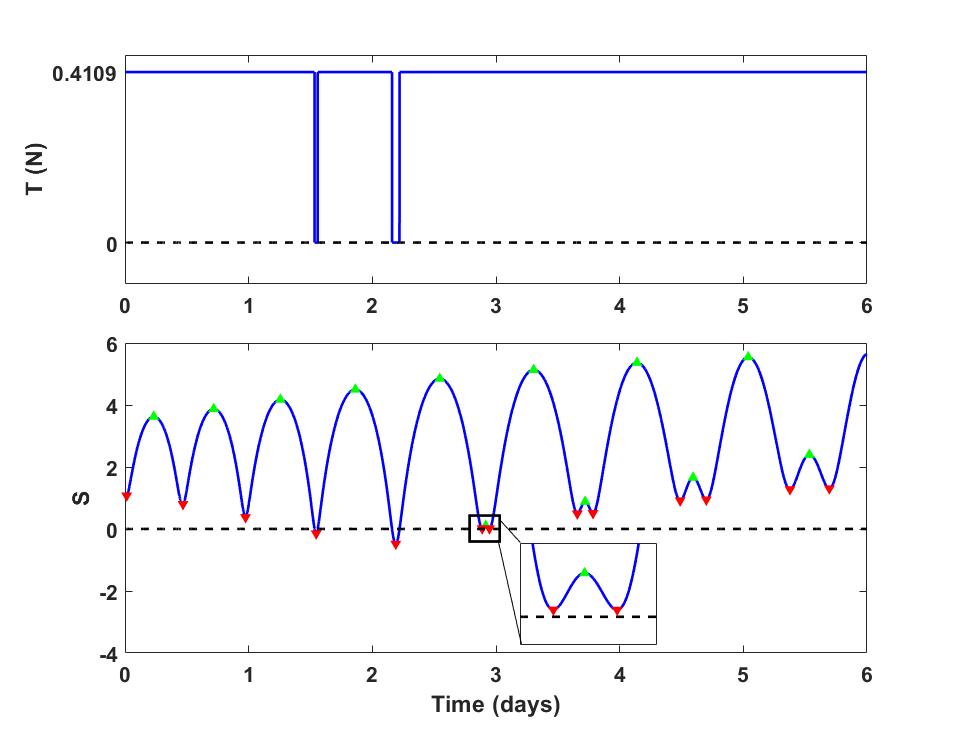}
\caption{GTO-to-GEO switching function with \\$T \approx  0.41090$ N.}
\label{fig:GTO2GEOCT2_SF}
\centering
\includegraphics[width = 0.5\textwidth]{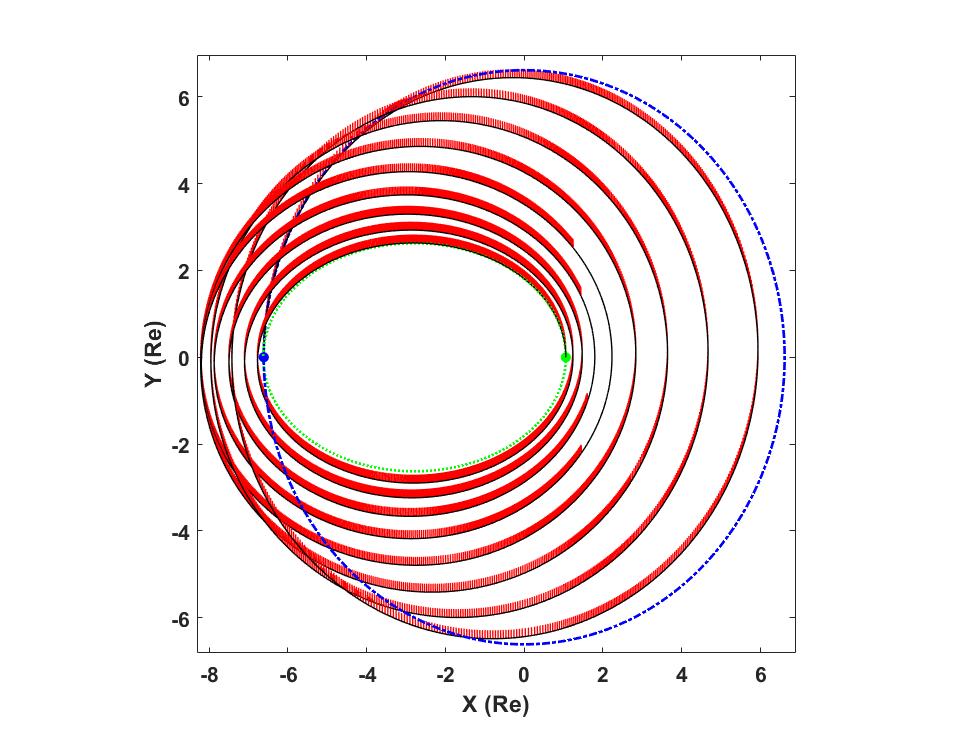}
\caption{GTO-to-GEO trajectory with $T \approx  0.41090$ N.}
\label{fig:GTO2GEOCT2_Traj}
\end{multicols}
\end{figure}

\begin{figure}[htbp!]
\begin{multicols}{2}
\centering
\includegraphics[width = 0.5\textwidth]{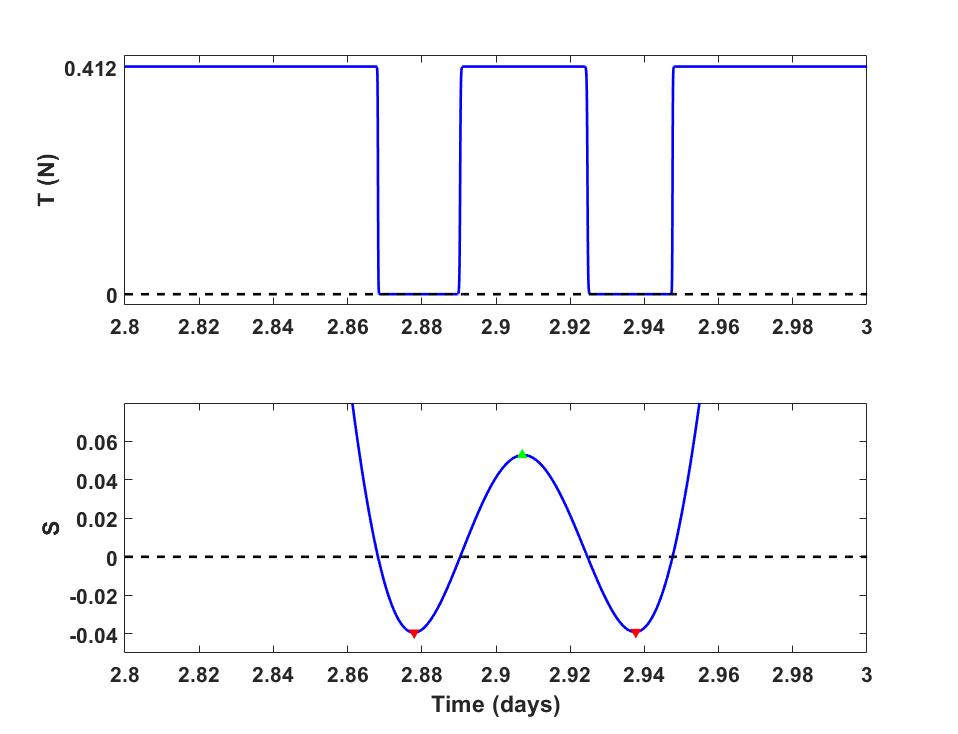}
\caption{GTO-to-GEO switching function with $T \approx  0.412$ N.}
\label{fig:GTO2GEOCT3_SF}
\centering
\includegraphics[width = 0.5\textwidth]{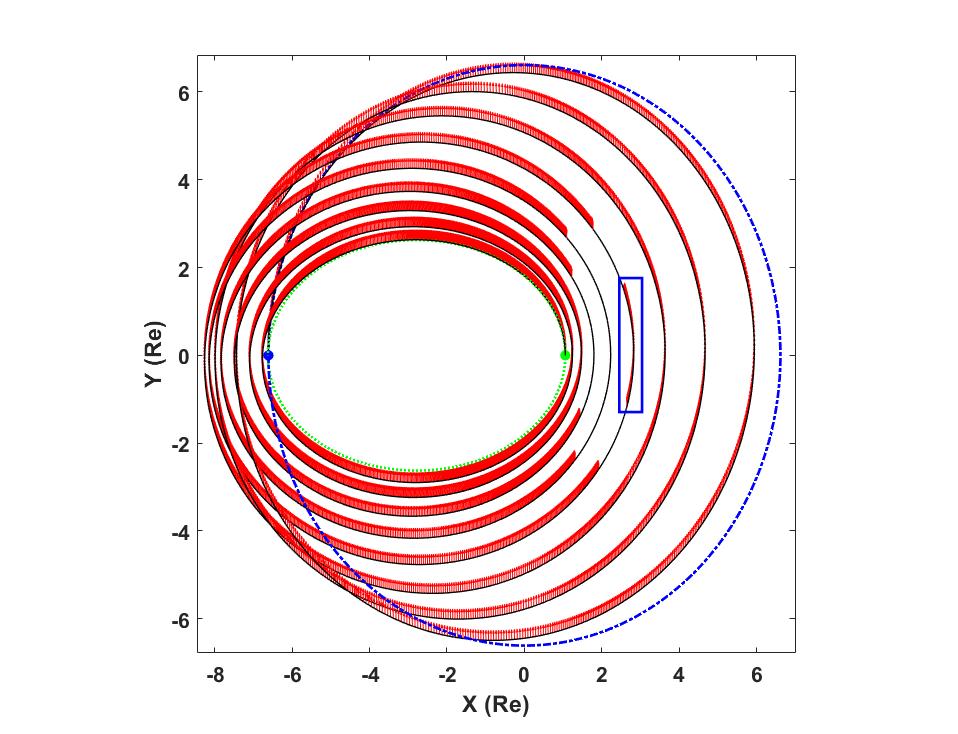}
\caption{GTO-to-GEO trajectory with $T \approx  0.412$ N.}
\label{fig:GTO2GEOCT3_Traj}
\end{multicols}
\end{figure}

Note, referring to Figures \ref{fig:GTO2GEOCT1_SF} and \ref{fig:GTO2GEOCT2_SF} we see $T_{\text{min}} = 0.4099$ N and we also see that there are four of these tooth-like features and each gets larger (in size) compared to the one to the left of it. This means that the perigee thrust arcs, at higher altitudes later in the spiral transfer, \textit{are those that will vanish last}. This makes sense physically as the perigee thrust arcs that occur at higher altitudes are more efficient (lower-velocity, less intense gravity and less gravity losses) compared to the low-altitude perigee thrust arcs (for a highly eccentric orbit). For instance, Figures \ref{fig:GTO2GEOCT3_SF} and \ref{fig:GTO2GEOCT3_Traj} show the switching function, thrust profile and trajectory for the thrust magnitude, $T \approx 0.412$ N at which the first tooth-like feature is about to touch $S = 0$ line. Note that two of the regular-formed perigee thrust arcs to its left have already disappeared (as thrust magnitude is increased). This fact can be seen clearly in the trajectory plot. 

Figure \ref{fig:GTO2GEOCT3_SF} shows the blown-up view of the switching function and its associated thrust profile. Figure \ref{fig:GTO2GEOCT3_Traj} shows the whole trajectory for the thrust magnitude, $T \approx 0.412$ N at which there is a perigee thrust arc that is surrounded by two coast arcs. A blue rectangle on the trajectory plot shows this perigee thrust arc. 

A very important aspect of these tooth-like local switch function features is that they may appear (created at other regular perigee thrust arcs) as the thrust magnitude is increased. For instance, the fact that at first, there are only four of these features (see Figure \ref{fig:GTO2GEOCT1_SF}) does not mean that their number will remain the same throughout the process of sweeping the thrust magnitudes. Put in another way, it is possible (for various boundary conditions and system parameter variations) that even perigee thrust arcs at the inner-most revolutions also reveal such tooth-like features; it is due to the fact that thrust magnitude is increased to a level that the optimality conditions demand short perigee thrusting structure. 

\begin{figure}[htbp!]
\centering
\includegraphics[width=5.0in]{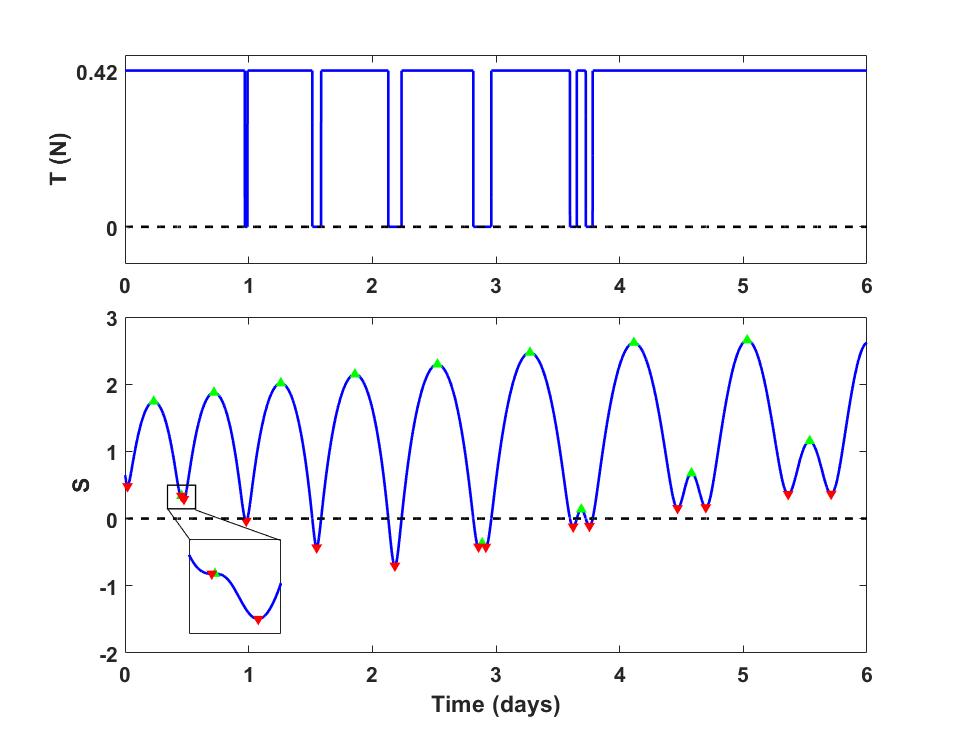}
\caption{Switching function for $T = 0.42$ N where a tooth-like feature is created.}
\label{fig:GTO2GEOCT4_SF}
\end{figure}

The creation of these tooth-like thrust level features qualitatively denote of an impending (or a tendency towards) short perigee thrust and coast arcs. For instance, Figure \ref{fig:GTO2GEOCT4_SF} shows the thrust magnitude, $T = 0.42$ N at which the previously believed-to-be regular perigee thrust arc is now modifying its shape to become a tooth-like feature! At this point, we now have quantitative clues on why a small (needle-like) thrust ridge exists on the lower-left region of the switching surface (at $t = 50$ TU). 

By inspecting Figure \ref{fig:GTO2GEOSS}, it is easy to confirm that the life of the osculating perigee thrust ridge, at $t \approx 500$ TU, is shorter than this needle-like thrust at a lower osculating perigee altitude since the former vanishes at a lower magnitude of the thrust. Eventually, the needle-like perigee thrust ridge also vanishes. The next important phenomenon occurs at a thrust magnitude at which the last thrust arc detaches from $t_f$, i.e., the time of flight is dictated by the time of the last thrust ridge. 

There is, however, an important distinction that is worthy of explanation. Unlike the previous optimal transfer example problems, the last presumably perigee thrust ridge remains, even for high thrust levels. In fact, the distinction between apogee and perigee thrust ridges is becoming mute at high altitude because of the fact that the trajectory is becoming near-circular, i.e., $e \approx 0$ to satisfy the final orbit boundary conditions. This last very small thrust ridge acts as the final thrust kick to circularize the orbit. For the sake of brevity, the details of each of the above critical thrust magnitudes and their respective switching functions and trajectories are omitted. After all of the perigee thrust arcs vanish, the remaining thrust arcs denote pure apogee thrust arcs.   

%The second one lies on the right-most point of the apogee thrust arc, which will become the main early-arrival boundary as the thrust magnitude is increased. Note also that the problem has an early-departure boundary. While it is shown that there could be multiple early-arrival boundaries, we have not encountered multiple late-departure boundaries for obvious reasons! Speculation on the existence of multiple late-departure boundaries is irrelevant since the spacecraft has only one chance of departure and any other thrust ridge is merely a thrust arc. In this context, it does not make sense to leave Earth twice!
An interesting aspect of low-thrust optimal solutions for the GTO-to-GEO problem that include plane change is that the apogee is initially raised to beyond that of the GEO orbit and a series of apogee and perigee thrust arcs are performed to achieve an optimal trajectory depending on the thrust magnitude. The results confirm also the fact that, significant portion of the thrust arcs at high thrust levels are performed near apogee to achieve simultaneous changes in the shape, and especially, the orientation of the orbit. 

\begin{figure}[htbp!]
\centering
\includegraphics[width=5.0in]{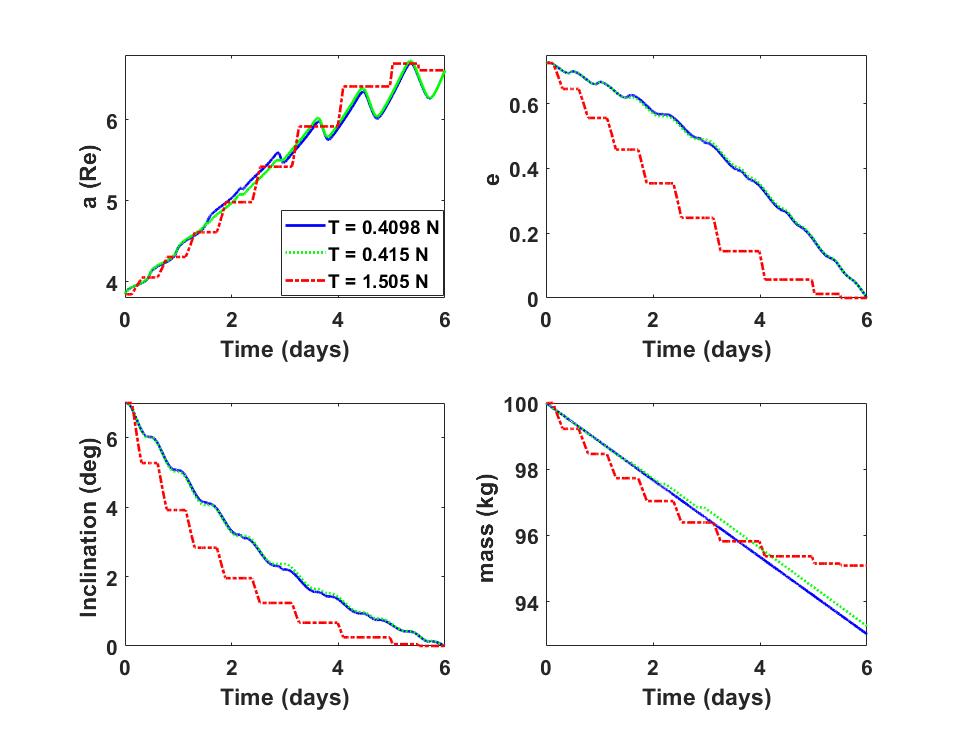}
\caption{Time history of a number of important orbit parameters ($i$, $e$, $a$) and mass for three different thrust magnitudes.}
\label{fig:GTO2GEO_ThreeCasesCOE}
\end{figure}

Figures \ref{fig:GTO2GEO_ThreeCasesCOE} shows the time history of three osculating classical orbit elements along with the mass of the spacecraft for three different thrust values. Figure \label{fig:GTO2GEO_ThreeCasesra} shows the variation of the osculating apogee versus time for three thrust magnitudes. The non-monotonic increase/decrease nature of osculating apogee variations versus time is remarkable, especially for low thrust values. Note that the apogee is raised to above GEO altitude where it is efficient to perform plane-change maneuvers. It is also interesting to see that the apogee value undergoes relatively pronounced nonlinear oscillatory trend in the second half of the maneuver, where the apogee is not reduced monotonically. In addition, the shortest transfer time belongs to a 8-impulse solution, where the impulses can be verified to occur centered on the first and last high-thrust ridges of Figure \ref{fig:GTO2GEOSS}. 

\begin{figure}[htbp!]
\centering
\includegraphics[width=5.0in]{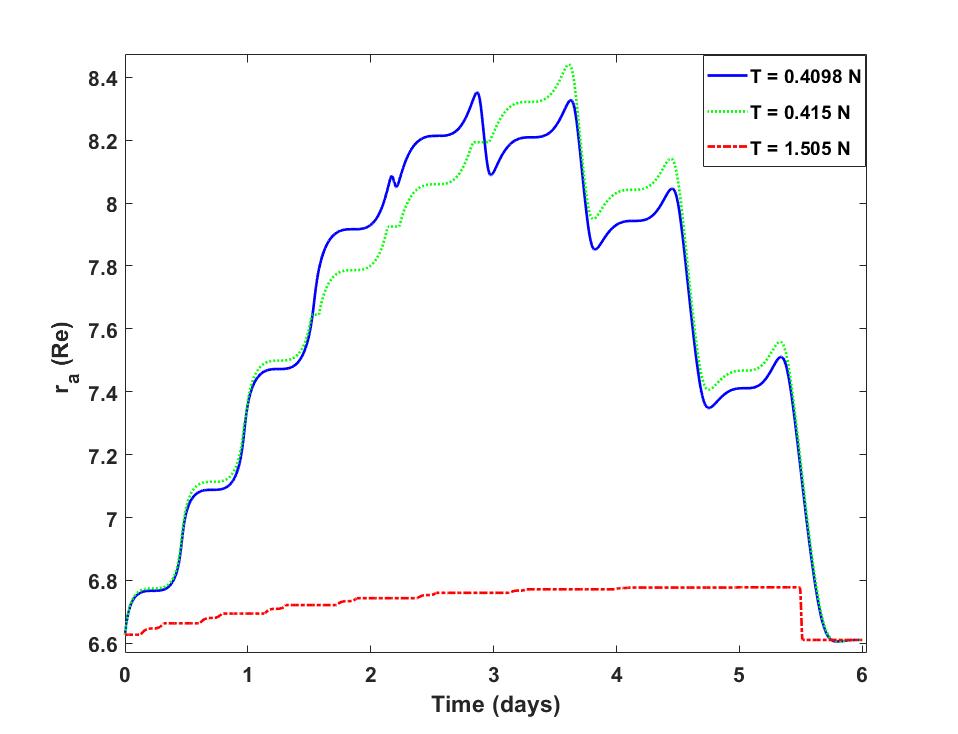}
\caption{Changes in the value of the apogee radius, $r_a = a(1+e)$ versus time of flight for three thrust magnitudes.}
\label{fig:GTO2GEO_ThreeCasesra}
\end{figure}

All of the explanations are presented to shed some quantitative light on why we see these results are based on the evident existence of these optimal solutions and our effort to interpret them through basic orbital mechanics principles. There is no direct way of guessing such subtle variations in the control structure; in fact, if this information was available and easy to predict, there was no need to solve an optimal control problem! However, the ultimate interpretation of the optimal maneuver results comes from the fact that they are based on first principles. Nonetheless, it is useful to make heuristic physical and geometrical arguments to appreciate the results. The main point we make here is that using the systematic analysis of the switching surface topography allows us to study the behavior of the family of optimal maneuvers in more detail than has been heretofore feasible. 

\begin{figure}[htbp!]
\begin{multicols}{2}
\centering
\includegraphics[width=0.5\textwidth]{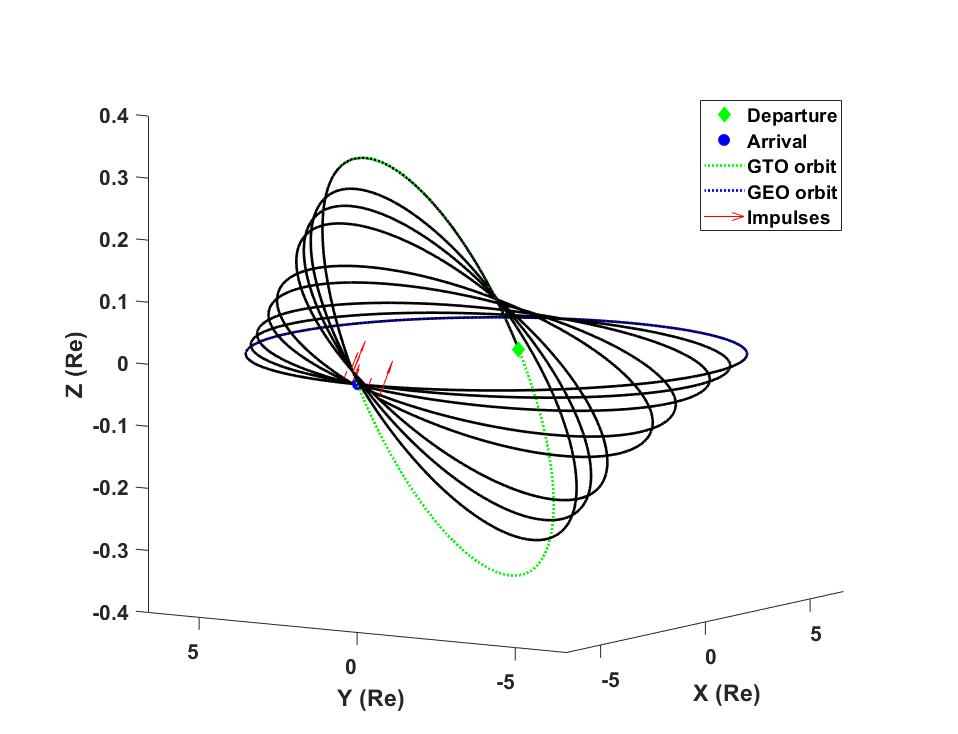}
\caption{Minimum-$\Delta v$ 8-impulse trajectory for the GTO-to-GEO problem.}
\label{fig:GTO2GEO_8ImpulseTraj}
\hfill
\centering
\includegraphics[width=0.51\textwidth]{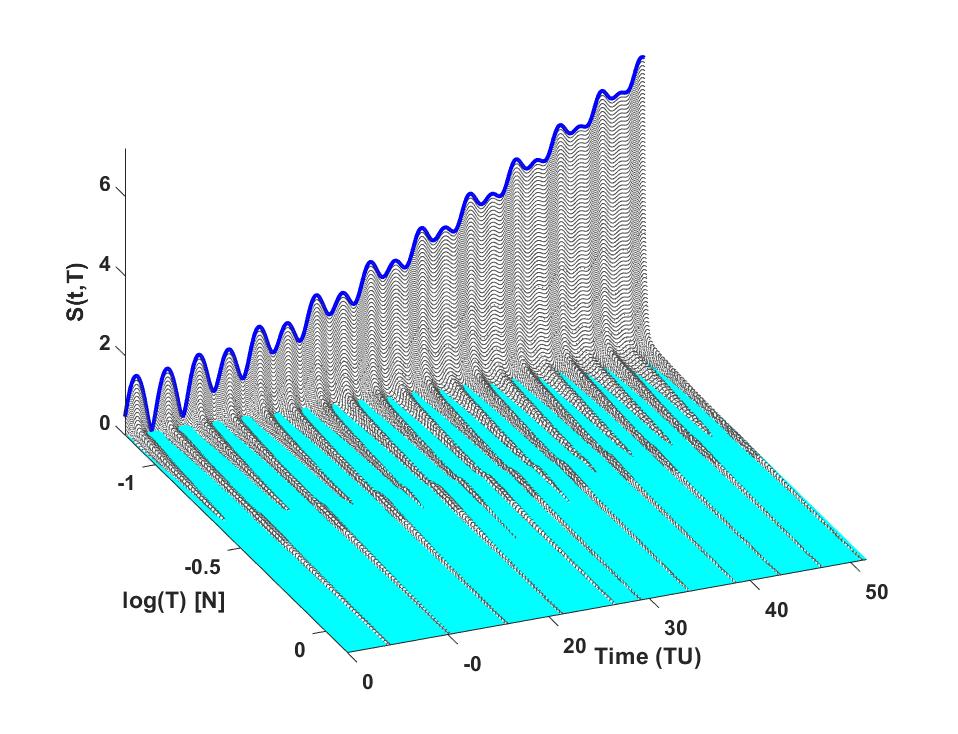}
\caption{Switching surface for the Earth-to-Venus problem with $N^*_{\text{rev}} = 10$.}
\label{fig:EV_SS_3}
\end{multicols}
\end{figure}

Table \ref{tab:impulsiveGTO2GEO} summarizes the times and magnitudes of the impulses. This problem is challenging where the conventional impulsive approaches encounter significant numerical difficulty in converging to the solution of such problems. Figure \ref{fig:GTO2GEO_8ImpulseTraj} shows the details of the minimum-$\Delta v$ 8-impulse trajectory for the GTO-to-GEO problem where all of the impulses are applied at the apogee of the intermediate elliptical orbits. These apogee kicks with judicious direction (primer vector) and magnitude are responsible for a simultaneous increase in the energy and change in the inclination. The impulses are applied very close to the descending node, which happens to be near apogee in this case.  

\begin{table}[h!] 
	\begin{center} 
		\caption{Summary of the minimum-$\Delta v$ 8-impulse solution for the GTO-to-GEO problem.}\label{tab:impulsiveGTO2GEO}
		{\small%\scriptsize
		\begin{tabular}{ c c c c c c c c c c}
        \hline
        \hline
         impulse \# &  1 & 2 & 3 & 4 & 5 & 6 & 7 & 8 & $\sum_1^8 \Delta v_i$\\
          \hline                         % &    & \cline{3-8}\\
         $t_i$ (days)        &  0.213642 & 0.69826 & 1.2187 & 1.80081 & 2.45636 & 3.20317 & 4.0564 & 5.0107 & -\\
         $\Delta v_i$ (km/s) &  0.299093 & 0.16197 & 0.1620 & 0.37607 & 0.14786 & 0.18396 & 0.0954 & 0.0705 & \textbf{1.49692}\\
         \hline
         \hline
        \end{tabular}
		}
	\end{center}
\end{table}

\subsection{Earth to Venus}
This example orbit transfer is to an inner planet of the Solar system where minimum-fuel multiple-revolution trajectories from the Earth to planet Venus are investigated. The following values are considered for the parameters of the spacecraft and its low-thrust propulsion system: $m_0 = 3000$ kg, and $I_{\text{sp}} = 3800$ s. The time of flight is chose as 3000 days.  The Earth position and velocity vectors at the departure are $\textbf{r}_{\oplus} = [145234429.88, 35542120.342, -249.986]^{\top}$ km, $\textbf{v}_{\oplus} = [-7.576, 28.83077, 0.00044765]^{\top}$ km/s, respectively. The Venus position and velocity vectors at the arrival are $\textbf{r}_{\text{\Venus}} = [-49025885.0573, 95580637.6882, 4137770.887]^{\top}$ km, $\textbf{v}_{\text{\Venus}} = [-31.278, -16.1786, 1.5837]^{\top}$ km/s, respectively. 

\begin{figure}[htbp!]
\begin{multicols}{2}
\centering
\includegraphics[width=0.52\textwidth]{figures/EVMinThrustNrevs.jpg}
\caption{Changes in $T_{\text{min}}$ and $m_f$ vs. $N_{\text{rev}}$ for \\the Earth-to-Venus problem.}
\label{fig:EVMinThrustNrevs}
\hfill
\centering
\includegraphics[width=0.50\textwidth]{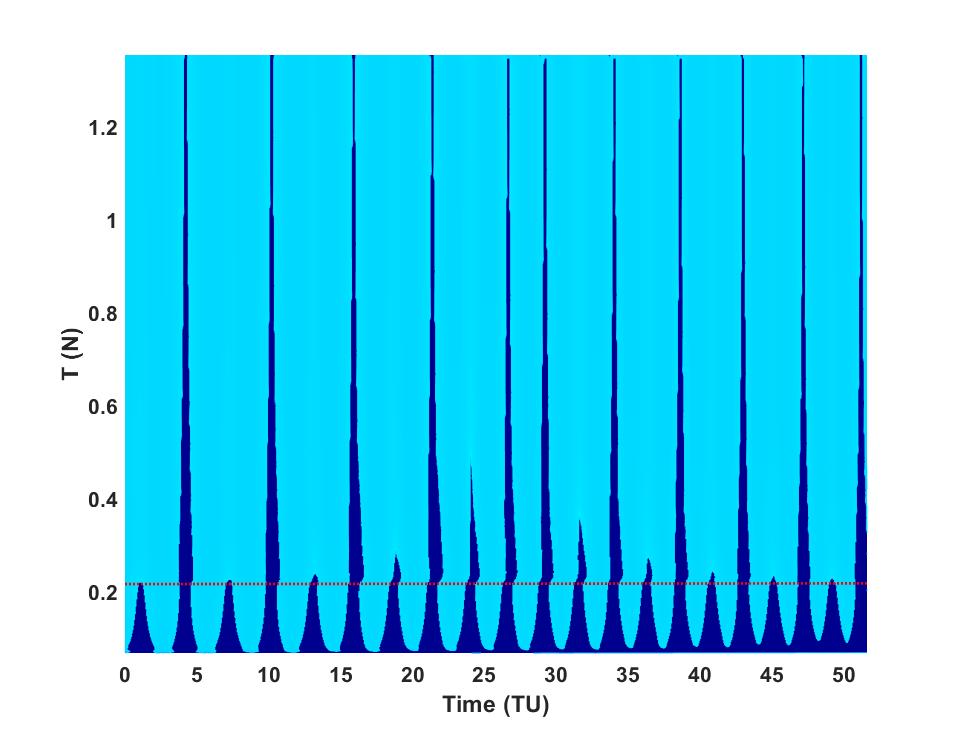}
\caption{Switching surface $S=0$ contour map for the Earth-to-Venus problem with $N^*_{\text{rev}} = 10$.}
\label{fig:EV_SS_Full}
\end{multicols}
\end{figure}

Figure \ref{fig:EV_SS_3} shows the fundamental switching surface for the Earth-to-Venus problem along with the switching function of the fundamental minimum-thrust solution. Figure \ref{fig:EVMinThrustNrevs} shows that fundamental minimum-thrust solution, which consists of ten revolutions around the Sun, $N^*_{\text{rev}} = 10$. Figure \ref{fig:EV_SS_Full} shows the switching surface $S = 0$ contour map of the Earth-to-Venus problem where it is easy to recognize the existence of both late-departure and early-arrival boundaries. Once again, we notice some irregular features in the low-thrust region but very regular behavior at hight thrust. There are eleven thrust ridges that remain for increasing thrust, approaching straight lines locating the times for impulsive thrusts.
% \begin{figure}[htbp!]
% 	\centering
%     \includegraphics[width=4.0in]{figures/EV_SS_3.eps}
% 	\caption{Switching surface for Earth-to-Venus problem for $N^*_{\text{rev}} = 10$.}
% 	\label{fig:EV_SS_3}
% \end{figure}

% \begin{figure}[htbp!]
% \centering
% \includegraphics[width=5.0in]{figures/EVMinThrustNrevs.eps}
% \caption{Changes in $T_{\text{min}}$ and $m_f$ vs. $N_{\text{rev}}$ for Earth-to-Venus problem.}
% \label{fig:EVMinThrustNrevs}
% \end{figure}
% \begin{figure}[htbp!]
% 	\centering
%     \includegraphics[width=4.0in]{figures/EV_SS_Full.eps}
% 	\caption{Switching surface for Earth-to-Venus problem.}
% 	\label{fig:EV_SS_Full}
% \end{figure}
Table \ref{tab:impulsiveEV} lists the times and magnitudes of the impulses for this many-revolution, many-impulse solution. The considered time of flight and number of revolutions is, of course, formidable. However, the outlined procedure and the proposed construct allows us to solve such challenging problems in a systematic manner.

\begin{figure}[htbp!]
\centering
\includegraphics[width=4.0in]{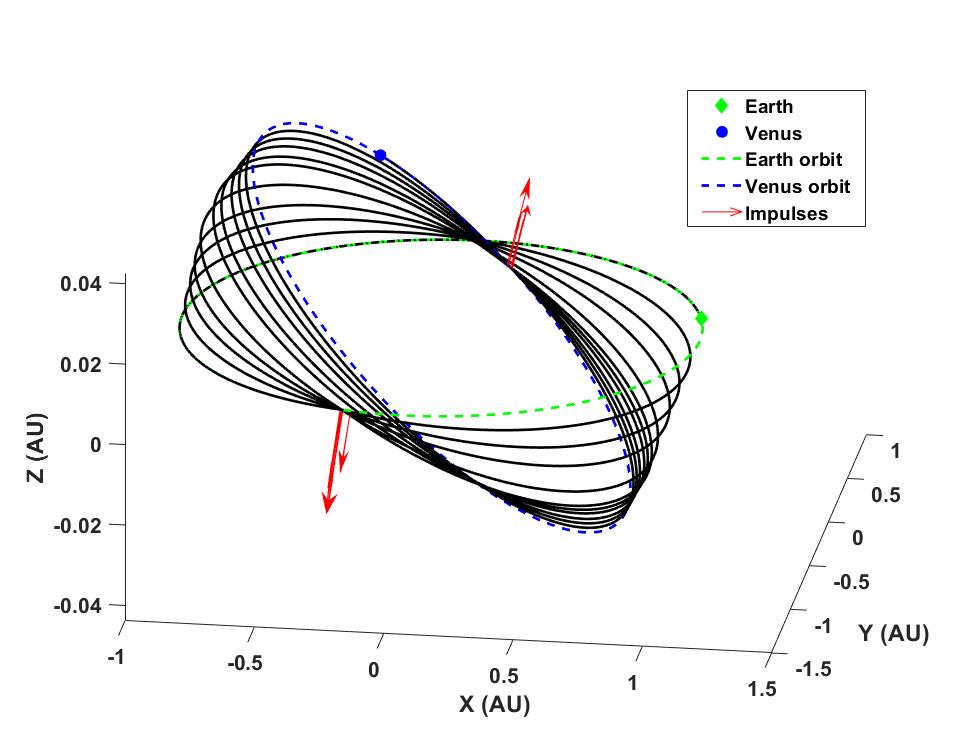}
\caption{Minimum-$\Delta v$ 11-impulse trajectory for the Earth-to-Venus problem.}
\label{fig:EV_11ImpulseTraj}
\end{figure}

\begin{table}[h!] 
	\begin{center} 
		\caption{Summary of the minimum-$\Delta v$ 11-impulse solution for the Earth-to-Venus problem.}\label{tab:impulsiveEV}
		{\small%\scriptsize
		\begin{tabular}{ c c c c c c c}
        \hline
        \hline
         \textbf{impulse} \# &  1 & 2 & 3 & 4 & 5 & 6 \\
          \hline                         
         $t_i$ (days)        & 243.6004 & 591.1218 & 923.212 & 1243.025 & 1548.372 &  1695.807 \\
         $\Delta v_i$ (km/s) & 0.562219 & 0.528433 & 0.4980  & 0.650717 & 0.361844 &  0.495807 \\
         \hline
         \textbf{impulse} \# &  7 & 8 & 9 & 10 & 11 & $\sum_1^{11} \Delta v_i$\\
          \hline                         
         $t_i$ (days)        &  1976.665 & 2244.158 & 2497.798 & 2741.055 & 2974.710 & - \\
         $\Delta v_i$ (km/s) &  0.447318 & 0.596977 & 0.462952 & 0.458392 & 0.470971 & \textbf{5.533633}\\
         \hline
         \hline
        \end{tabular}
		}
	\end{center}
\end{table}
Figure \ref{fig:EV_11ImpulseTraj} shows the details of the minimum-$\Delta v$ 11-impulse trajectory for the Earth-to-Venus problem where the impulses are split between the ascending and descending nodes of the line of nodes of the orbits of the Earth and Venus. The arc denoted by a different color denotes the separatix which acts as a boundary where the impulses switch from being applied only at aphelion to being applied only at perihelion passages. 

\section{Remarks On $N^*_{\text{\MakeLowercase{rev}}}$ and Optimal Number of Impulses} \label{sec:OpN}
In the studied test cases, the fundamental switching surface is used for determining the number of impulses. The switching surface corresponding to the fundamental minimum-thrust solution is important since the extremal field map is guaranteed to be the minimum of minima for continuous-thrust minimum-fuel trajectories. 

However, there is no guarantee that the impulsive solution obtained from the fundamental switching surface  has the \textit{smallest} $\Delta v$! For instance, for the Earth-to-Dionysus problem, we considered a high-thrust minimum-fuel solution for different number of revolutions, i.e., $N_{\text{rev}} \in \{4,5,6,7\}$ (case with $N_{\text{rev}} = 5$ is already considered) and used their associated switching function to generate impulsive solutions. The results in Table \ref{tab:impulsiveED_DiffNrev} clearly indicate that there are multiple solutions with (exact to seven digits) the same total required $\Delta v$ where it is also possible to have four, six and seven impulses. Note that the solution with $N_{\text{rev}}=4$ and only 4 impulses has the shortest time of flight, $t_f = 1840.314$ days, which means that the spacecraft can reach the asteroid 4.637 years earlier with the same $\Delta v$. Obviously, this 4-impulse solution is the preferred solution due to the shorter time of flight. 

\begin{figure}[htbp!]
\begin{multicols}{2}
\centering
\includegraphics[width=0.5\textwidth]{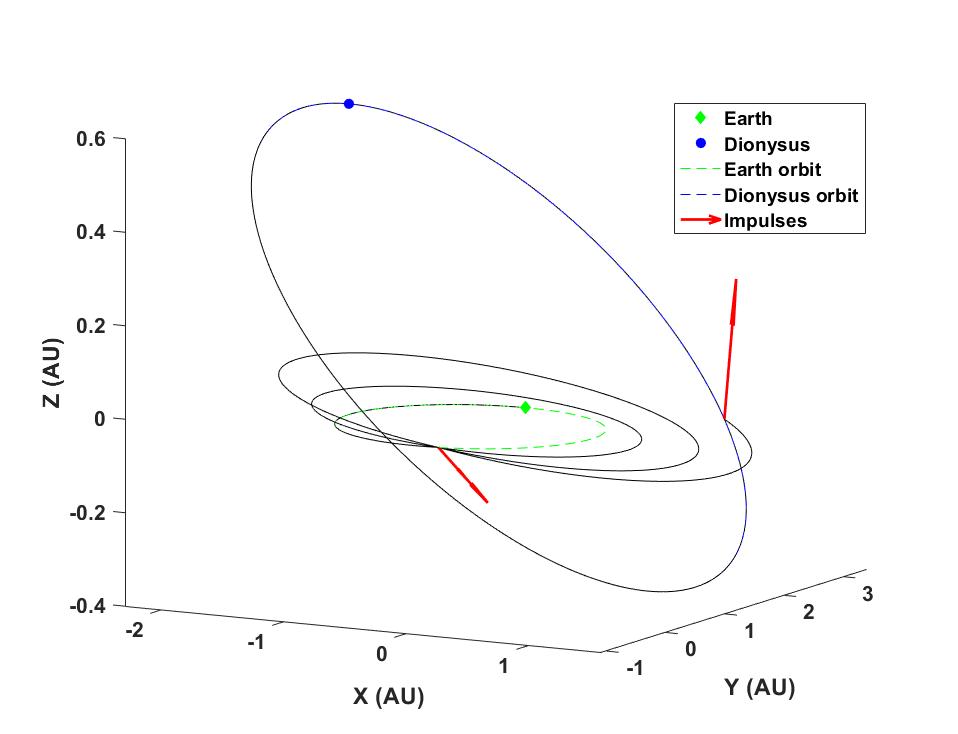}
\caption{Impulsive trajectory of the Earth-to-Dionysus problem with four impulses.}
\label{fig:EDTrajNrev4}
\hfill
\centering
\includegraphics[width=0.51\textwidth]{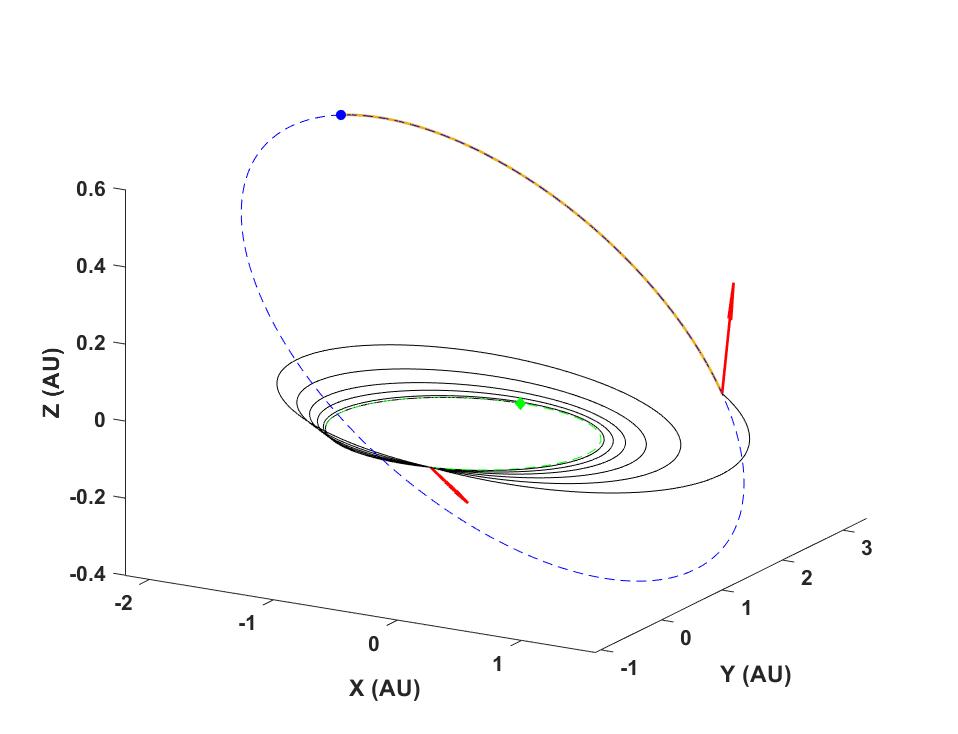}
\caption{Impulsive trajectory of the Earth-to-Dionysus problem with seven impulses.}
\label{fig:EDTraj_Nrev6}
\end{multicols}
\end{figure}

\begin{figure}[htbp!]
\centering
\includegraphics[width=4.0in]{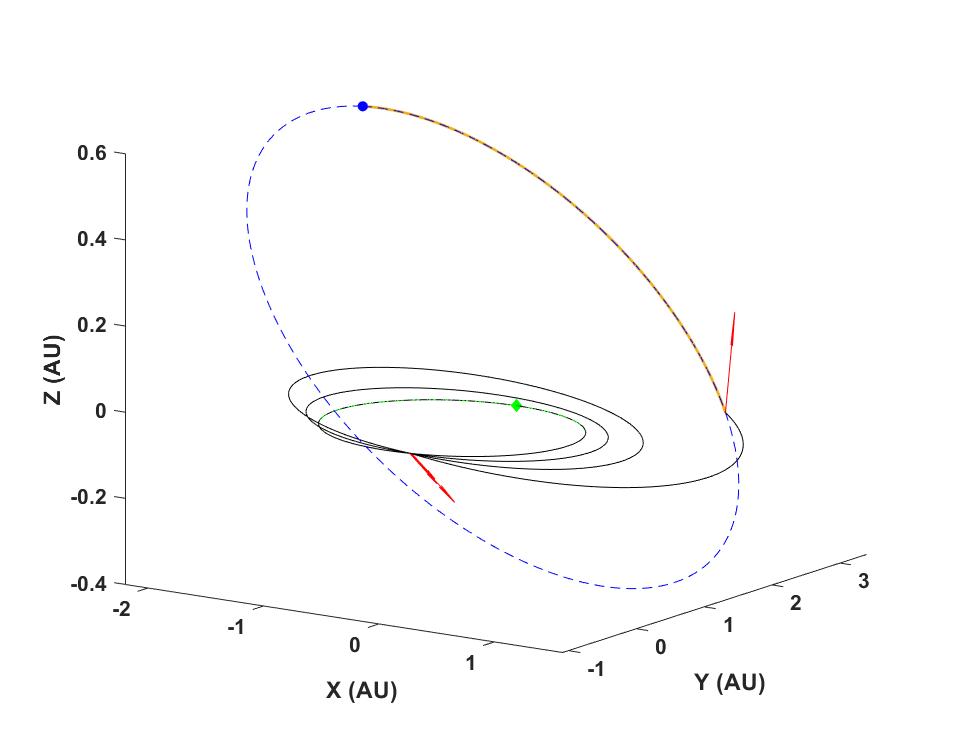}
\caption{Impulsive trajectory of the Earth-to-Dionysus problem with five impulses.}
\label{fig:EDTraj_Nrev7}
\end{figure}

Figures \ref{fig:EDTrajNrev4} and \ref{fig:EDTraj_Nrev6} show the trajectories associated with the impulsive solutions with four and seven impulses (the first and the third rows in Table \ref{tab:impulsiveED_DiffNrev}). Figure \ref{fig:EDTraj_Nrev7} shows the impulsive solution with five impulses.

\begin{table}[h!] 
	\begin{center} 
		\caption{Summary of the minimum-$\Delta v$ impulsive solutions for the Earth-to-Dionysus problem with different values for $N_{\text{rev}}$.}\label{tab:impulsiveED_DiffNrev}
		{\small%\scriptsize
		\begin{tabular}{c c c c c c c c c c}
        \hline
        \hline
        $N_{\text{rev}}$ & impulse \# &  1 & 2 & 3 & 4 & 5 & 6 & 7 & $\sum \Delta v_i$\\
          \hline                         % &    & \cline{3-8}\\
         \multirow{2}{*}{4} &$t_i$ (days)        & 193.246 &   695.199 &  1491.851 & 1840.314 & - & - & - & -\\
                            &$\Delta v_i$ (km/s) & 2.88331 &   2.94489 &  1.69332  & 2.38590 & - & - &  - & \textbf{9.907425} \\
         \hline                   
         \multirow{2}{*}{5} &$t_i$ (days)        & 193.246 &   624.164 &  1147.393 & 1810.202 & 2683.730 & 3032.192 &  - & - \\
                            &$\Delta v_i$ (km/s) & 1.61045 &   1.58896 &  1.594121 & 1.496731 & 1.231273 & 2.385884 &  - & \textbf{9.907425} \\
         \hline                   
         \multirow{2}{*}{6} &$t_i$ (days)        & 193.251 &   569.037 &  976.565 & 1429.931 & 1970.849 & 2683.731 &  3032.192 & - \\
                            &$\Delta v_i$ (km/s) & 0.29671 &   0.79962 &  0.95828 & 1.389612 & 1.778099 & 2.299207 &  2.385884 &\textbf{9.907425}\\
         \hline                   
         \multirow{2}{*}{7} &$t_i$ (days)        & 1291.121 &   1654.287 &  2091.291 & 2683.729 & 3032.186 & - &  - & - \\
                            &$\Delta v_i$ (km/s) & 0.000005 &   1.736133 &  2.341697 & 3.443678 & 2.385911 & - &  - &\textbf{9.907425}\\                            
         \hline
         \hline
        \end{tabular}
		}
	\end{center}
\end{table}

The results clearly indicate that there are multiple impulsive solutions with different number of impulses, whereas the total $\Delta v$ is the same. In this problem, we conclude that 4 local extrema (quite distinct multi-impulse transfer orbits!) give the identical minimizing total $\Delta v$, therefore we have answered Edelbaum's question, but the answer is not unique! Note that each solution has a different number of impulses. 

This curious result is consistent with some planar cases in the literature \cite{luo2010interactive,shen2015global}. However, the results in our paper correspond to significantly more difficult multiple-revolution 3D cases and it is shown that same exact $\Delta v$ (up to seven significant digits) exist with different number of impulses. Thus, Edelbaum's question ``How many impulses?'' does not generally have a unique answer, if the total $\Delta v$ is the only performance metric. However, there are generally secondary considerations, and it is obvious that the early arrival associated with $N_{\text{rev}} = 4$ would be appealing in many circumstances. 
\begin{table}[h!] 
	\begin{center} 
		\caption{Comparison of the final mass (kg) and rendezvous time (days), $t_{\text{R}}$ for the Earth-to-Dionysus problem with different values for $N_{\text{rev}}$ and two values of thrust $T = 1$ and $T = 1.4$ N; $\rho_{\text{min}} = 9.68 \times 10^{-6}$.}\label{tab:ED_Finalmass}
		{\small%\scriptsize
		\begin{tabular}{c c c | c c | c c | c c}
        \hline
        \hline
         \multirow{3}{*}{$T$ (N)} & \multicolumn{8}{c}{$N_{\text{rev}}$} \\
                              \cline{2-9}      
                 &  \multicolumn{2}{c}{4} & \multicolumn{2}{c}{5} & \multicolumn{2}{c}{6} & \multicolumn{2}{c}{7} \\
                 %\cline{2-3} \cline{4-5} \cline{6-7} \cline{8-9}
                 &   $m_f$ & $t_{\text{R}}$ &  $m_f$ & $t_{\text{R}}$ & $m_f$ & $t_{\text{R}}$ & $m_f$ & $t_{\text{R}}$ \\
                \cline{2-9}                         % &    & \cline{3-8}\\
          1.0    &      2815.128 & 3116.09 & \textbf{2842.908} & 3089.65 & 2841.049 & 3091.25 & 2812.402 & 3119.65\\
          1.4  &      2835.229 & 1895.65 & \textbf{2849.535} & 3070.20 & 2848.675 & 3070.50 & 2834.979 & 3087.15\\
         \hline
         \hline
        \end{tabular}
		}
	\end{center}
\end{table}
It must be noted, however, that minimum-$\Delta v$ \textit{does not} equate to minimum fuel for any rocket engine with a finite thrust level, a given initial mass, and some \textit{specific} specific impulse. In fact, when the thrust level is specified, then the minimum fuel criterion, as will be shown herein, eliminates the lack of uniqueness associated with attempting to minimize $\Delta v$.
\begin{figure}[htbp!]
\centering
\includegraphics[width=5.0in]{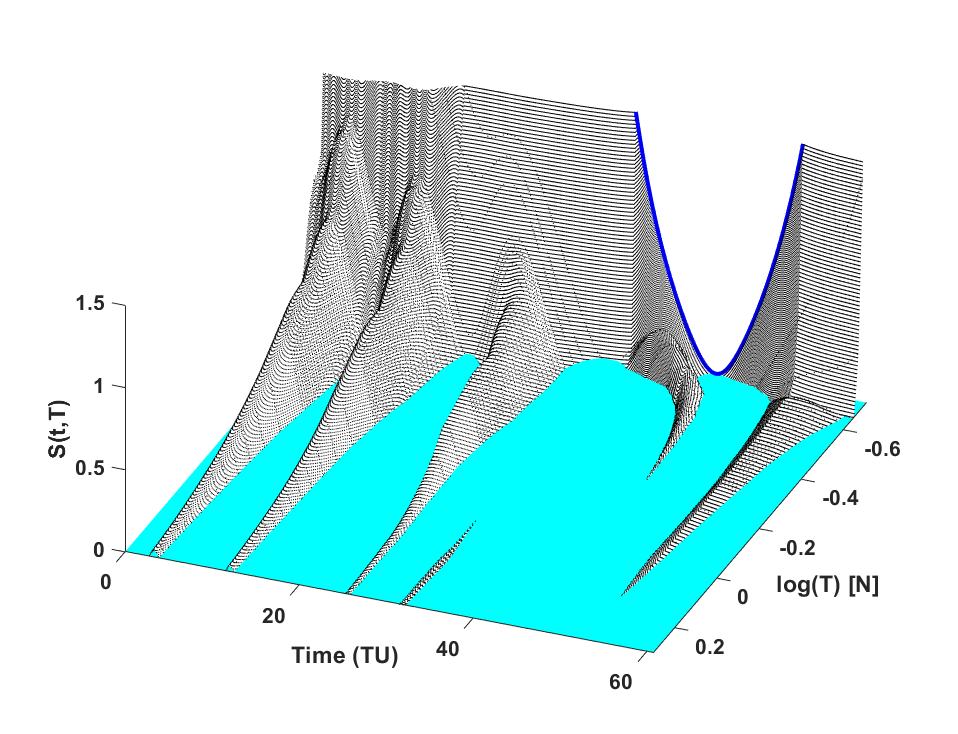}
\caption{Enlarged view of the switching surface for the Earth-to-Dionysus problem with $N_{\text{rev}}=4$.}
\label{fig:EDSS_Nrev4_2}
\end{figure}

\begin{figure}[htbp!]
\centering
\includegraphics[width=5.0in]{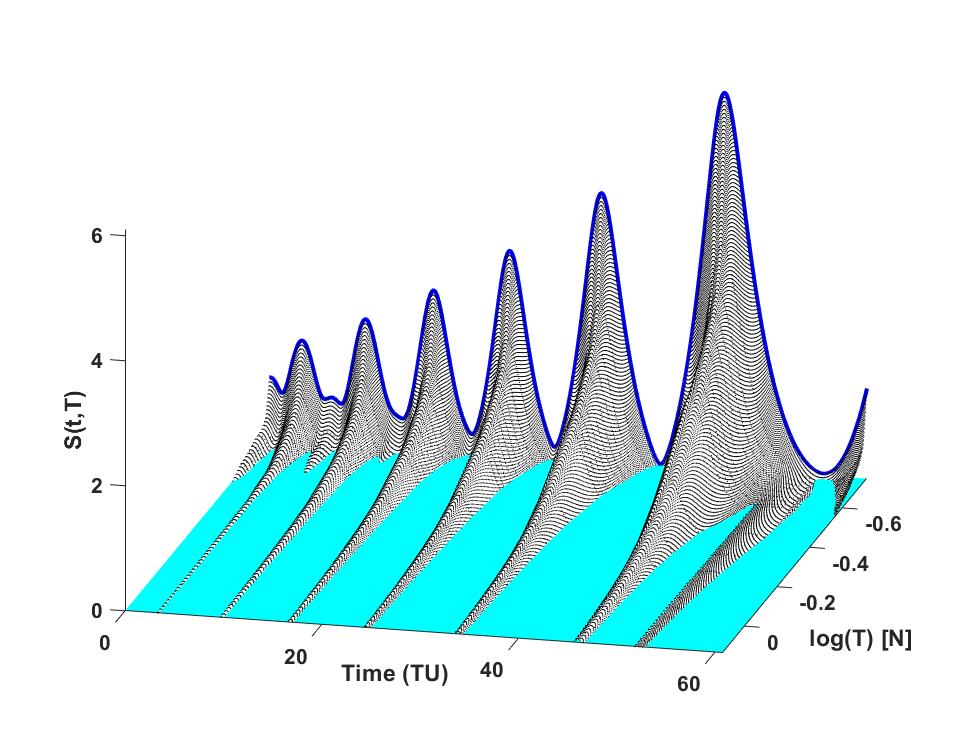}
\caption{Switching surface for the Earth-to-Dionysus problem with $N_{\text{rev}}=6$.}
\label{fig:EDSS_Nrev6}
\end{figure}

While the impulsive thrust is an idealized realization of maneuvers that ignores mass variation associated with any actual rocket motor. Observe the earlier thrust arcs have to accelerate much larger mass than do the later thrust arcs, and this physical truth allows us to find the unique optimal from the point of view of maximizing payload mass (or equivalently, minimizing propellant consumption). While the impulsive idealization remains useful in early design studies, at some point, the mission design invariably must focus on optimality for a specific payload, propellant, and engine design (or a parametric range of designs).

\begin{figure}[htbp!]
\begin{multicols}{2}
\centering
\includegraphics[width=0.50\textwidth]{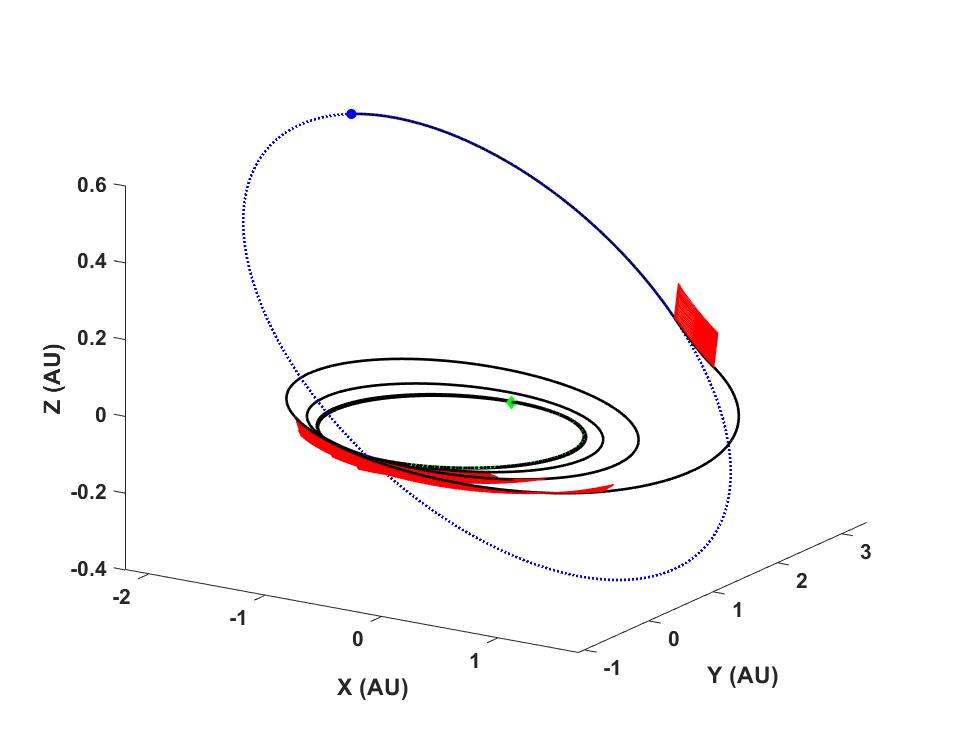}
\caption{Trajectory of the Earth-to-Dionysus problem for $T = 1$ N with $N_{\text{rev}}=4$.}
\label{fig:EDTraj_Nrev4_T1}
\hfill
\centering
\includegraphics[width=0.50\textwidth]{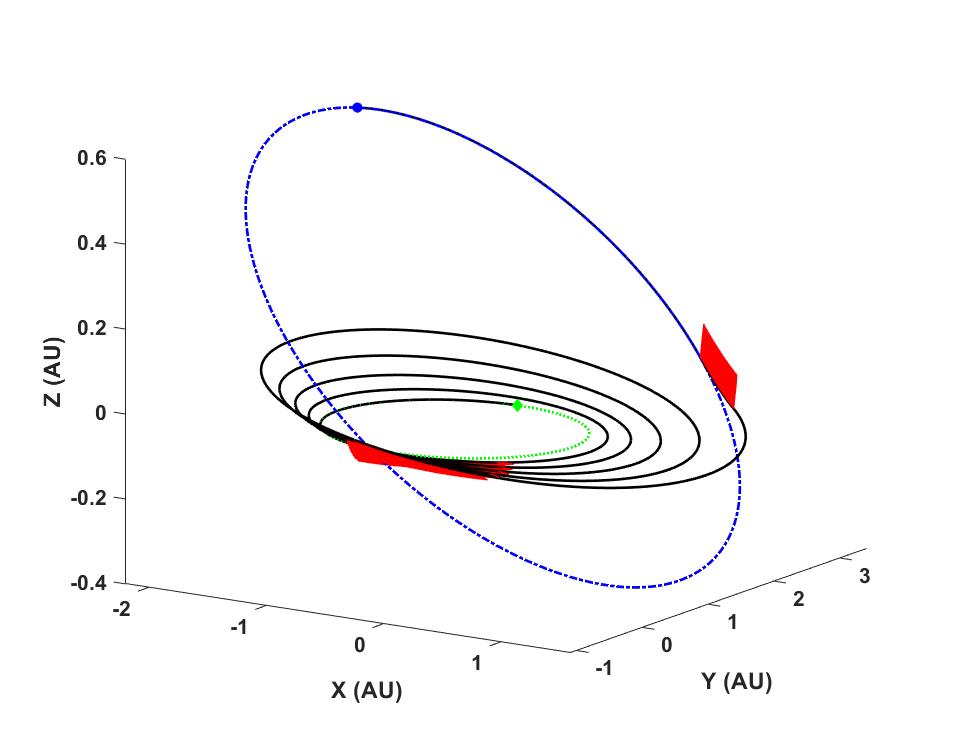}
\caption{Trajectory of the Earth-to-Dionysus problem for $T = 1$ N with $N_{\text{rev}}=5$.}
\label{fig:EDTraj_Nrev5_T1}
\end{multicols}
\end{figure}

\begin{figure}[htbp!]
\begin{multicols}{2}
\centering
\includegraphics[width=0.50\textwidth]{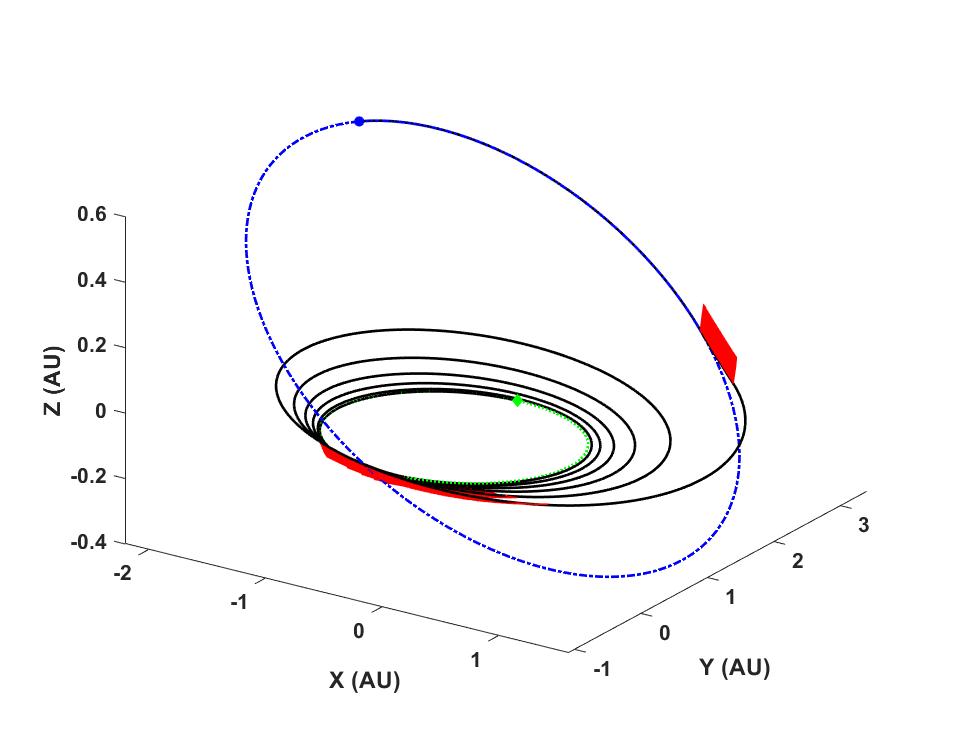}
\caption{Trajectory of the Earth-to-Dionysus problem for $T = 1$ N with $N_{\text{rev}}=6$.}
\label{fig:EDTraj_Nrev6_T1}
\hfill
\centering
\includegraphics[width=0.50\textwidth]{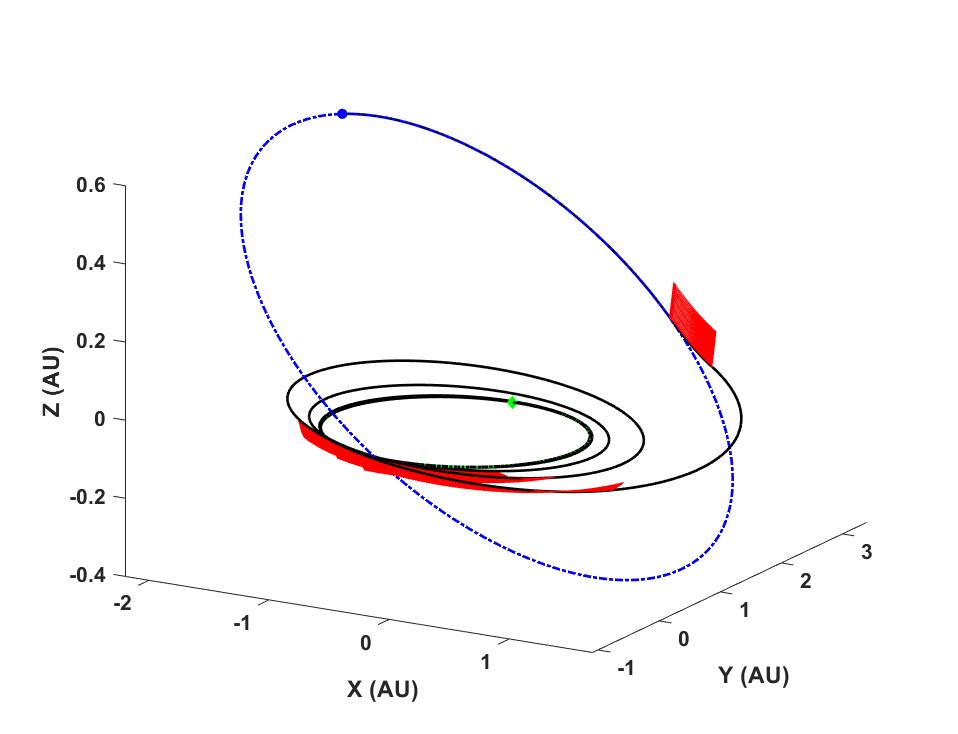}
\caption{Trajectory of the Earth-to-Dionysus problem for $T = 1$ N with $N_{\text{rev}}=7$.}
\label{fig:EDTraj_Nrev7_T1}
\end{multicols}
\end{figure}

Table \ref{tab:ED_Finalmass} compares the final mass and rendezvous time of different minimum-fuel solutions for the Earth-to-Dionysus problem with different values of $N_{\text{rev}}$ when two different thrust magnitudes are used, $T = 1$ N and $T = 1.4$ N. As we expect, at each particular thrust level, the largest final mass corresponds to a unique solution, in this case, the solution with $N_{\text{rev}} = 5$. Recall the $N_{\text{rev}} = 5$ is also the fundamental extremal which also results in the smallest thrust level that can reach the final BCs. However, when $N_{\text{rev}} = 4$, an increase in the thrust magnitude by 0.4 Newtons reduces the rendezvous time significantly, whereas only 20 more kg of additional propellant is required (see Figure \ref{fig:EDSS_Nrev4_2} for the last thrust arc, which is a result of a bifurcation). This particular test case, shows the importance of the existence of such switching surfaces. Figures \ref{fig:EDTraj_Nrev4_T1}-\ref{fig:EDTraj_Nrev7_T1} show the finite-thrust trajectories when $T = 1$ N for different values of $N_{\text{rev}}$.

\begin{figure}[htbp!]
\centering
\includegraphics[width=4.0in]{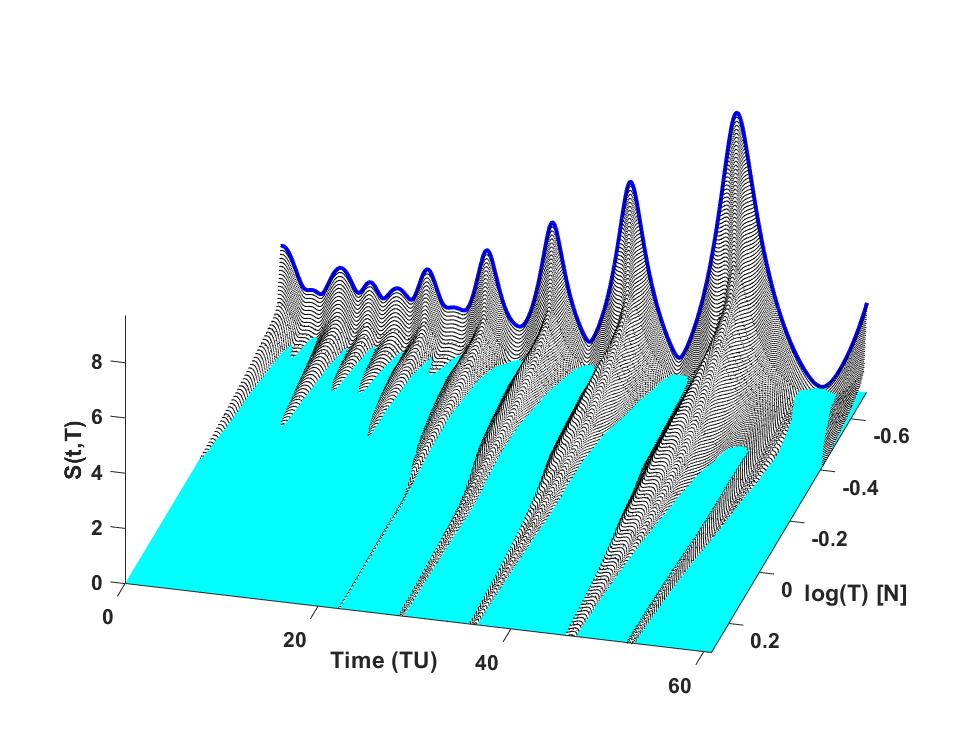}
\caption{Switching surface for the Earth-to-Dionysus problem with $N_{\text{rev}}=7$.}
\label{fig:EDSS_Nrev7}
\end{figure}

\begin{figure}[htbp!]
\centering
\includegraphics[width=4.0in]{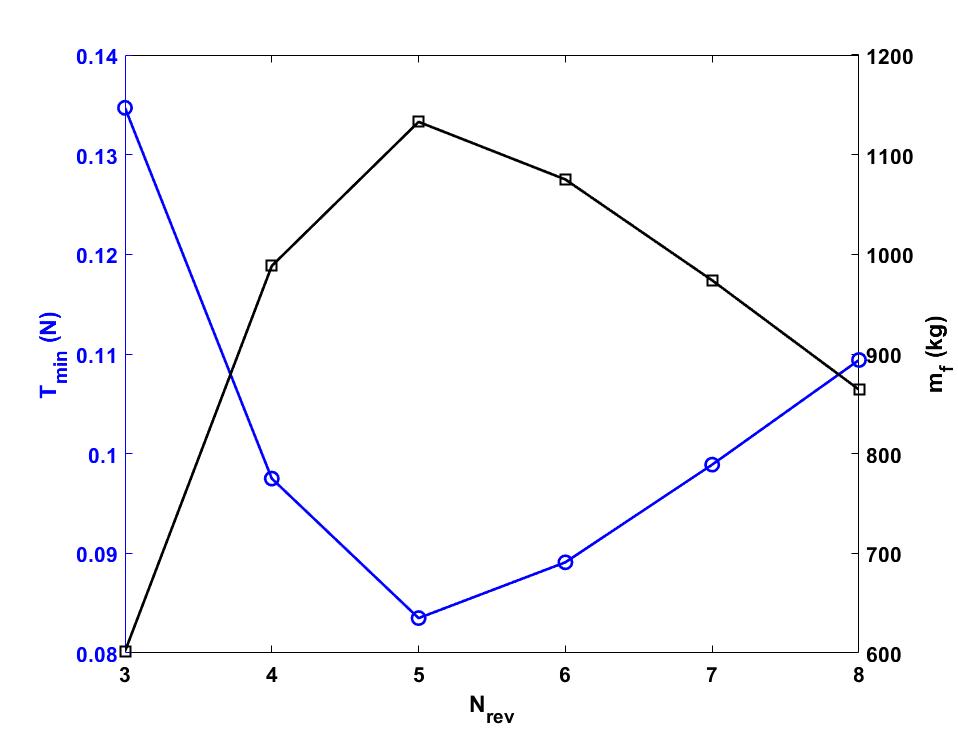}
\caption{Changes in $T_{\text{min}} $ and $m_f$ vs. $N_{\text{rev}}$ for the Earth-to-Dionysus problem; $m_0 = 2000$ kg.}
\label{fig:EDMinThrustNrevs_m02000}
\end{figure}

A reasonable question to ask is:  How does the number of revolutions change if the initial mass of the spacecraft is modified? In order to answer this question, the initial mass of the Earth-to-Dionysus problem is halved, $m_0 = 2000$ kg. Then, the minimum-thrust algorithm is invoked to satisfy Eq.~\eqref{eq:terminalconmintime} for this problem . The results are plotted in Figure \ref{fig:EDMinThrustNrevs_m02000} and show that the fundamental number of revolutions, $N^*_{\text{rev}}$ does not change. In comparison with the results in Figure \ref{fig:EDMinThrustNrevs}, the minimum-thrust values have decreased in a near-linear manner. 

To gain further confidence and perspective, we consider the corresponding extremals for other choices of $N_{\text{rev}}$. Figure \ref{fig:EDSS_Nrev4_2} shows an enlarged view of the switching surface for the Earth-to-Dionysus problem with $N_{\text{rev}} = 4$. Note that there exist a bifurcation at high thrust levels that creates an ``island'' thrust ridge. The solid blue line in the figure corresponds to the switching function of the minimum-thrust for the particular number of revolutions and touches $S = 0$ at an individual point. Figure \ref{fig:EDSS_Nrev6} shows the switching surface for the Earth-to-Dionysus problem with $N_{\text{rev}} = 6$. Figure \ref{fig:EDSS_Nrev7} shows the switching surface for the Earth-to-Dionysus problem with $N_{\text{rev}} = 7$. Based on this result and other not reported in detail here, we are confident that the 4 extremals (Table \ref{tab:ED_Finalmass}) represent the non-unique answer to Edelbaum's question (minimum total $\Delta v$), and the $N_{\text{rev}} = 5$ extremal (Table \ref{tab:ED_Finalmass}) represents the specific extremal that correspond to the spacecraft design with initial mass $m_0 = 4000$ kg, and $I_{\text{sp}} = 3000$ s and thrust of $T = 1$ N. The family of neighboring designs with variable thrust levels underly Figures \ref{fig:EDSS_3} and \ref{fig:EDSS_Side}.

Table \ref{tab:impulsiveGTO2GEO_DiffNrev} also summarizes the impulsive solutions for the GTO-to-GEO problem  for different values of $N_{\text{rev}}$. It shows that the solutions differ with negligible change in the total $\Delta v$. The case with $N_{\text{rev}} = 9$ consists of ten impulses: $t_1 = 0.220847$, $t_2 = 0.72598$, $t_3 = 1.16227$, $t_4 = 1.704234$, $t_5 = 2.26117$, $t_6 = 2.80307$, $t_7 = 3.53088$, $t_8 = 4.247474$, $t_9 = 5.02344$, $t_{10} = 5.53719$. $\Delta v_1 = 0.18774$, $\Delta v_2 = 0.00007$, $\Delta v_3 = 0.40844$, $\Delta v_4 = 0.04923$, $\Delta v_5 = 0.058367$, $\Delta v_6 = 0.203009$, $\Delta v_7 = 0.175899$, $\Delta v_8 = 0.12519$, $\Delta v_9 = 0.29177$, $\Delta v_{10} = 0.005166$. The total impulse is $\sum_{k=1}^{10} \Delta v_k = 1.50489$ km/s. In Table \ref{tab:impulsiveGTO2GEO_DiffNrev}, once again, we have three extremals with negligible difference in required $\Delta v$, so we are free to invoke secondary criteria to choose a solution. Notice $N_{\text{rev}} = 6$ results in 5 impulses, but the $4^{th}$ impulse is so small that we can consider it a 4-impulse maneuver. On the other hand, if we seek to minimize for a final thrust, we will find $N_{\text{rev}} = 8$ to be optimal because (Figure \ref{fig:GTO2GEOMinThrustNrevs}), this maneuver belongs to the fundamental minimum-fuel surface of Figure \ref{fig:GTO2GEOSS_3}.
\begin{table}[h!] 
	\begin{center} 
		\caption{Summary of the minimum-$\Delta v$ impulsive solutions for the GTO-to-GEO problem with different values for $N_{\text{rev}}$.}\label{tab:impulsiveGTO2GEO_DiffNrev}
		{\small%\scriptsize
		\begin{tabular}{c c c c c c c c c c c}
        \hline
        \hline
        $N_{\text{rev}}$ & impulse \# &  1 & 2 & 3 & 4 & 5 & 6 & 7 & 8 & $\sum \Delta v_i$\\
          \hline                         % &    & \cline{3-8}\\
          \multirow{2}{*}{6} &$t_i$ (days)        & 0.0    & 0.2316 &   0.926354 &  1.80382            & 5.510785 & - & - &  - &  \\
                             &$\Delta v_i$ (km/s) & 0.0122 & 1.0279 &   0.441644 &  $3 \times 10^{-8}$ & 0.016472 & - & - &  - & \textbf{1.49827} \\
                            \hline
         \multirow{2}{*}{7} &$t_i$ (days)        & 0.23336 &   0.72617 &  1.314734 & 2.056413 & 2.981276 & 5.505702 &  - & - & \\
                            &$\Delta v_i$ (km/s) & 0.40388 &   0.31568 &  0.434553 & 0.252489 & 0.078449 & 0.002732 &  - & - &\textbf{1.487785} \\
         \hline                   
         \multirow{2}{*}{8} &$t_i$ (days)        &  0.213642 & 0.69826 & 1.2187 & 1.80081 & 2.45636 & 3.20317 & 4.0564 & 5.0107 & - \\
                            &$\Delta v_i$ (km/s) &  0.299093 & 0.16197 & 0.1620 & 0.37607 & 0.14786 & 0.18396 & 0.0954 & 0.0705 & \textbf{1.49692}\\
                                  
         \hline
         \hline
        \end{tabular}
		}
	\end{center}
\end{table}
\vspace{-5mm}
In the literature, it is hypothesized that for multiple-revolution non-coplanar cases, a large number of $N$ impulses (with unknown times, impulse magnitudes and directions) may become more convenient for optimal rendezvous [\citen{prussing1986optimal,shen2015global}]. However, finding the optimal $N$ impulses by NLP becomes a complex task, for large $N$, through traditional approaches as discussed in [\citen{jezewski1968efficient}]. 

As documented in the examples above, the systematic approach of the current paper has been successful in solving problems in which the number of ``optimal'' impulses vary from $N^* = 3$ up to $N^* = 11$. The key is that the high thrust limit of the family of extremals that underlie the switching surface approach with high precision the impulsive limit, so we have excellent starting iteratives and know the number of extremals. We need only to refine the number of significant digits in the approximate starting solutions. The number of impulses is dictated by the high thrust asymptotic structure (number of surviving thrust ridges including those due to bifurcation) of the switching surface. According to the reported results and our experience, the optimal number of impulses depend strongly on the types of orbits, time of flight, angular phasing and the orbits’ relative orientation, and furthermore, we find that the minimum $\Delta v$ frequently can be achieved by more than one set of impulses. Thus, again, we emphasize the answer to Edelbaum's question is not unique, but answering this question usually reveals attractive solutions. 

%It would be quite interesting, of course, if one could find an ultimate explicit expression for the number of impulses. However, it seems certain that any such expression does not exist, because it must inevitably take into account the generally unknown optimal number of intermediate revolutions and associated bifurcations in the solution associated changes with variations in prescribed maximum thrust, initial and final phasing in the orbit and time of flight. 
Our studies indicate that location and number of solution bifurcations are important features of the switch surface at critical thrust values and lead to creation of thrust ridges that alter the number of persistent optimal thrust ridges (in the switching function), which affects the number of optimal thrust arcs. In our numerical examples, a thrust arc, which is created due to bifurcation at a critical value of thrust is persistent (i.e., does not vanish) with increasing maximum thrust specification. We digress briefly to make qualitative statements inferred from our analysis of the physics and numerics in the problems studied to date, using the switching surfaces as the focal point for analyzing families of extremals. We have found that the bifurcations of the switching surfaces are frequently associated with regions where either the accelerations due to the gravitational physical forces are nearly in balance ( for example near L1 point in the restricted three-body problem , not reported here), or in regions where the local ``effectiveness'' of the optimal controls is weak, in that the local variations of the thrust have a small effect on the performance index and the terminal BCs (see Earth-to-Mars problem). 

Our hypothesis is that, with all BCs fixed and for a fixed takeoff and arrival time, the optimal number of impulses $N^*$, the optimal number of revolutions, $N^*_{\text{rev}}$ associated with the fundamental minimum-thrust solution, and the number of bifurcations, $N_{\text{bifur}}$, are all properties of the switching surface; these features are discoverable by actually generating the switching surface increasing thrust from $T_{\text{min}}$ to find the high thrust asymptotic behavior. For instance, Table \ref{tab:Nop} summarizes the results of the five test cases studied in this paper. 
%It may be possible to establish an empirical relation between $N^*$, $N^*_{\text{rev}}$, and $ N_{\text{bifr}}$, however, more studies would have to be conducted to gain confidence about those relations.

While we list a specific value for $N^*$ in Table \ref{tab:Nop}, the minimum $\Delta v$ is generally achieved by several values of $N$ and the listed value is associated with one of several minimizing extremals. To choose a specific $N$ and the associated extremal, we can invoke a secondary consideration such as minimum fuel (for any chosen finite thrust engine and propellant $I_{\text{sp}}$), or time of arrival when the final thrust occurs and rendezvous is achieved at a time before the specified $t_f$. In this table, we chose $N^*$ that corresponds to the fundamental minimum-fuel solution for a family of constant thrust engine designs, with $T_{\text{min}} < T < T_{\text{max}}$.
\begin{table}[h!] 
	\begin{center} 
		\caption{Summary of the number of impulses associated with the fundamental switching surface for all considered test cases.}\label{tab:Nop}
		{\small%\scriptsize
		\begin{tabular}{ l c c c }
        \hline
        \hline
         \textbf{Test Case} &  $N^*_{\text{rev}}$ & $ N_{\text{bifr}}$ & $N^*$ \\
          \hline                         
         Earth-to-Mars       & 1 & 1 & 3 \\
         Earth-to-1989ML     & 1 & 0 & 3 \\
         Earth-to-Dionysus   & 5 & 0 & 6 \\
         GTO-to-GEO          & 8 & 0 & 8 \\
         Earth-to-Venus      & 10 & 0 & 11 \\
         \hline
         \hline
        \end{tabular}
		}
	\end{center}
\end{table}

Importantly, the possibility of terminal impulses is governed by the prescribed time of flight. As a result, the optimal number of impulses depends also on the departure and arrival times, which are in turn associated with the angular phase in the initial and target orbits. Note that the optimal number of revolutions itself depends on the time of flight, BCs, and the ratio of the acceleration of the propulsion to that of the central body. We conclude that a generally applicable, explicit, unique algebraic answer to Edelbaum's question does not exist, because we have conclusively shown that the answer is not unique (frequently a set containing several extremals with different number of impulses have the same total $\Delta v$). At some point in the mission design process we will have to consider specific propulsion systems and we abandon the idealized notion that $\Delta v$ should be minimized and instead recognize that the minimization of fuel consumed associated with a finite thrust engine is a more meaningful index, we are generally led to the true optimal trajectory of interest (or a specific family of trajectories when considering a family of spacecraft/engine designs).

As a consequence, we believe that Edelbaum's question should be more specific and refined to read:  ``For a given spacecraft initial mass and engine with a specified $I_{\text{sp}}$ and maximum thrust, what is the optimal sequence of thrusts and coasts to minimize fuel consumed when transferring from  orbit A to orbit B?''  We have developed and demonstrated a method that answers this more specific question, and in the process, we can identify a set of extremals that minimize $\Delta v$. 

On the other hand, the results clearly indicate that the high-thrust short arcs of thrusting as well as limiting case of impulses are predominantly applied near peripasis/apoapsis and/or the ascending or descending nodes, as might be anticipated for orbit raising/lowering and changes in inclination [\citen{vinh1988optimal}]. We also note that, multiple examples show that the time, direction and magnitude of impulses are accurately approximated by the high thrust limits of the corresponding minimum-fuel switching surfaces. So in that sense, we have established a means to answer, even if not uniquely, the original minimum-impulsive velocity maneuver question raised by Edelbaum in 1967. 

The impulsive idealization has always served as a first step in assessing mission feasibility and to establish candidate approximate trajectories. The methodology of this paper can be used to establish multiple minimum impulsive velocity extremals. More importantly, we provide the means to establish an associated family of minimum-fuel extremals, considering a specific spacecraft and a family of engine designs. For a maximum thrust specification, these extremals approach the impulsive limit.
\vspace{3mm}
\begin{tcolorbox}[colback=green!5,colframe=green!40!black,title= Minimum-fuel local extreme conjecture]
Based on the results, we have the following conjecture: for \textbf{rendezvous-type, fixed-time, minimum-fuel trajectories in a Newtonian gravitational field}, i.e., $1/r^2$, there is at most one local extremal for each integer specification of the number of en-route revolutions $N_{rev} \in \{0,1,2,\cdots\}$ made during the transfer orbit.
\end{tcolorbox}
\subsection{Choice of Coordinates/Elements And Its Impact On Switching Surfaces} \label{sec:choiceofset}
It is important to emphasize that switching surfaces are a concatenation (or an ensemble) of switching functions as the parameter of interest is swept over the defined range. Thus, switching functions are the building blocks that form the topology of the switching surfaces. The underlying optimal state and co-state trajectories that minimize fuel consumption constitute an extremal field map (actually a hyper-surface, in this case, of dimension 14). 

On the other hand, these switching functions are associated with solutions that are guaranteed to be extremals since they are systematically generated to satisfy the Pontryagin's necessary conditions. The optimal thrust is a  independent identity and described as a physical vector function associated with each optimal trajectory in the extremal field map. Irrespective of the chosen set of coordinates or elements used for modeling the motion of spacecraft, the thrusting regions of space are invariant under coordinate transformation. This means that switching surfaces are independent of the choice of coordinates. Furthermore, each extremal trajectory for any choice of coordinates can be mapped to another coordinate choice through a coordinate transformation applied at each time instant.  
\begin{figure}[htbp!]
\centering
\includegraphics[width=3.5in]{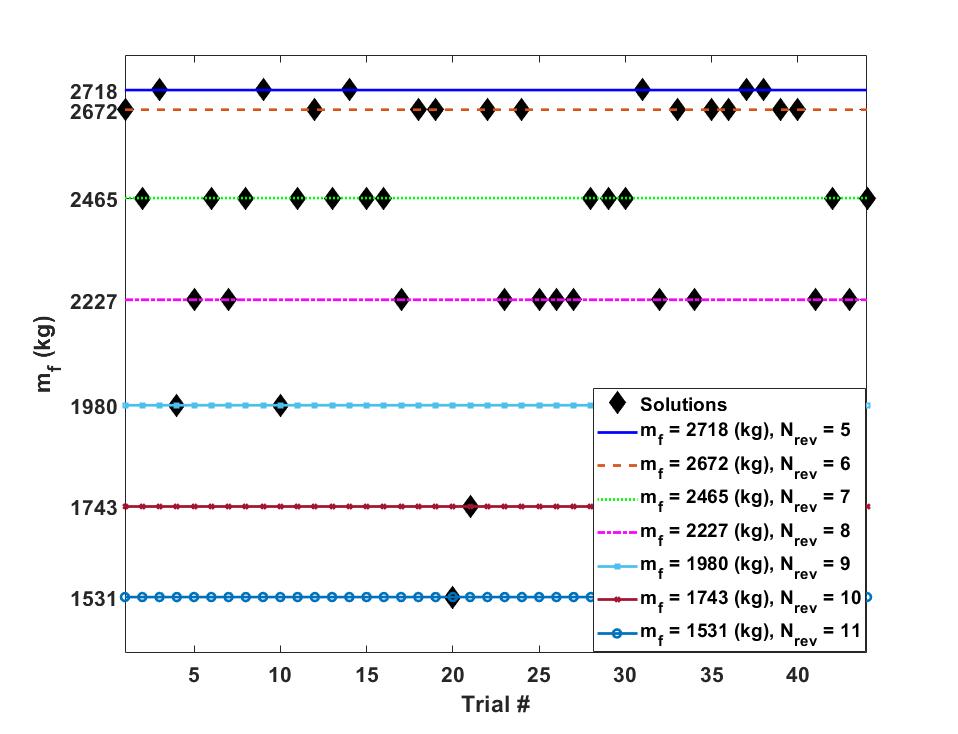}
\caption{Earth-to-Dionysus seven local solutions using Cartesian coordinates.}
\label{fig:ED_sol_levels}
\end{figure}

On the other hand, as explored in [\citen{taheri2016enhanced, junkins2018exploration}], it is known that the choice of coordinates and/or elements is important when it comes to solving TPBVPs. Set of Cartesian coordinates is shown to produce sub-optimal solutions. For instance, for the same boundary conditions and with a thrust magnitude of $T = 0.32$ N, the problem of minimum-fuel trajectories from Earth to asteroid Dionysus is studied in [\citen{taheri2016enhanced}] when the equations of motion are written in terms of Cartesian coordinates and MEEs. Fifty TPBVPs are solved where the solution procedure involves with a random initialization of the unknown initial elements of the co-state vector. A similar continuation procedure is used for solving each TPBVP; only 43 of the trials converged. The converged trials reveal the existence of seven levels of extremal solutions that satisfy the optimality conditions, where six of these levels correspond to sub-optimal solutions. Figure \ref{fig:ED_sol_levels} shows these solution levels (values of the final mass).

\begin{figure}[htbp!]
\begin{multicols}{2}
\centering
\includegraphics[width = 0.5\textwidth]{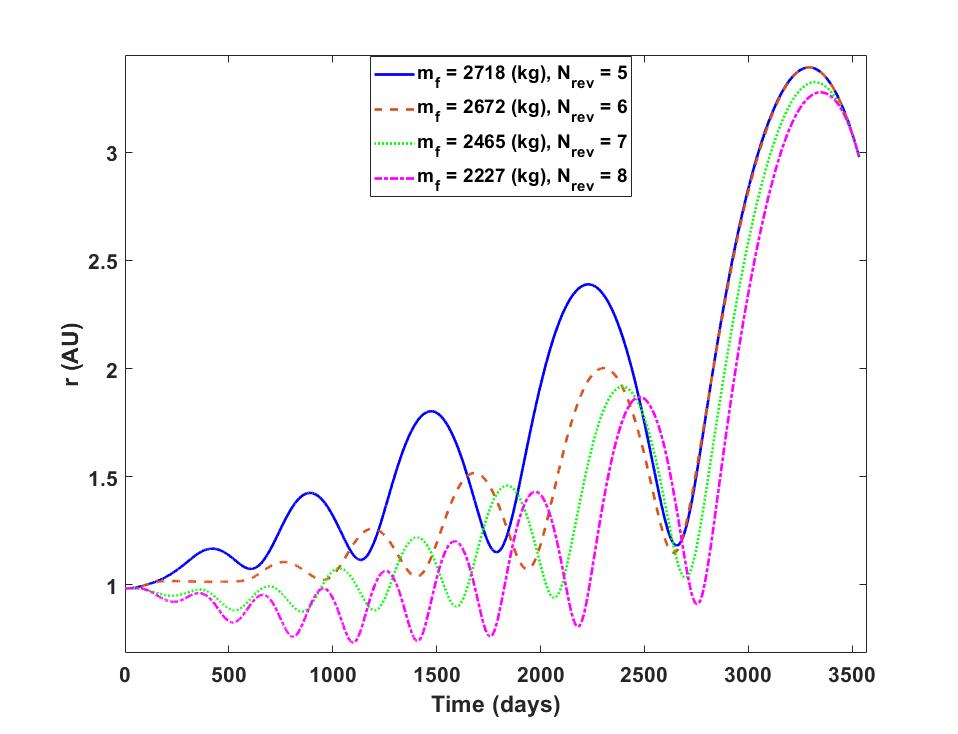}
\caption{Profile of the radius versus time of flight for the local solutions \#1 to \#4.}
\label{fig:ED_levels1to4}
\centering
\includegraphics[width = 0.5\textwidth]{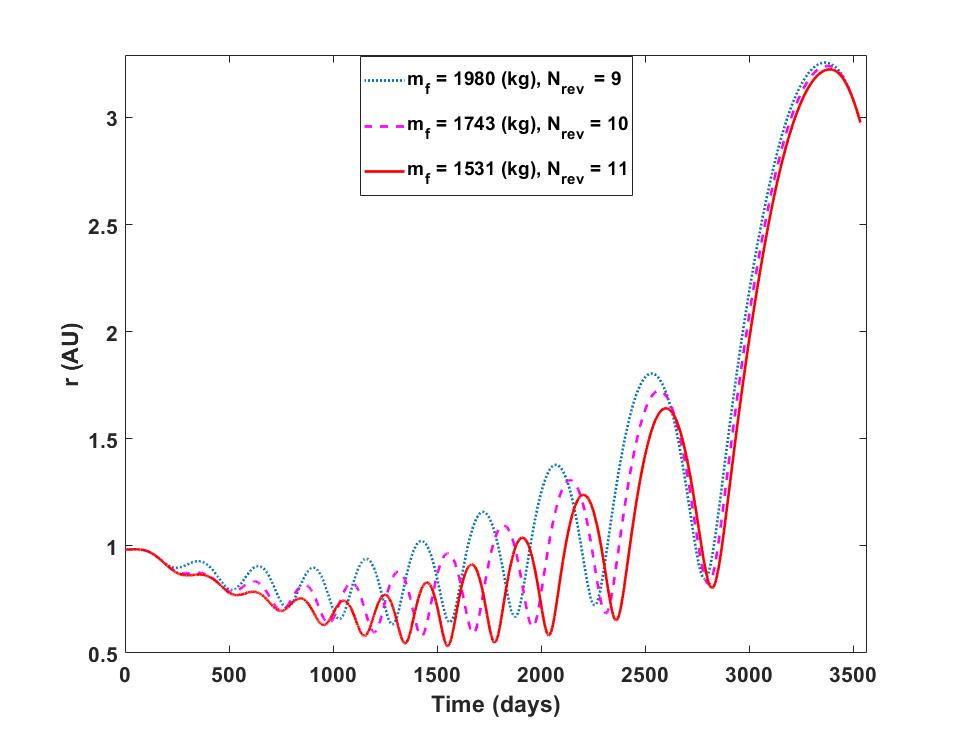}
\caption{Profile of the radius versus time of flight for sub-optimal levels \# 5 to \# 7.}
\label{fig:ED_levels5to7}
\end{multicols}
\end{figure}

There are actually two main reasons for obtaining different local sub-optimal solutions. The first one is related to how TPBVPs are typically formulated when Cartesian coordinates are used. Unlike the other choices of coordinates or elements, the number of revolutions, $N_{\text{rev}}$ does not appear explicitly as one of the parameters in the formulation of the TPBVP when Cartesian coordinates are used. The terminal constraints are on the position and velocity coordinates, and possibly the final value of the co-state associated with mass, $\lambda_m(t_f)$. This can be interpreted as a positive aspect of using Cartesian coordinates as we usually do not know the correct value of the number of revolutions. At the same time, this freedom has a negative consequence, which is the main cause of getting sub-optimal solutions (the number of revolutions may be known from prior studies). This freedom can be interpreted to be positive and negative simultaneously. For any systematic study, the mentioned freedom has to be harnessed. 
\begin{figure}[htbp!]
\begin{multicols}{2}
\centering
\includegraphics[width = 0.5\textwidth]{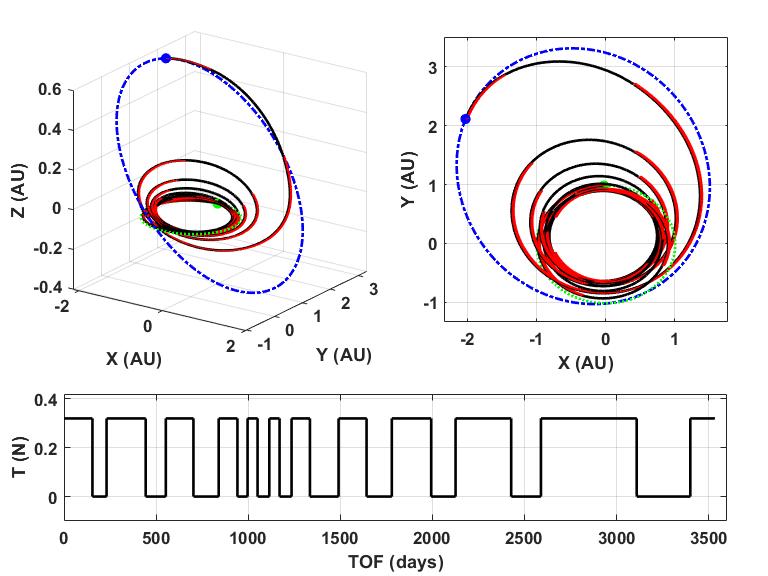}
\caption{Earth-to-Dionysus sub-optimal solution with $m_f = 1980$ kg.}
\label{fig:ED_Sub1}
\centering
\includegraphics[width = 0.5\textwidth]{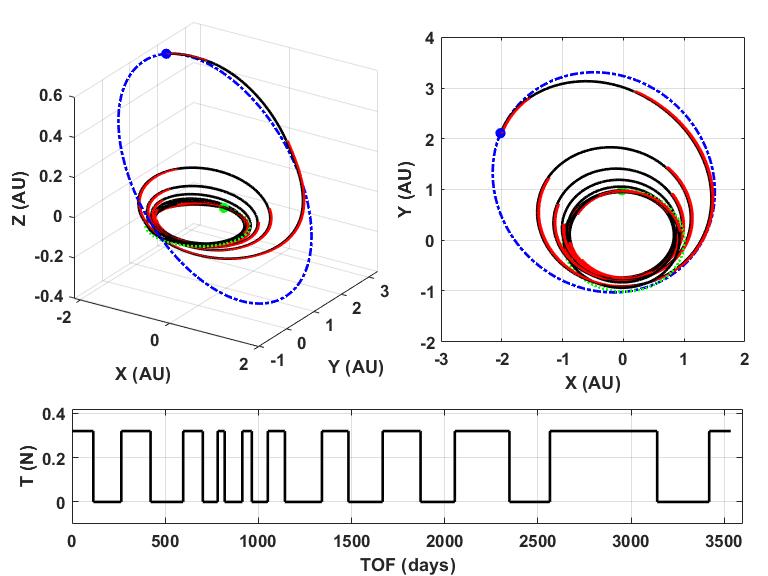}
\caption{Earth-to-Dionysus sub-optimal solution with $m_f = 2227$ kg.}
\label{fig:ED_Sub2}
\end{multicols}
\end{figure}
For minimum-fuel problems with a large time of flight, the optimal trajectories have a tendency to get close to the central body (Sun in our problem) to utilize its gravitational potential to early on gain kinetic energy. Efficient increase of the energy of orbits is achieved at perihelion passages of elliptical orbits where the velocity takes its maximum value. In other words, for some optimal trajectories, the optimal solution initially dives toward the central body. This is a characteristic of multi-revolution trajectories and impacts the solution when the number of revolutions are significantly large. But, for problems with relatively small number of revolutions, this is frequently encountered. Figures \ref{fig:ED_levels1to4} and \ref{fig:ED_levels5to7} show the profile of the magnitude of the radius vector, $r = \| \textbf{r}\|$ versus time of flight for representative solutions of all of the seven levels of solutions. Clearly, those solutions with a higher number of revolutions (compared to the optimal solution with the lowest number of revolutions) require more fuel. The optimal solution makes only five revolutions around the Sun, whereas the worst sub-optimal solution makes eleven revolutions around the Sun. Figure \ref{fig:ED_levels5to7} shows clearly what we stated earlier where the majority of the time of flight is spent at lower radius, and closer to the Sun. Our intent is not to perform an extensive analysis on the number of levels of sub-optimal solutions, but these results indicate that formulation of minimum-fuel optimal control problems using Cartesian coordinates will result, most of the time, in sub-optimal solutions. 

\begin{figure}[htbp!]
\begin{multicols}{2}
\centering
\includegraphics[width = 0.5\textwidth]{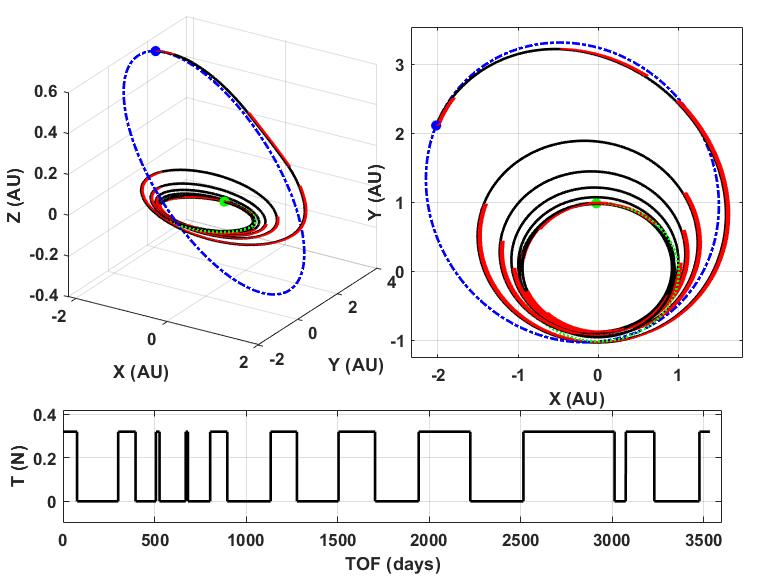}
\caption{Earth-to-Dionysus sub-optimal solution with $m_f = 2465$ kg.}
\label{fig:ED_Sub3}
\centering
\includegraphics[width = 0.5\textwidth]{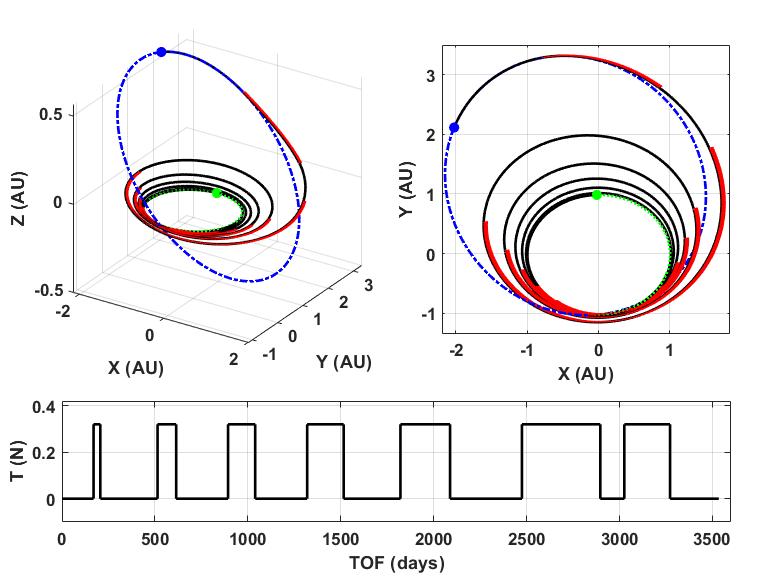}
\caption{Earth-to-Dionysus sub-optimal solution with $m_f = 2672$ kg.}
\label{fig:ED_Sub4}
\end{multicols}
\end{figure}

The existence of these sub-optimal solutions (associated with $N_{\text{rev}}$ en route revolutions) is a direct consequence of large time of flight and thrust acceleration magnitude and the additional freedom to make as many revolutions as possible, within the given time. Figures \ref{fig:ED_Sub1}-\ref{fig:ED_Sub4} show trajectories and thrust profiles of the four sub-optimal solutions. %When this additional freedom of the number of revolutions is removed from the problem and is fixed, like what we have implemented when using MEEs, any \text{converged} solution is the ``optimal'' solution.  

\subsection{GTO To A Halo Orbit Around L1}\label{sec:HaloProblem}
This example illustrates how to design an optimal low-thrust transfer from a Geostationary Transfer Orbit (GTO) to an L1 halo orbit in the restricted three-body problem associated with the Earth-Moon system. The equations of motion are given in [\citen{taheri2018genericconf}]. The departure of the spacecraft is from the perigee of an elliptical orbit with perigee $\times$ apogee altitude ratio of $400 \times 35,864$ km (GTO). The spacecraft with an initial mass of $m_0 = 1500$ kg uses an electric engine with a specific impulse of $I_{sp} = 3000$ seconds.
\begin{figure}[h!]
\centering
\includegraphics[width = 0.5\textwidth]{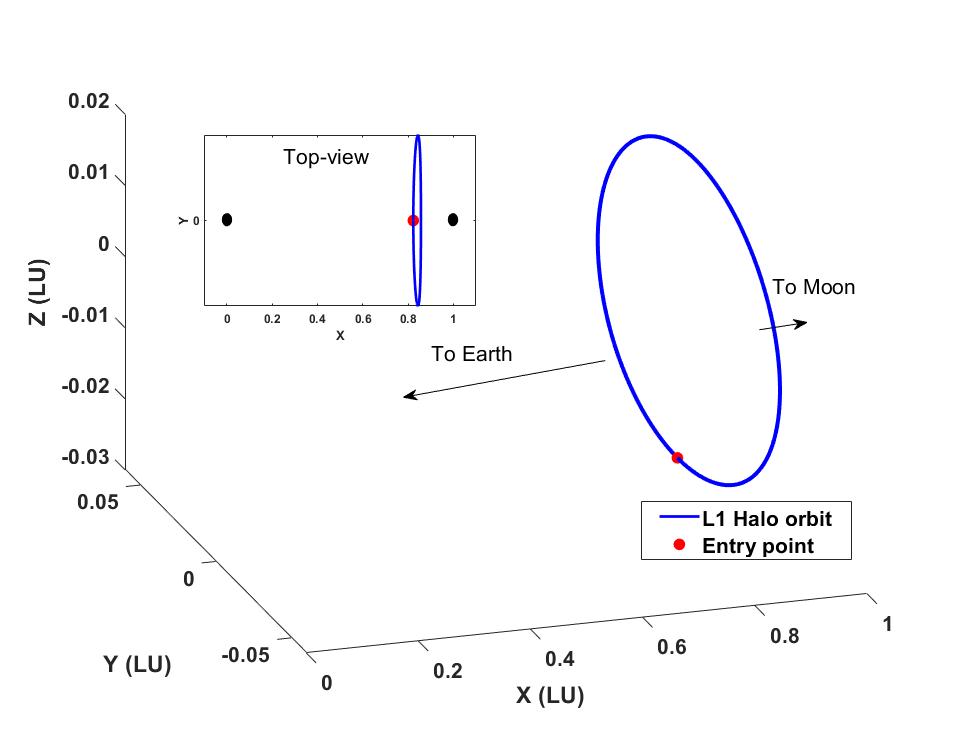}
\caption{Halo orbit and point of entrance to the orbit.}
\label{fig:Haloorbit}
\end{figure}
The target orbit is a halo orbit with an out-of-plane amplitude of 8000 km [\citen{martin2010optimal}]. Figure \ref{fig:Haloorbit} shows the halo orbit and the point on orbit at which the spacecraft will enter the orbit. The parameters used in numerical simulations are : $g_0 = 9.80665$ (m/s$^2$), Length unit (LU) = 3.84405000 $\times 10^5$ (km), Time unit (TU) = $3.75676967 \times 10^5$ (s), Velocity unit (VU)  = $1.02323281$ (km/s), $\mu = 1.21506683 \times 10^{-2}$ (LU$^3$/TU$^2$). LU in actually the reference semi-major axis of the Moon's orbit relative to the Earth (the distance between Earth and Moon).
\begin{figure}[htbp!]
	\centering
    \includegraphics[width=4.0in]{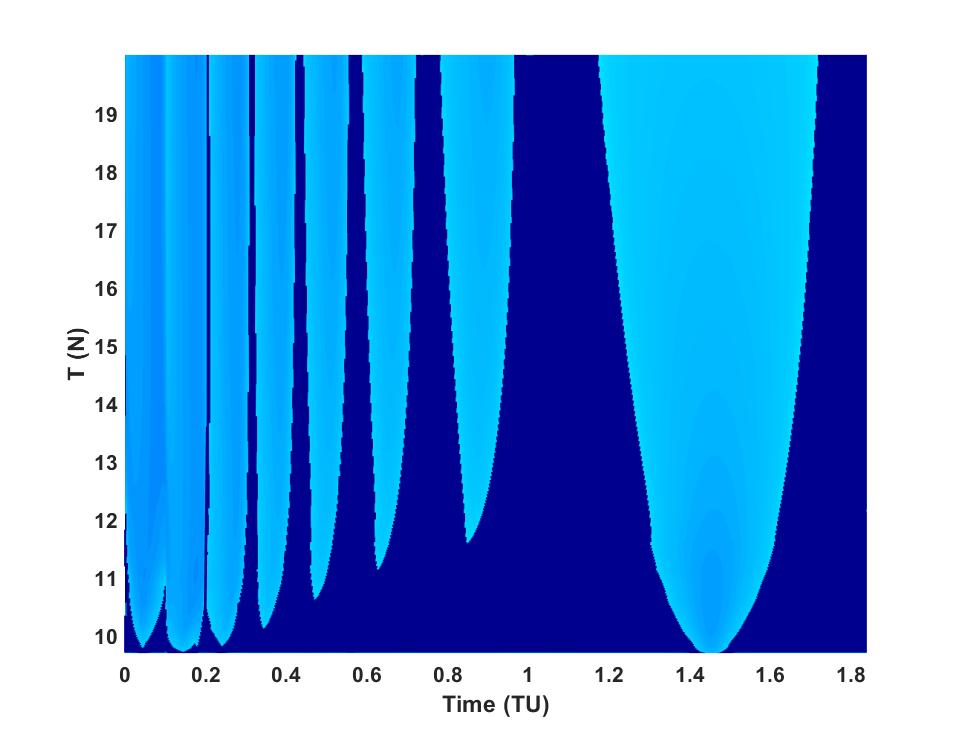}
	\caption{Switching surface for sub-optimal trajectories for GTO-to-Halo problem with $N_{\text{rev}} = 6$.}
	\label{fig:HaloSS}
\end{figure}

This problem is more difficult compared to the previous cases where the dynamics were simpler. For the three-body dynamical models, convergence is found to be difficult to achieve using a finite difference approach for calculating sensitivities. Therefore, the STM method (with hyperbolic tangent smoothing) is used for the sensitivities according to the procedure described in [\citen{taheri2016enhanced}]. When integrating the state and co-state equations, we make use of hyperbolic tangent smoothing to enhance the convergence through a continuation procedure.

\begin{figure}[htbp!]
\begin{multicols}{2}
\centering
\includegraphics[width = 0.5\textwidth]{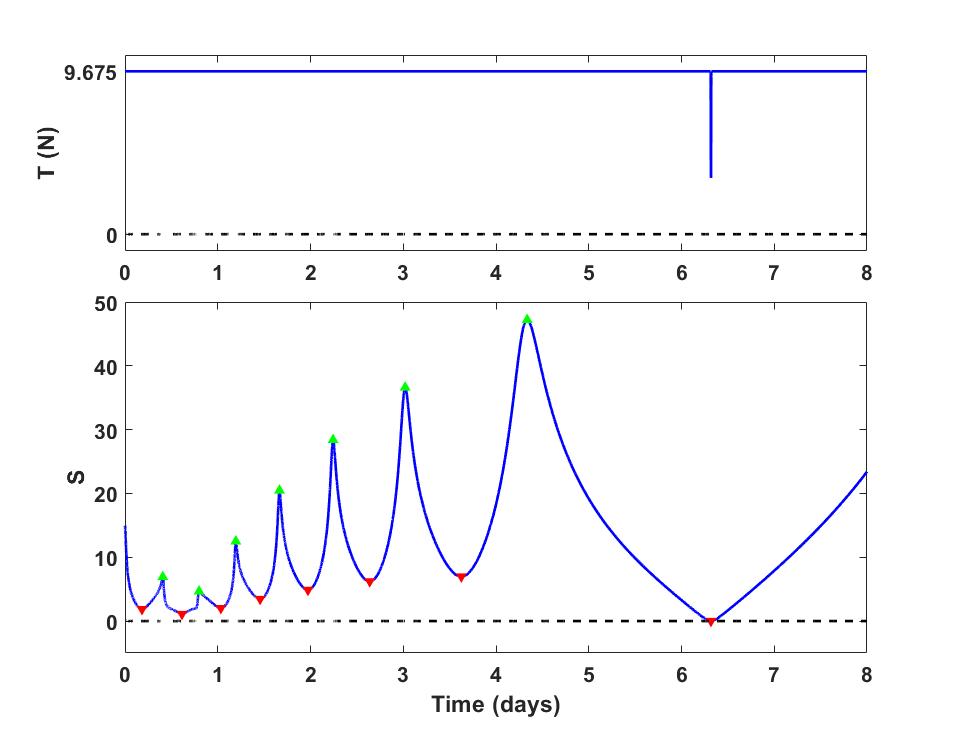}
\caption{GTO-to-Halo switching function with $T \approx  9.675$ N and $N_{\text{rev}} = 6$.}
\label{fig:HaloCT1_SF}
\centering
\includegraphics[width = 0.5\textwidth]{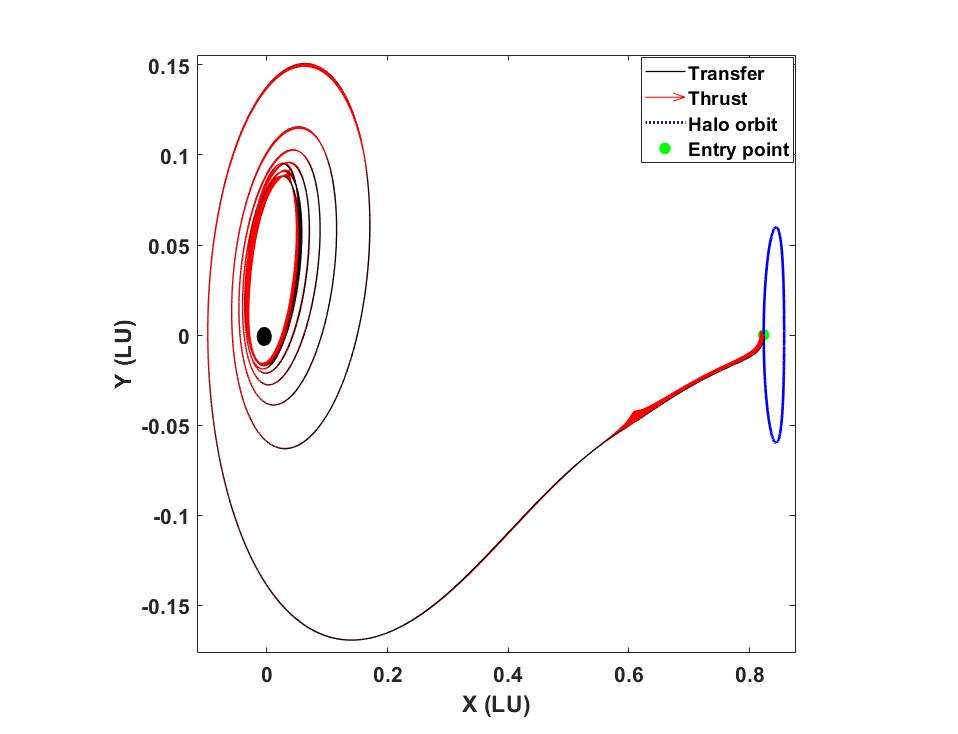}
\caption{GTO-to-Halo trajectory with $T \approx  9.675$ N and $N_{\text{rev}} = 6$.}
\label{fig:HaloCT1_Traj}
\end{multicols}
\end{figure}

Finding the minimum-thrust solution is the first step for constructing the switching surface. However, the minimum-thrust solution is not unique, which is an important factor. For instance, Figures \ref{fig:HaloCT1_SF} and \ref{fig:HaloCT1_Traj} show the switching function and thrust profile and trajectory for a minimum-thrust solution with $T_{\text{min}} \approx 9.675$ N with $N_{\text{rev}} = 6$. Clearly, the switching function is positive along the whole trajectory. However, this is not the fundamental minimum-thrust solution since the number of revolutions are greater than the global solution. If this local extremal solution is considered as the base solution and the thrust magnitude is swept, Figure \ref{fig:HaloSS} shows the resulting switching surface for $T \in [9.675, 20]$ N with $N_{\text{rev}} = 6$. Note that when Cartesian elements are used for formulating the TPBVPs, the number of revolutions does not appear explicitly. There appears to be seven main thrust ridges that remain at high thrust levels. There is a late-departure boundary, whereas the last thrust ridge is attached to the terminal time, which indicates that there is no early-arrival boundary. 

\begin{figure}[htbp!]
\begin{multicols}{2}
\centering
\includegraphics[width = 0.5\textwidth]{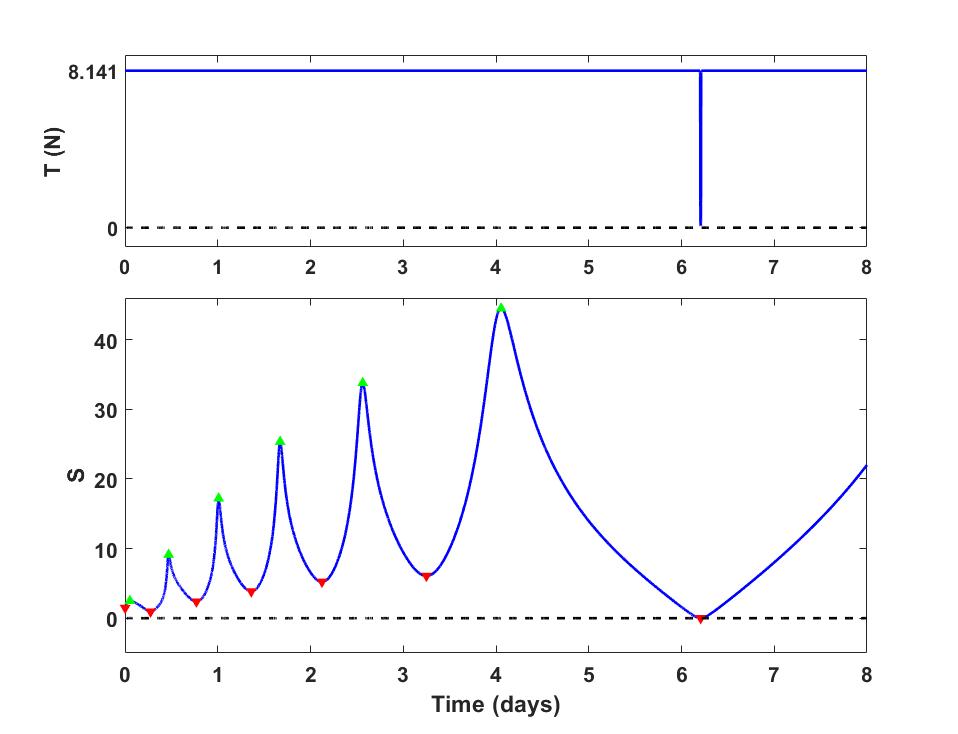}
\caption{GTO-to-Halo switching function with \\$T \approx  8.141$ N.}
\label{fig:GTO2HaloCT1_SF_Op}
\centering
\includegraphics[width = 0.5\textwidth]{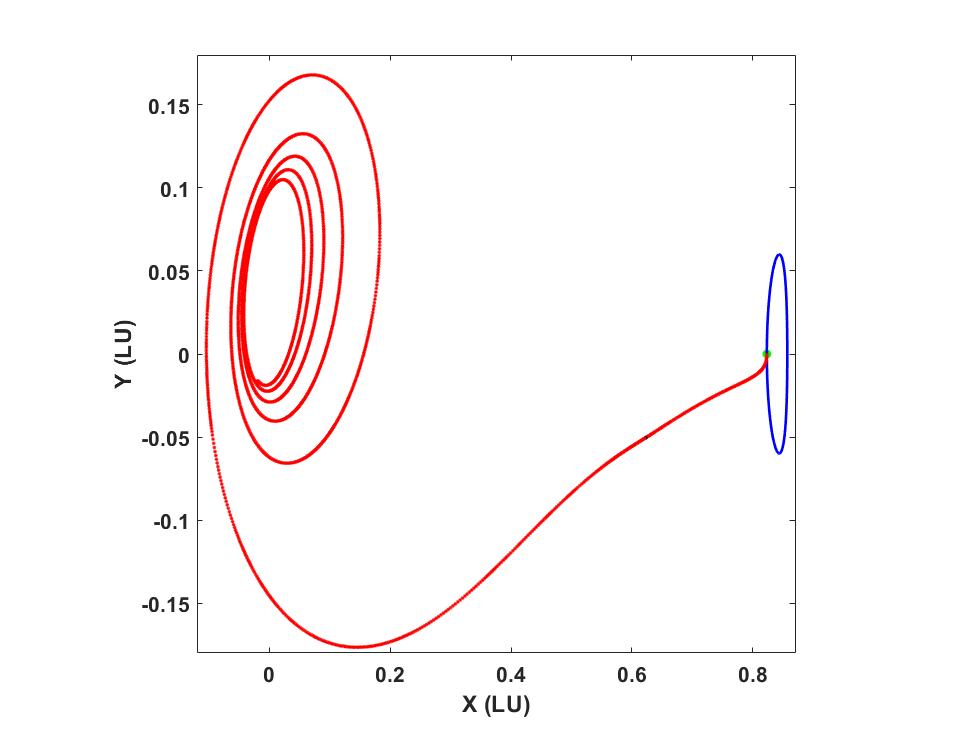}
\caption{GTO-to-Halo trajectory with $T \approx  8.141$ N.}
\label{fig:GTO2HaloCT1_Traj_Op}
\end{multicols}
\end{figure}

However, the global minimum-thrust solution of this problem has been found to correspond to $T_{\text{min}} \approx 8.141$ N with $N_{\text{rev}} = 5$. Figures \ref{fig:GTO2HaloCT1_SF_Op} and \ref{fig:GTO2HaloCT1_Traj_Op} show the switching function and thrust profile and trajectory for the fundamental minimum-thrust solution with $T_{\text{min}} \approx 8.141$ N. It is interesting to note that in Figure \ref{fig:HaloCT1_SF}, the profile of the switching function has a different behavior during the time interval $t \in [0.5, 1]$ compared to the smooth one in Figure \ref{fig:GTO2HaloCT1_SF_Op}. 
%\textcolor{red}{feature may be used as a criterion for figuring out that the ``optimal'' minimum-thrust case has been found or not}. 

\begin{figure}[htbp!]
\begin{multicols}{2}
\centering
\includegraphics[width = 0.5\textwidth]{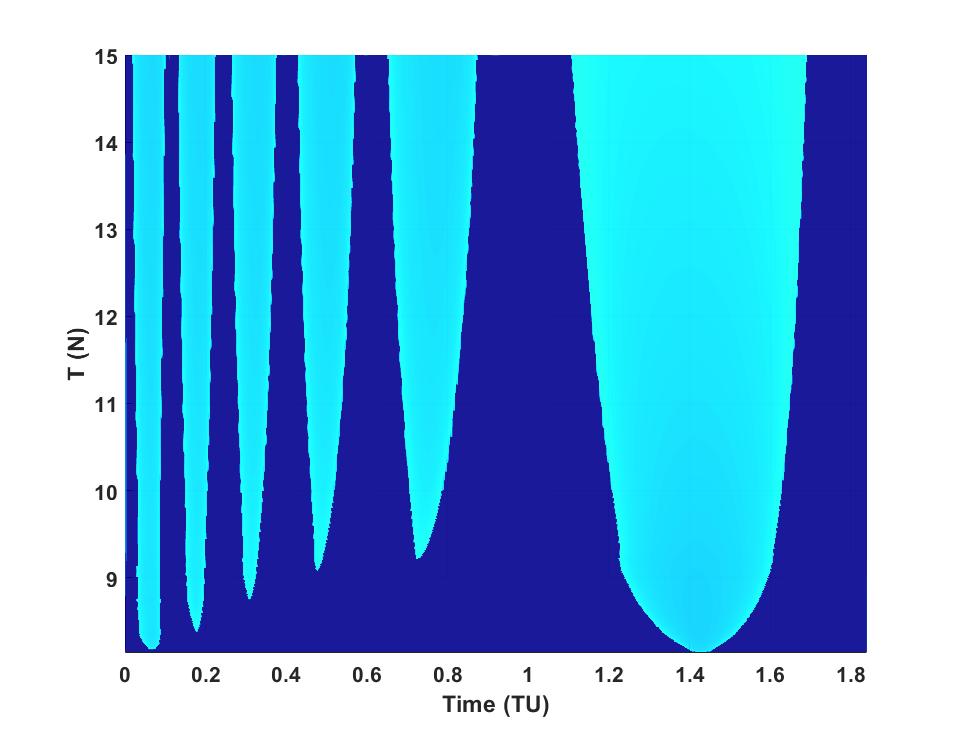}
\caption{Switching surface and thrust profile for optimal trajectories, $N_{\text{rev}} = 5$ for GTO-to-Halo problem with $T \in [8.141, 15]$ N.}
\label{fig:HaloSS_Op}
\centering
\includegraphics[width = 0.5\textwidth]{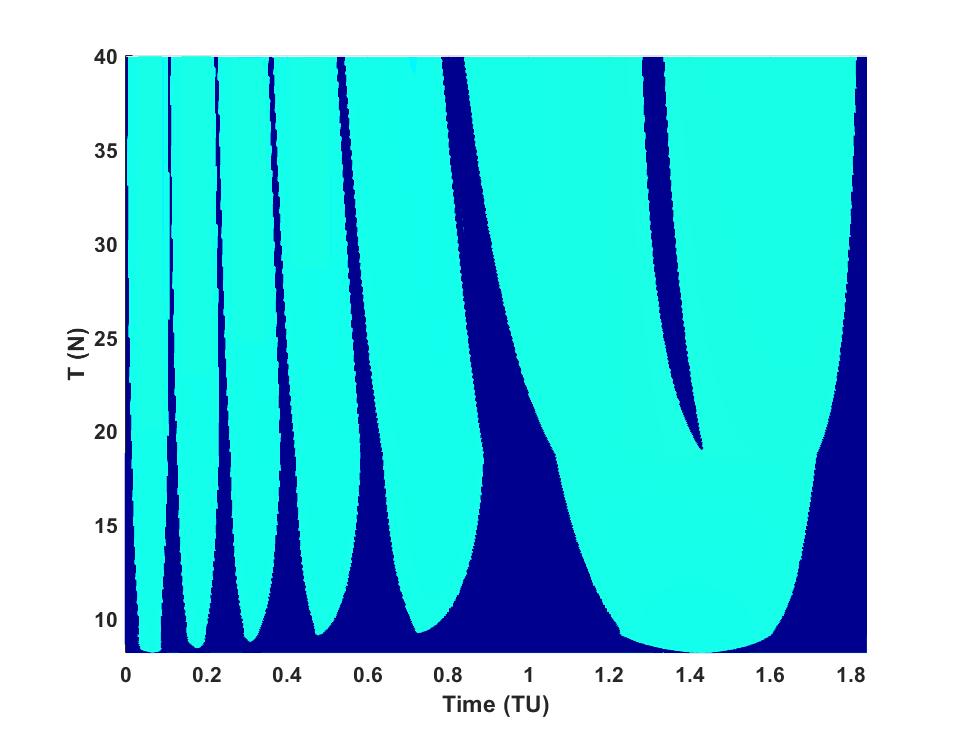}
\caption{Switching surface for optimal trajectories, $N_{\text{rev}} = 5$ for GTO-to-Halo problem $T \in [8.141, 40]$ N.}
\label{fig:HaloSS_Op2}
\end{multicols}
\end{figure}

Figure \ref{fig:HaloSS_Op} shows the resulting optimal switching surface when $T \in [8.141, 15]$ N is considered. There appears to be only seven main thrust ridges that remain at ``high'' thrust levels. The boundary conditions result in a switching surface without early-arrival or later-departure boundaries. Note that the thrust is not necessarily high enough as the width of the thrust ridges are still large.

Figure \ref{fig:HaloSS_Op2} shows the optimal switching surface for $T \in [8.141, 40]$ N when $\rho_{\text{min}} = 9.15 \times 10^{-6}$. Note that in practice, a relatively small value is chosen for the smoothing parameter in order to facilitate the generation of the switching surface. For particular regions of interest, it is possible to use a smaller value for $\rho_{\text{min}}$.

\begin{figure}[htbp!]
\begin{multicols}{2}
\centering
\includegraphics[width = 0.5\textwidth]{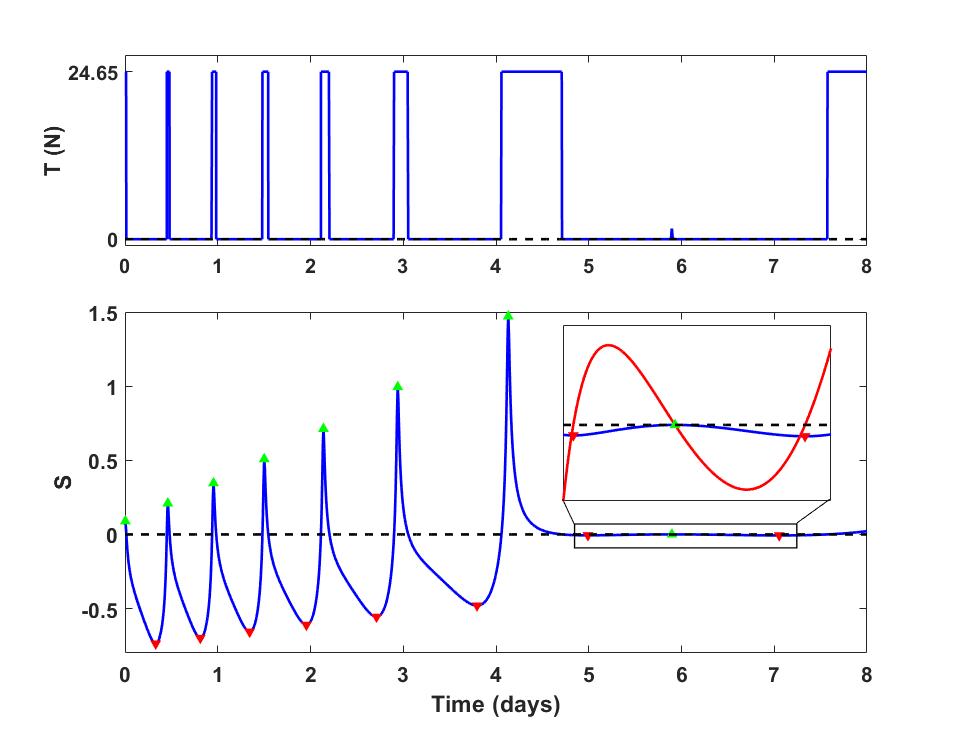}
\caption{GTO-to-Halo switching function and thrust profile with $T \approx 24.65$ N.}
\label{fig:HaloBifur_SF}
\centering
\includegraphics[width = 0.5\textwidth]{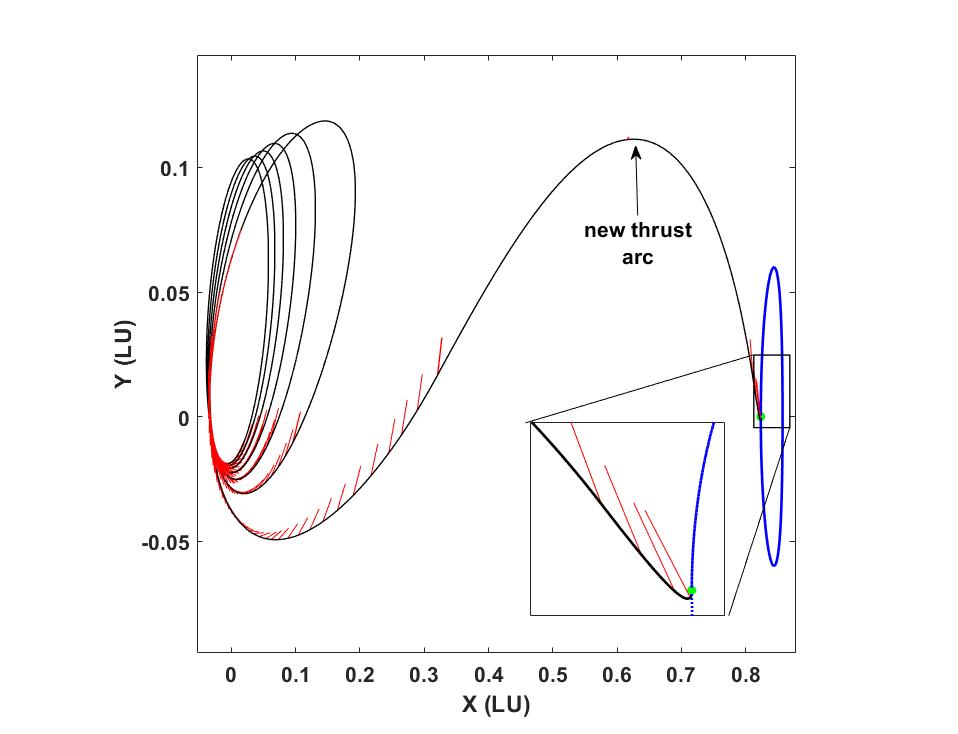}
\caption{GTO-to-Halo trajectory with $T \approx 24.65$ N.}
\label{fig:HaloBifur_Traj}
\end{multicols}
\end{figure}
There is an interesting bifurcation phenomenon at $T \approx 19$ N, see Figure \ref{fig:HaloSS_Op2}. The nature of this bifurcation point is completely different from the previously shown bifurcations in the two-body dynamics (where a corner point was the beginning point of the bifurcation and the beginning of the thrust arc was attached to a thrust ridge, see Fig. \ref{fig:SSEM_lower}). In addition, in the two-body dynamics $\ddot{S} >0$ in the vicinity of the bifurcation point, whereas here the sign of $\ddot{S}$ changes with $S=\dot{S}=\ddot{S}=0$ at the bifurcation; hence, this bifurcation appears in the middle of a coast canyon. Clearly, the creation of this thrust ridge has also visibly correlated to sharp features on the two closest thrust ridges to its right and left where two corner points are noticeable and their time duration (width) shrinking rate with increasing thrust is increased. 

\begin{figure}[htbp!]
\centering
\includegraphics[width = 0.5\textwidth]{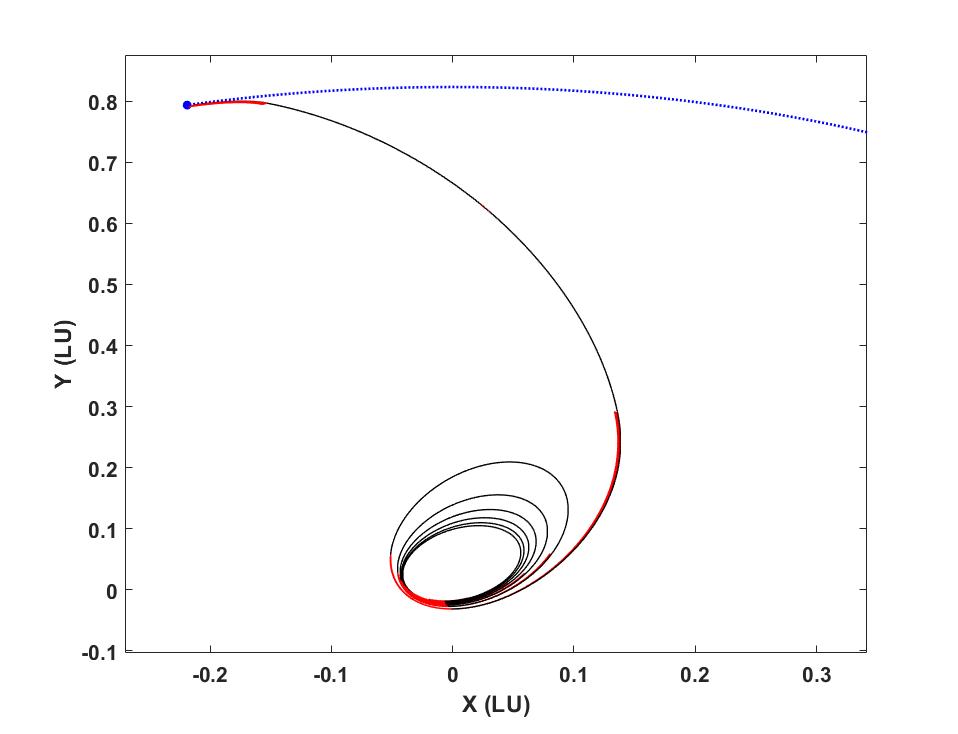}
\caption{Inertial trajectory for GTO-to-Halo problem with $T \approx 24.65$ N.}
\label{fig:EHalo_Bifur_InertialTraj}
\end{figure}

Figures \ref{fig:HaloBifur_SF} and \ref{fig:HaloBifur_Traj} show the switching function and the trajectory for this bifurcation point. In order to be able to isolate shorter thrust arcs more accurately, we set $\rho_{\text{min}} = 1.0 \times 10^{-6}$, which shows a small increase in the thrust magnitude at which this bifurcation occurs. The blown-up part shows the time when the switching function (blue line) nearly touches the $S = 0$ line, whereas the red line shows the derivative of the switching function with respect to time, $\dot{S}$. Figure \ref{fig:HaloBifur_Traj} shows the trajectory where the point of the creation of the new thrust arc is shown. The final phase of the flight consists of a thrust arc where the trajectory undergoes a sudden turn an reaches the entry point to the L1 orbit with a final mass of $m_f = 1394.85$ kg. The trajectory looks different when plotted in a non-rotating frame as is shown in Figure \ref{fig:EHalo_Bifur_InertialTraj}, where the blue line denotes the path made by the entry point along the time of flight. It also shows a 3D view of the trajectory in which the out-of-plane motion is noticeable. 

\section{Computational Effort} \label{sec:comp}
%\vspace{-3mm}
In this section, computational effort is discussed. The proposed methodology relies on solving TPBVPs to generate switching functions, which, as thrust is swept, gives rise to neighboring extremals that underly the optimal switching surface. The first step has to do with the determination of the minimum-thrust solution for each $N_{\text{rev}}$. The second step, deals with generating the corresponding switching surfaces. A third step, if one desires to determine the idealized impulsive trajectories, is to perform $N$-impulse optimization. 

Solving TPBVPs efficiently and reliably in the presence of high nonlinearity, and high dimensionality is a challenging task [\citen{betts1998survey,woolley2019cargo}] and there are several approaches to enhance convergence when discontinuous inputs are forcing the differential equations [\citen{taheri2016enhanced}]. Our experience is that for a great variety of optimal trajectory cases, the hyperbolic tangent smoothing is quite helpful and alleviates a number of difficulties associated with solving the TPBVPs associated with switching-type controls [\citen{woollands2019efficient}]. On the other hand, if our goal is to determine a set of impulsive solutions, in order to generate the results similar to those reported in Table \ref{tab:impulsiveED_DiffNrev}, it is not (usually) required to sweep over the thrust magnitudes, in a dense incremental fashion. Instead, it is possible to choose a relatively large thrust level and solve the corresponding minimum-fuel TPBVP directly. Even though, we use the hyperbolic tangent smoothing, only the final trajectory needs to approach the very sharp control switches. As mentioned earlier, the impulsive solution reveals itself at high-thrust levels encountered in the family of optimal trajectory problems under consideration.

This strategy has been adopted to generate the results in Table \ref{tab:impulsiveED_DiffNrev}. For the Earth-to-Dionysus problem, the computation wall-time to obtain an impulsive solution is usually only 30 to 60 seconds. This time takes into account the time it takes to solve minimum-thrust problem according to the methodology outlined in in this paper while relying on hyperbolic tangent smoothing, converging with a large value of thrust that solves the corresponding minimum-fuel TPBVP, including the optimal switching function, and the time it takes to solve the $N$-impulse optimization problem. The $N$-impulse optimization is the least computationally demanding part of the overall scheme. Of course, harder problems with highly nonlinear trajectories and more impulses require more computation time. The GTO-to-GEO orbit transfer problem proved to be the most computationally demanding case when two-body dynamics was considered. Generating the GTO-to-GEO impulsive solutions took several minutes. %Generation of the switching surfaces associated with the trajectories in the restricted three-body dynamics was also computationally demanding. Generation of the switching surfaces associated with the trajectories in the restricted three-body dynamics was the most demanding task.

All computations were performed on a Core i7 desktop machine with two 3.4-GHz processors and 16 GB of RAM and the propagation of the differential equations was achieved using a compiled version of MATLAB's \textit{ode45} function to speed up the numerical simulations. The absolute and relative integration convergence tolerances were set to $1.0 \times 10^{-10}$ in all simulations. For MATLAB's \textit{fmincon} optimizer, the function and constraint tolerances were set to $TolFun = 10 \times 10^{-6}$ and $TolCon = 10 \times 10^{-7}$, respectively.
\section{Conclusion} \label{sec:conclusion}
%\vspace{-3mm}
We used indirect methods of optimal control theory along with Lawden's primer vector theory to introduce illuminating switching surfaces for large and low-thrust minimum-fuel trajectories. These surfaces establish the connection between impulsive and continuous-thrust trajectories. The impulsive limits tell us ``how many impulses'' are associated with each application of $N_{\text{rev}}$, the number of revolutions. 

We have shown that there exists a \textit{fundamental} minimum-thrust solution that plays a pivotal role in determining fundamental switching surface, and an associated $N^*_{\text{rev}}$, which in turn reveals the ``optimal'' number of thrust arcs for any finite specification of maximum thrust. The fundamental switching surface is generated by a homotopic sweep of the maximum thrust away from the the minimum-thrust extremal among all minimum-thrust solutions, considering all feasible multi-revolution solutions. This fundamental minimum-thrust solution is used to construct the switching surfaces that reveal the number and approximate time, direction and magnitudes of the impulses associated with the optimal impulsive solution as a by-product. A number of optimal orbit transfers with different number of revolutions for the fundamental extremal are studied and the corresponding impulsive solution with as many as eleven impulses are found, which demonstrate the capability of the proposed construct. According to these results and our experience, the optimal number of thrusting arcs at high thrust levels depend strongly on the size and shape of the initial and final orbits, time of flight, angular phasing and the orbits’ relative orientation. 

While the switching surface can be generated for very high thrust values to approach with arbitrary precision the impulsive limit, we find that after the elapsed thrust-on time of all short thrust arcs is less than some tolerance (e.g., $<0.001 \times$ (time of flight)), then the impulsive approximation can be invoked with a good accuracy. A direct method can then be initiated accurately to converge to the impulsive solution and satisfy the force model and all boundary conditions to a high precision. It is vital to note that when mass variation is accounted for, considering any real propulsion system, that minimizing $\Delta v$ does not minimize fuel consumption (or maximize payload mass for a given engine). We showed two examples where minimizing $\Delta v$ led to multiple distinct local extremals, with identical minimum $\Delta v$ values; the local extremals $\Delta v$ matched to 7 digits! Only one of these extremals belonged to the fundamental switching surface that minimizes fuel consumption. So, we conclude that Edelbaum's question does not have a unique answer, in general. However, invoking the minimum propellant consumption for a specific engine model gives a unique extremal.

We showed that other secondary criteria, such as time of flight (for cases with early arrival) can also be invoked, for the same total $\Delta v$ and modest fuel increases. The results clearly indicate the intuitively reasonable truth that the optimal short thrusting arcs as well as limiting case of impulses are predominantly applied near peripasis/apoapsis and/or the ascending or descending nodes, as might be anticipated for orbit raising/lowering and changes in inclination.  

Finally, we note that the switching surfaces provide a unified means for optimizing low thrust, high thrust, and impulsive maneuvers, and enable important global insights in mission design. In principle, the methodology proposed in this paper can be pursued for many engineering systems if/when a systematic study of the resulting switching surfaces are performed by sweeping various important parameters of interest.  

%Numerical results indicate that location and number of extremal solution bifurcations, which occur occasionally for critical thrust values and lead to creation of thrusts ridges, alter the number of optimal thrust arcs (that exist at high thrust), which ultimately dictate the number of impulses. Our hypothesis is that, the number of impulses depends on the optimal number of revolutions associated with the fundamental minimum-time solution, and the number of bifurcations; these however, are all properties of the extremal field map switching surface associated with the particular boundary conditions and force model. 

\section*{Conflict of Interest Statement}
\vspace{-1mm}
On behalf of all authors, the corresponding author states that there is no conflict of interest. 
\vspace{-1mm}
\section*{Acknowledgments}
We are pleased like to thank our sponsors: AFOSR (Stacie Williams), AFRL (Alok Das et al.) and ADS (Matt Wilkins), for their support and collaborations under various contracts and grants. We are indebted to many colleagues, we especially appreciate the interactions we have had on the ideas in this paper with Manoranjan Majji.
 
%\bibliographystyle{aiaa}

%\bibliographystyle{AAS_publication}   % Number the references.
%\bibliography{References.bib}   % Use references.bib to resolve the labels.
% \bibliography{References}

\end{document}